\documentclass[12pt]{article}

\usepackage{mathtools}
\usepackage[colorlinks=true,linkcolor=blue,citecolor=blue,urlcolor=blue]{hyperref}
\usepackage{amsmath}
\usepackage{amsthm}
\usepackage[noabbrev]{cleveref}
\usepackage{amsfonts}
\usepackage{amssymb}
\usepackage{url}

\usepackage{graphicx}
\usepackage{float}
\usepackage[makeroom]{cancel}
\usepackage{framed}
\usepackage{scalerel}
\usepackage{color}
\usepackage{yfonts}

\usepackage{caption}

\newcommand{\beaa}{\begin{eqnarray*}}
\newcommand{\eeaa}{\end{eqnarray*}}
\newcommand{\bea}{\begin{eqnarray}}
\newcommand{\eea}{\end{eqnarray}}

\newcommand{\Pro}{\noindent\textit{Proof.}\ \ }

\def\dd{\hskip 1pt{\rm{d}}}

\def\Var{\mathop{\rm Var\,}\nolimits}
\def\Cov{{\mathop{\rm Cov}}}

\def\bC{\mathbb{C}}
\def\bN{\mathbb{N}}
\def\bR{\mathbb{R}}
\newcommand\mi{\mathrm{i}\hskip 1pt}	

\numberwithin{equation}{section}
\theoremstyle{plain}
\newtheorem{theorem}{Theorem}[section]

\newtheorem{lemma}[theorem]{Lemma}
\newtheorem{proposition}[theorem]{Proposition}
\newtheorem{remark}[theorem]{Remark}
\newtheorem{conjecture}[theorem]{Conjecture}
\newtheorem{example}[theorem]{Example}
\newtheorem{assumption}[theorem]{Assumptions}

\theoremstyle{definition}

\begin{document}

\title{\bf \large
Integral Transform Methods in Goodness-of-Fit Testing, I: The Gamma Distributions
}

\author{{Elena Hadjicosta}{\hskip1pt}\thanks{Department of Statistics, Pennsylvania State University, University Park, PA 16802, U.S.A. E-mail address: \href{mailto:exh963@psu.edu}{exh963@psu.edu}.
\endgraf
\ $^\dag$Department of Statistics, Pennsylvania State University, University Park, PA 16802, U.S.A. E-mail address: \href{mailto:richards@stat.psu.edu}{richards@stat.psu.edu}.
\endgraf
\endgraf
}
\ \ and \ {Donald Richards}{\hskip1pt}$^\dag$
}

\maketitle

\date

\begin{center}
\vskip-10pt
{\em\large Dedicated to Professor Norbert Henze, on the occasion \\ of his 67th birthday}
\end{center}

\medskip

\begin{abstract}
We apply the method of Hankel transforms to develop goodness-of-fit tests for gamma distributions with given shape parameter and unknown rate parameter, thereby extending results of Baringhaus and Taherizadeh (2010) on the exponential distributions.  We derive the limiting null distribution of the test statistic as an integrated squared Gaussian process, obtain the corresponding covariance operator, and oscillation properties of its eigenfunctions.  We show that the eigenvalues of the operator satisfy an interlacing property, and we apply that property in approximating critical values of the test statistic in one of the two applications to data considered.  Further, we establish the consistency of the test.  In studying properties of the test statistic under a variety of contiguous alternatives, we obtain the asymptotic distribution of the test statistic for gamma alternatives with varying rate or shape parameters and for a class of contaminated gamma models.  We investigate the approximate Bahadur slope of the test statistic under local alternatives, and we establish the validity of the Wieand condition under which the approaches through the approximate Bahadur efficiency and the Pitman efficiency are in accord.
\end{abstract}

\tableofcontents

\normalsize

\parindent=25pt

\section{Introduction}
\label{introduction}
\setcounter{equation}{0}

The topic of goodness-of-fit testing for numerous classes of distributions has, in recent years, been the subject of intense activity.  As a consequence, there now exists a comprehensive body of results developed, by Henze and other authors, using test statistics based on integral transforms of Fourier, Laplace, Mellin, and related types, and making astute use of related differential equations and distributional characterizations.  The resulting test statistics have been shown to be superior in various ways to classical goodness-of-fit statistics, notably in comparisons of power, consistency, and in its behavior with respect to contiguous alternatives.  The depth and breadth of these results therefore provides motivation for further analysis and application of these methods.  

On reviewing the literature on goodness-of-fit tests, we were motivated by theoretical and applications considerations to consider the development, for certain multivariate exponential families \cite{BarLev} , tests that would be based on multidimensional integral transforms.  For such a program to be feasible, it was clear that the first step would be to determine whether such results can be obtained for the classical gamma distributions.  

In this paper, we apply Hankel transform methods to develop goodness-of-fit tests for gamma distributions with given shape parameter and unknown rate parameter, thereby extending the results of Baringhaus and Taherizadeh \cite{BT2010} and Taherizadeh \cite{ref27} for the exponential distributions.  We remark that we were fortuitous to have at hand the results of \cite{BT2010,ref27} as a constant guide in our investigations.

We remark that the gamma distributions with known shape parameters and unknown rate parameters arise as statistical models in numerous areas of research.  The most well-known example is the exponential distribution, corresponding to the shape parameter $\alpha = 1$ and which arises in many fields, e.g., queueing theory \cite{allen}.  For other values of $\alpha$, the gamma distribution with known $\alpha$ arises in studies of ion channel activation \cite{kass}; in the analysis of engineering equipment breakdowns \cite{czaplicki,sturgul}; and in the calculation of insurance premiums for sea ports and maritime transportation \cite{postan}.  

Turning to a description of the results in the paper, we begin by summarizing in Section \ref{propertieshankel} results from the theory of Bessel functions, confluent hypergeometric functions, and generalized Laguerre polynomials.  We refer wherever possible to the extant literature, and we provide proofs of results that we have not located elsewhere.  

In Section \ref{goodnessoffittests}, we apply the empirical Hankel transform to construct the test statistic, $T_n^2$.  We obtain the limiting null distribution of $T_n^2$ as an integral of the square of a centered Gaussian process $Z$.  We derive the covariance operator $\mathcal{S}$ corresponding to $Z$ and we obtain some oscillation properties of the eigenfunctions of $\mathcal{S}$.  We show that the eigenvalues of the operator satisfy an interlacing property, and we apply that property in approximating critical values of the test statistic in one of the two applications to data considered herein.  Further, we establish the consistency of the test.

In Section \ref{contiguous}, we consider the test statistic under various contiguous alternatives to the null hypothesis.  In particular, we obtain the asymptotic distribution of $T_n^2$ under gamma alternatives with varying rate or shape parameters and for a class of contaminated gamma models.  

Finally, in Section \ref{sec:efficiency}, we establish the Bahadur and Pitman efficiency properties of the statistic $T_n^2$.   We investigate the approximate Bahadur slope of $T_n^2$ under local alternatives and we establish the validity of the Wieand condition, under which the approaches through the approximate Bahadur efficiency and the Pitman efficiency are in accord.

In deriving these results, we discover two unsolved problems, one on an integral operator on a Hilbert space of square-integrable functions and a second on generalized Lambert series; we view the emergence of these problems as an on-going illumination of the symbiotic relationship between mathematical statistics and mathematical analysis.  

In closing this introduction, we acknowledge the many contributions of Professor Norbert Henze to probability, mathematical statistics, and statistics exposition.  To give a complete account of Professor Henze's work would require in and of itself a publication devoted solely to that goal.  Nevertheless, we note that his contributions to goodness-of-fit testing, as exemplified by the subset \cite{henze17,henze92,henze08,henze18,henze002,henze82,henze83,henze96,henze14,henze02,henze00} of his publications on that subject, have spanned four decades, covered broad classes of distributions, pioneered numerous distinct approaches, were published across several continents, and motivated the research of many of his co-authors.  Throughout, and by means of these activities, Professor Henze has set high standards for us all.

\section{Bessel Functions and Hankel Transforms}
\label{propertieshankel}
\setcounter{equation}{0}

\subsection{Preliminary notation and results}
\label{appendix}

Throughout the paper, all needed results on the classical special functions can be found in the books by Erd\'elyi, et al. \cite{ref5} and Olver, et al. \cite{ref1}, and we will conform to the notation used by them.  Thus, we use the standard notation, 
$$
\Gamma(\alpha) = \int_0^\infty x^{\alpha-1} e^{-x} \dd x,
$$
Re$(\alpha) > 0$, for the gamma function.   For $\alpha \in \mathbb{C}$ and $k \in \bN_0$, the set of nonnegative integers, we will make frequent use of the {\it rising factorial}, 
$$
(\alpha)_k = \frac{\Gamma(\alpha+k)}{\Gamma(\alpha)} = \alpha (\alpha+1)\cdots(\alpha+k-1),
$$

We write $X \sim Gamma (\alpha,\lambda)$ whenever a random variable $X$ is gamma-distributed with shape parameter $\alpha > 0$, rate parameter $\lambda > 0$, and probability density function, 
$$
f(x) = \frac{\lambda^\alpha \, x^{\alpha-1} \, e^{-\lambda x}}{\Gamma(\alpha)},
$$
$x > 0$.  The following result is a well-known characterization of the gamma distribution; cf., Gut \cite{gut2002} and the references given there.    

\medskip

\begin{lemma}
\label{lemma2}
Suppose that $X \sim Gamma(\alpha,1)$.  Then $E(X^k)=(\alpha)_k$, $k \in \mathbb{N}_0$.

Conversely, let $X$ be a nonnegative random variable with moments $E(X^k) = (\alpha)_k$, $k \in \bN_0$, where $\alpha > 0$. Then, $X \sim Gamma(\alpha,1)$.
\end{lemma}

\subsection{Bessel functions and generalized Laguerre polynomials}
\label{Bessel_functions}

For $\nu \in \bR$ and $-\nu \notin \bN$, the \textit{Bessel function of the first kind of order $\nu$} is defined as the infinite series, 
\begin{eqnarray}
\label{besselseriesdefinition}
J_{\nu}(z) = \sum_{j=0}^{\infty} \frac{(-1)^j}{j! \, \Gamma(\nu+1+j)} \, (z/2)^{2j+\nu},
\end{eqnarray}
$z \in \bC$.  We refer to Erd\'elyi, et al. \cite[Chapter 7]{ref5} and Olver, et al. \cite[Chapter 10]{ref1} for extensive accounts of these functions.  In particular, the series (\ref{besselseriesdefinition}) is continuous, converges absolutely for all $z \in \bC$, and converges uniformly on compact subsets of $\bC$.  

For the special case in which $\nu = -\tfrac12$, it follows from (\ref{besselseriesdefinition}) that, for $x \in \bR$, 
\begin{equation}
\label{besselminushalf}
x^{1/2} \, J_{-1/2}(x) = \Big(\frac{2}{\pi}\Big)^{1/2} \cos x,
\end{equation}
For $\nu > -1/2$, the Bessel function is also given by the \textit{Poisson integral}, 
\begin{equation}
\label{besselintegraldefinition}
J_{\nu}(x) = \frac{(x/2)^{\nu}}{\pi^{1/2} \,  \Gamma\big(\nu+\frac{1}{2}\big)} \int^\pi_0 \cos(x\cos\theta) (\sin\theta)^{2\nu} \, \dd\theta,
\end{equation}
$x \in \bR$; see \cite[(7.12(9))]{ref5}, \cite[(10.9.4)]{ref1}.  This result can be proved by expanding $\cos(x\cos\theta)$ as a power series in $x \cos (\theta)$ and integrating term-by-term.  

The Bessel function $J_\nu$ also satisfies the inequality, 
\begin{equation}
\label{besselineq1}
| J_{\nu}(z)| \le \frac{1}{\Gamma(\nu+1)} \, |z/2 |^{\nu} \, \exp(\rm{Im}(z)),
\end{equation}
$\nu \ge -1/2$, $z \in \bC$; see \cite[(7.3.2(4))]{ref5} or \cite[(10.14.4)]{ref1}.  

Henceforth, we assume that $\nu \ge -1/2$.  For $t, x \ge 0$, we set $z = 2(tx)^{1/2}$ in (\ref{besselineq1}) to obtain the inequality 
\begin{equation}
\label{besselineq2}
\big|(tx)^{-\nu/2}J_{\nu}\big(2(tx)^{1/2}\big)\big| \le \frac{1}{\Gamma(\nu+1)}.
\end{equation}

Although the following two results may be known, we have not been able to find them stated explicitly in the literature.

\begin{lemma}
\label{lemma5}
For $\nu \ge -1/2$ and $t \ge 0$, 
\begin{eqnarray}
\label{inequalitylemma5}
\big|t^{-\nu}J_{\nu+1}(t)\big| \le \frac{1}{2^{\nu} \pi^{1/2} \Gamma\big(\nu+\frac{3}{2}\big)}.
\end{eqnarray}
\end{lemma}

\Pro 
By \cite[(10.6.6)]{ref1}, 
\begin{equation}
\label{besselderiv}
t^{-\nu}J_{\nu+1}(t) = -\big(t^{-\nu}J_{\nu}(t)\big)',
\end{equation}
$t \ge 0$.  For $\nu > -1/2$, it follows by differentiating the Poisson integral (\ref{besselintegraldefinition}) that 
\begin{align*}
|t^{-\nu}J_{\nu+1}(t)| &= \frac{1}{2^\nu \pi^{1/2} \Gamma\big(\nu+\frac{1}{2}\big)} \bigg| \int^\pi_0 {\cos\theta \ \sin(t\cos\theta) \ (\sin\theta)^{2\nu}} \, \dd\theta \bigg| \\
&\le \frac{1}{2^\nu \pi^{1/2} \Gamma\big(\nu+\frac{1}{2}\big)} \int^\pi_0 {| \cos \theta | \ |(\sin\theta)^{2\nu} |} \, \dd\theta \\
&= \frac{2}{2^\nu \pi^{1/2} \Gamma\big(\nu+\frac{1}{2}\big)} \int^{\pi/2}_0 {|\cos\theta| \ |(\sin\theta)^{2\nu} |} \, \dd\theta.
\end{align*}
By a substitution, $s = \sin^2 \theta$, the latter integral reduces to a {\em beta integral},
$$
\int_0^1 s^{a-1} (1-s)^{b-1} \dd s = \frac{\Gamma(a) \, \Gamma(b)}{\Gamma(a+b)},
$$
$a, b > 0$.  This produces the result stated in (\ref{inequalitylemma5}).  

For $\nu = -1/2$, it follows from (\ref{besselderiv}) and (\ref{besselminushalf}) that 
\begin{equation}
\label{Jhalf}
t^{1/2} J_{1/2}(t) = (2/\pi)^{1/2} \, \sin t;
\end{equation}
cf. \cite[(10.16.1)]{ref1}.  Then, 
$|t^{1/2} J_{1/2}(t)| \le (2/\pi)^{1/2},$ 
as stated in (\ref{inequalitylemma5}).  
$\qed$

\begin{remark}{\rm 
Substituting $\nu = 0$ in Lemma \ref{lemma5}, we obtain $|J_1(t)| \le 2/\pi$, $t \ge 0$.  This bound is sharper than a bound given in \cite[(10.14.1)]{ref1}, {\it viz.}, $|J_1(t)| \le 2^{-1/2}$, $t \ge 0$. 
}\end{remark}

\begin{lemma}
\label{lemma_lipschitz}
For $\nu \ge -1/2$, the function $t^{-\nu}J_{\nu+1}(t)$, $t \ge 0$, is Lipschitz continuous, satisfying for $u, v \in \bR$, the inequality 
\begin{eqnarray}
\label{inequalitylemmalipschitz}
\big|u^{-\nu}J_{\nu+1}(u) - v^{-\nu}J_{\nu+1}(v)\big| \le \frac{1}{2^{\nu+1} \Gamma(\nu+2)} \, |u-v|.
\end{eqnarray}
\end{lemma}

\Pro 
As before, we apply (\ref{besselderiv}) and the Poisson integral (\ref{besselintegraldefinition}).  For $\nu > -1/2$, we apply the triangle inequality to obtain 
\begin{align*}
\big|u^{-\nu} & J_{\nu+1}(u) - v^{-\nu}J_{\nu+1}(v)\big| \\
&= \frac{1}{2^\nu \pi^{1/2} \Gamma\big(\nu+\frac{1}{2}\big)} \bigg| \int^\pi_0 \big[\sin(u\cos\theta)- \sin(v \cos\theta)\big] (\cos\theta) \ (\sin\theta)^{2\nu} \, \dd\theta \bigg| \\
&\le \frac{1}{2^\nu \pi^{1/2} \Gamma\big(\nu+\frac{1}{2}\big)} \int^\pi_0 |\sin(u\cos\theta)- \sin(v \cos\theta)| \  |\cos \theta| \  (\sin\theta)^{2\nu}  \, \dd\theta.
\end{align*}
By a well-known trigonometric identity, and the inequality $|\sin t| \le |t|$, $t \in \bR$, 
\begin{align}
|\sin(u\cos\theta)- \sin(v \cos\theta)| &= 2 \big|\sin\big(\tfrac12 (u-v) \cos \theta\big) \ \cos\big(\tfrac12 (u+v) \cos \theta\big)\big| \nonumber \\
&\le |u-v| \ |\cos \theta| \ \big|\cos\big(\tfrac12 (u+v) \cos \theta\big)\big| \nonumber \\
&\le |u-v| \ |\cos \theta|. \label{sinLipschitz}
\end{align}
Therefore, 
\begin{align*}
\big|u^{-\nu}J_{\nu+1}(u) - v^{-\nu}J_{\nu+1}(v)\big| &\le \frac{1}{2^\nu \pi^{1/2} \Gamma\big(\nu+\frac{1}{2}\big)} |u-v| \int^\pi_0 (\cos \theta)^2 \ (\sin\theta)^{2\nu} \, \dd\theta \\
&= \frac{2}{2^\nu \pi^{1/2} \Gamma\big(\nu+\frac{1}{2}\big)} |u-v| \, \int^{\pi/2}_0 (\cos \theta)^2 \ (\sin\theta)^{2\nu} \, \dd\theta.
\end{align*}
The latter integral can be reduced to a beta integral by the substitution $t = \sin^2 \theta$, and then we obtain (\ref{inequalitylemmalipschitz}).  

Finally, for $\nu = -1/2$, we apply (\ref{Jhalf}) to obtain 
\begin{align*}
\big|u^{1/2} J_{1/2}(u) - v^{1/2} J_{1/2}(v)\big| &= (2/\pi)^{1/2} \ |\sin u - \sin v| \\
&\le (2/\pi)^{1/2} \ |u-v|,
\end{align*}
where the latter inequality follows by substituting $\theta = 0$ in (\ref{sinLipschitz}).   Then, we obtain (\ref{inequalitylemmalipschitz}) for $\nu = -1/2$.   
$\qed$

\bigskip

The \textit{modified Bessel function of the first kind of order $\nu$} is defined for $-\nu \notin \bN$ as 
\begin{eqnarray*}
I_{\nu}(x) = \sum_{j=0}^{\infty} \frac{1}{j! \, \Gamma(\nu+1+j)} (x/2)^{2j+\nu},
\end{eqnarray*}
$x \in \bR$.  Letting $\mi = \sqrt{-1}$, it follows from (\ref{besselseriesdefinition}) that $I_\nu(x) = \mi^{-\nu} \, J_{\nu}(\mi x)$, $x \in \bR$.  Applying (\ref{besselineq1}), it follows that 
\begin{equation}
\label{besselineq3}
| \Gamma(\nu+1) \, (x/2)^{-\nu} \, I_\nu(x)| \le 1.
\end{equation}

Let $a, b \in \mathbb{R}$, where $-b \notin \bN_0$.  The {\em confluent hypergeometric function} is defined, for $x \in \bR$, as the infinite series, 
\begin{equation}
\label{confluenthgf}
{_1}F_1(a;b;x)=\sum_{j=0}^{\infty} \frac{(a)_j}{(b)_j} \, \frac{x^j}{j!}.
\end{equation}
We refer to \cite[Chapter 6]{erdelyi1} and \cite[Chapter 13]{ref1} for detailed accounts of this function.  
Especially, we will make repeated use of \textit{Kummer's formula}: For $x \in \bR$, 
\begin{equation}
\label{kummer}
{_1}F_1(a;b;x) = e^{x} \, {_1}F_1(b-a;b;-x).
\end{equation}

For $n \in \mathbb{N}_0$ and $\alpha > 0$, the (\textit{generalized}) \textit{Laguerre polynomial} of order $\alpha-1$ and degree $n$ is defined by 
\begin{align*}
L_n^{(\alpha-1)}(x) &= \frac{1}{n!} \, {}_1F_1(-n;\alpha;x) \\
&= \sum_{k=0}^n \frac{(\alpha+k)_{n-k}}{(n-k)!} \, \frac{(-x)^k}{ k!},
\end{align*}
$x \in \bR$; see Olver, et al. \cite[Chapter 18]{ref1} or Szeg\"o \cite[Chapter 5]{ref16}.  The \textit{normalized (generalized) Laguerre polynomial} of order $\alpha-1$ and degree $n$ is defined by 
\begin{equation}
\label{laguerre2}
\mathcal{L}_n^{(\alpha-1)}(x):=
\left(\frac{n!}{(\alpha)_n}\right)^{1/2} L_n^{(\alpha-1)}(x), 
\end{equation}
$x \in \bR$.  It is well-known (see \cite[Chapter 18.3]{ref1} or \cite[Chapter 5.1]{ref16}) that the polynomials $\mathcal{L}_n^{(\alpha-1)}$ are orthonormal with respect to the $Gamma(\alpha,1)$ distribution: 
\begin{equation}
\label{laguerreorthog}
\int^\infty_0 {\mathcal{L}_n^{(\alpha-1)}(x)\mathcal{L}_m^{(\alpha-1)}(x)\frac{x^{\alpha-1}e^{-x}}{\Gamma(\alpha)}} \, \dd x =
\begin{cases}
1, & \hbox{if } n = m \\
0, & \hbox{if } n \neq m
\end{cases}
\end{equation}

\medskip

\begin{lemma}
\label{lemma6}
For $v > 0$ and $\alpha > 0$, 
$$
\int^\infty_0 x^{\alpha} e^{-v x} L_n^{(\alpha-1)}(x) \, \dd x = \frac{\Gamma(\alpha+n)}{n!} (v-1)^{n-1} v^{-(\alpha+n+1)} \big(\alpha(v-1)-n\big).
$$
\end{lemma}

\Pro 
Starting with the known integral \cite[(18.17.34)]{ref1},
$$
\int^{\infty}_0 x^{\alpha-1}e^{-v x} L_n^{(\alpha-1)}(x) \, \dd x = \frac{\Gamma(\alpha+n)}{n!} \, (v-1)^n \, v^{-(\alpha+n)},
$$
we differentiate both sides with respect to $v$.   Simplifying the resulting expression, we obtain the stated result.
$\qed$

\subsection{Hankel transforms and some of their properties}
\label{subsechankel}

Let $X$ be a nonnegative random variable with probability density function $f(x)$.  For $\nu \ge -1/2$, we define the {\it Hankel transform of} $X$ as the function 
\begin{equation}
\label{hankeltransformdefinition}
\mathcal{H}_{X, \nu}(t) = \Gamma(\nu+1) \int^\infty_0 (tx)^{-\nu/2} \, J_\nu\big(2(tx)^{1/2}\big) \, f(x) \, \dd x,
\end{equation}
$t \ge 0$.  
The Hankel transform satisfies the following properties:

\begin{lemma}  
\label{existencehankel} For $\nu \ge -1/2$, 

\noindent
(i) $|\mathcal{H}_{X, \nu}(t)| \le 1$ for all $t \ge 0$.

\noindent
(ii) $\mathcal{H}_{X, \nu}(0)=1$.

\noindent
(iii) $\mathcal{H}_{X, \nu}(t)$ is a continuous function of $t$.
\end{lemma}

\Pro (i) By the inequality (\ref{besselineq2}) for $J_{\nu}(x)$, 
$$
\Gamma(\nu+1) \big|(tx)^{-\nu/2}J_{\nu}(2\sqrt{tx})\big| \le 1
$$
for all $x, t > 0$.  Therefore, by the triangle inequality, 
$$
|\mathcal{H}_{X, \nu}(t)| \le \int^\infty_0 f(x) \, \dd x = 1.
$$

(ii) By the series expansion (\ref{besselseriesdefinition}) of $J_\nu(x)$, we have 
\begin{equation}
\label{Hankelkernel}
\Gamma(\nu+1) \, (tx)^{-\nu/2} \, J_{\nu}\big(2(tx)^{1/2}\big) = \sum_{j=0}^{\infty} \frac{(-1)^j}{j!\, (\nu+1)_j} \, (tx)^j,
\end{equation}
for $x, t > 0$.  Therefore, 
$$
\Gamma(\nu+1) (tx)^{-\nu/2} J_{\nu}\big(2(tx)^{1/2}\big)\Big|_{t=0} = 1,
$$
so we obtain 
$$
\mathcal{H}_{X, \nu}(0) = \int^\infty_0 f(x) \, \dd x = 1.
$$

(iii) As the function $(tx)^{-\nu/2}J_{\nu}(2\sqrt{tx})$ is a power series in $tx$, it is continuous in $t \ge 0$ for every fixed $x \ge 0$.  As it is also bounded, then $\Gamma(\nu+1)(tx)^{-\nu/2}J_{\nu}(2\sqrt{tx}) f(x)$ is bounded by the Lebesgue integrable function $f(x)$ for all $x, t \ge 0$.  Therefore, the conclusion follows from the Dominated Convergence Theorem.

We also remark that an alternative proof of the continuity of $\mathcal{H}_{X, \nu}(t)$ may be obtained using an approach given by Zemanian \cite[p. 142]{ref10}.  
$\qed$

\medskip

\begin{example} 
\label{hankeltransformgammadistn}
{\rm 
Let $X \sim Gamma (\alpha,\lambda)$, where $\alpha, \lambda > 0$.  For $t \ge 0$, it follows from the definition (\ref{hankeltransformdefinition}) of the Hankel transform that 
\begin{align*}
\mathcal{H}_{X, \nu}(t) &= \frac{\Gamma(\nu+1)}{\Gamma(\alpha)} \lambda^{\alpha} \int^\infty_0 \, (tx)^{-\nu/2} \, J_{\nu}(2\sqrt{tx}) \, x^{\alpha-1} e^{-\lambda x} \, \dd x.
\end{align*}
Writing $(tx)^{-\nu/2} \, J_{\nu}(2\sqrt{tx})$ as a power series in $tx$, integrating term-by-term, and simplifying the resulting series, we obtain \begin{equation}
\label{Hankel_gamma}
\mathcal{H}_{X, \nu}(t)={_1}F_1 (\alpha;\nu+1;-t/\lambda).
\end{equation}

For the case in which $\nu = \alpha - 1$, (\ref{Hankel_gamma}) reduces to 
$$
\mathcal{H}_{X, \nu}(t)={_1}F_1(\alpha;\alpha;-t/\lambda)  = e^{-t/\lambda},
$$
$ t \ge 0$.  In particular, if $\alpha = 1$, so that $X$ has an exponential distribution with rate parameter $\lambda$, then $\mathcal{H}_{X,0}(t) = e^{-t/\lambda}$, $t \ge 0$, as shown by Baringhaus and Taherizadeh \cite[Example 2.1]{BT2010}.
}\end{example}

\medskip

\begin{example} 
\label{hankeltransformgammamixture}
{\rm 
Let $Z \sim Gamma(\alpha,1)$ and $X \ge 0$ be a random variable that is independent of $Z$.  For $t \ge 0$, 
$$
\mathcal{H}_{XZ, \nu}(t) = E_{X} \, \big[ \,{_1}F_1(\nu+1-\alpha;\nu+1; tX)\big].
$$
To prove this result, we again apply (\ref{hankeltransformdefinition}), and the independence of $X$ and $Z$, obtaining 
$$
\mathcal{H}_{XZ, \nu}(t) = E_X \, E_Z \Big[\Gamma(\nu+1) (tXZ)^{-\nu/2} J_{\nu}\big(2(tXZ)^{1/2}\big)\Big].
$$
Applying Example \ref{hankeltransformgammadistn} to calculate the expectation with respect to $Z$, we obtain 
$$
\mathcal{H}_{XZ, \nu}(t) = E_{X} \big[ \, {_1}F_1(\alpha;\nu+1;-tX)\big].
$$
In particular, if $\nu = \alpha - 1$ then $\mathcal{H}_{XZ, \nu}(t) = E_X \big[e^{-tX}\big]$, the Laplace transform of $X$, a result shown for $\nu = 0$ in \cite[Example 2.2]{BT2010}.
}\end{example}

\medskip

The following example, which provides the Hankel transform of a function related to the gamma density, will be needed repeatedly in the sequel.  

\begin{example}{\rm 
\label{lemma4}
Suppose that $X \sim Gamma(\alpha,1)$.  Then, for $t \ge 0$, 
\begin{equation}
\label{lemma4expectation}
E \big[(tX/\alpha)^{1-(\alpha/2)} \, J_\alpha\big(2(tX/\alpha)^{1/2}\big)\big] = \frac{1}{\Gamma(\alpha+1)} \, t \, e^{-t/\alpha}.
\end{equation}
Here again, we write $(tX/\alpha)^{1-(\alpha/2)} \, J_\alpha\big(2(tX/\alpha)^{1/2}\big)$ as a power series in $tX/\alpha$, integrate term-by-term, and simplify the resulting series to obtain (\ref{lemma4expectation}).
}\end{example}

\medskip

The following Hankel transform inversion theorem is a classical result which can be obtained from many sources, e.g., Sneddon \cite[p. 309, Theorem 1]{ref18}.  

\medskip

\begin{theorem}[Hankel Inversion]
\label{inversion}
Let $X$ be a non-negative, continuous random variable with probability density function $f(x)$ and Hankel transform $\mathcal{H}_{X, \nu}$.  For $x > 0$, 
$$
f(x) = \frac{1}{\Gamma(\nu+1)} \int^\infty_0 (tx)^{\nu/2} J_{\nu}(2\sqrt{tx}) \ \mathcal{H}_{X,\nu}(t) \, \dd t,
$$
\end{theorem}

\medskip

As a consequence of the Hankel inversion formula, we obtain the uniqueness of the Hankel transforms of random variables.  

\medskip

\begin{theorem}[Hankel Uniqueness]
\label{uniqueness}
Let $X$ and $Y$ be non-negative random variables with corresponding Hankel transforms $\mathcal{H}_{X, \nu}$ and $\mathcal{H}_{Y, \nu}$. Then $\mathcal{H}_{X, \nu} = \mathcal{H}_{Y, \nu}$ if and only if $X \stackrel{d}{=} Y$.
\end{theorem}

\smallskip

The next result, on the continuity of the Hankel transform, is analogous to Theorem 2.3 of Baringhaus and Taherizadeh \cite{BT2010}.  Therefore, we will omit the proof.  

\smallskip

\begin{theorem}[Hankel Continuity]
Let $\{X_{n}, n \in \mathbb{N}\}$ be a sequence of nonnegative random variables with corresponding Hankel transforms $\{\mathcal{H}_{n}, n \in \mathbb{N}\}$.  If there exists a non-negative random variable $X$, with Hankel transform $\mathcal{H}$, such that $X_{n} \xrightarrow{d} X$, then
\begin{eqnarray}
\label{continuityhankel}
\lim_{n \rightarrow \infty} \mathcal{H}_{n}(t) = \mathcal{H}(t)
\end{eqnarray}
for each $t \ge 0$.

Conversely, suppose there exists a function $\mathcal{H}: [0,\infty) \rightarrow \mathbb{R}$ such that $\mathcal{H}(0) = 1$, $\mathcal{H}$ is continuous at $0$, and (\ref{continuityhankel}) holds.  Then $\mathcal{H}$ is the Hankel transform of a nonnegative random variable $X$, and $X_{n} \xrightarrow{d} X$.
\end{theorem}

\medskip

The next result constitutes a characterization of the gamma distributions using Hankel transforms of arbitrary order $\nu$, where $\nu \ge -1/2$.  The result enables the extension, to the gamma case, of the results of Baringhaus and Taherizadeh \cite{ks_exp} on a supremum norm test statistic.  

\begin{theorem}
\label{gamma_characterization}
Let $X$ be a non-negative random variable with Hankel transform $\mathcal{H}_{X, \nu}$.  If there exist $\epsilon > 0$ and $\alpha > 0$ such that $\mathcal{H}_{X, \nu}(t) = {_1}F_1(\alpha;\nu+1;-t)$ for all $t \in [0, \epsilon]$, then $X \sim Gamma(\alpha,1)$.
\end{theorem}

We refer the reader to Hadjicosta \cite{hadjicosta19}, where three proofs of this result are given.  One proof is analogous to the proof in \cite{BT2010} for the case in which $\nu = 0$; as that proof is somewhat lengthy, we will provide the two briefer proofs.  Our first proof uses the Dominated Convergence Theorem and the moment problem for the gamma distributions; the second proof uses the principle of analytic continuation.  

\bigskip

\noindent
{\it First Proof of Theorem \ref{gamma_characterization}}.  
As the function ${}_1F_1(\alpha;\nu+1;-t)$ has a power series expansion on the interval $(0,\epsilon)$ then so does the function $\mathcal{H}_{X,\nu}(t)$.  For $\nu \ge -1/2$, the function $t^{-\nu} J_\nu(t)$, $t \in \bR$, is bounded; therefore, as the density function $f$ is integrable, we insert the series expansion (\ref{Hankelkernel}) into (\ref{hankeltransformdefinition}), and apply the Dominated Convergence Theorem to interchange the integral and sum.  As a consequence, we deduce that all moments of $X$ exist.  

Again applying (\ref{hankeltransformdefinition}) and (\ref{Hankelkernel}) we obtain, for $t \in (0,\epsilon)$, 
\begin{align*}
\sum_{j=0}^{\infty} \frac{(-1)^j}{j!\, (\nu+1)_j} \, t^j \, E(X^j) &= E \big[\Gamma(\nu+ 1) \, (tX)^{-\nu/2} \, J_\nu\big(2(tX)^{1/2}\big)\big] \\
&= \mathcal{H}_{X,\nu}(t) \\
&= {}_1F_1(\alpha;\nu+1;-t) \\
&= \sum_{j=0}^\infty \frac{(-1)^j}{j! \, (\nu+1)_j} \, (\alpha)_j \, t^j.
\end{align*}
Comparing coefficients of $t^j$, we obtain $E(X^j) = (\alpha)_j$, $j \in \bN$; therefore, $X$ has the same moment sequence as the $Gamma(\alpha,1)$ distribution.  By Lemma \ref{lemma2}, the gamma distribution is uniquely determined by its moments, so we conclude that $X \sim Gamma(\alpha,1)$.
$\qed$

\bigskip

\noindent
{\it Second Proof of Theorem \ref{gamma_characterization}}.  
The Hankel transform, $\mathcal{H}_{X,\nu}(t)$, of $X$ is real-analytic in $t$. Also, the hypergeometric function ${}_1F_1(\alpha;\nu+1;-t)$ is real-analytic in $t$. Since these two functions agree in the open neighborhood $(0,\epsilon)$ then, by analytic continuation, they agree wherever they both are well-defined. Since they both are well-defined everywhere then we conclude that $\mathcal{H}_{X,\nu}(t) = {}_1F_1(\alpha;\nu+1;-t)$ for all $t \ge 0$. By Theorem \ref{uniqueness}), the uniqueness theorem for Hankel transforms, it follows that $X \sim Gamma(\alpha,1)$.  
$\qed$

\section{Goodness-of-Fit Tests for Gamma Distributions}
\label{goodnessoffittests}
\setcounter{equation}{0}

\subsection{The test statistic}

Let $X_1, \dotsc, X_n$ be mutually independent, identically distributed (i.i.d.), nonnegative, continuous random variables with an unknown distribution $\mathcal{P}$.  We wish to test the null hypothesis, 
$
H_0: \mathcal{P} \in \{Gamma(\alpha,\lambda), \lambda>0\}
$
against the alternative hypothesis, 
$
H_1: \mathcal{P} \not\in \{Gamma(\alpha,\lambda), \lambda>0\},
$
where $\alpha$ is known.  

Since $\lambda$ is unspecified by $H_0$, the random variables $X_1,\ldots,X_n$ cannot be used directly to construct a test statistic.  Thus, with $\overline{X}_{n} = n^{-1} \sum_{i=1}^n X_j$ denoting the sample mean, define 
$$
Y_j = \begin{cases}
X_j/\overline{X}_{n}, & \hbox{if } \overline{X}_{n} \neq 0 \\
0, & \hbox{if } \overline{X}_n = 0
\end{cases}
$$
for $j=1,\dotsc,n$.  Under $H_0$, the distribution of $Y_1,\ldots,Y_n$ does not depend on $\lambda$, so a test statistic can be based on them.  Let $P_0$ denote the cumulative distribution function of the $Gamma(\alpha,1)$ distribution.  For $\nu \ge -1/2$, define the \textit{empirical Hankel transform of order $\nu$} of $Y_1,\ldots,Y_n$ as 
\begin{eqnarray}
\label{empiricalhankel}
\mathcal{H}_{n, \nu}(t)=\Gamma(\nu+1) \frac{1}{n}\sum_{j=1}^{n} (tY_{j})^{-\nu/2} J_{\nu} (2\sqrt{tY_{j}}), \end{eqnarray}
$t \ge 0$, and define the test statistic, 
\begin{eqnarray}
\label{statisticc}
T^2_{n,\nu} = n \int^{\infty}_0 \big[\mathcal{H}_{n, \nu}(t) - {_1}F_1(\alpha; \nu+1; -t/\alpha)\big]^2 
\, \dd P_0(t).
\end{eqnarray}

To provide motivation for this statistic, suppose that $H_0$ is valid; then $E(X_1) = \alpha/\lambda$ and, for large $n$, we can expect to have $Y_j = X_j/\overline{X}_n \simeq \lambda X_j/\alpha$, almost surely.  In turn, by the Continuous Mapping Theorem (see Billingsley \cite[p. 31, Corollary 1]{ref21} or Chow and Teicher \cite[p. 67, Corollary 2]{ref2}), we can expect that for each $t \ge 0$ and for sufficiently large $n$, the sequence of random variables 
$(tY_1)^{-\nu/2} J_{\nu}(2\sqrt{tY_1})$, $j=1,\dotsc,n$, 
approximates the i.i.d. sequence 
$(\lambda t X_j/\alpha)^{-\nu/2} J_{\nu} (2 (\lambda t X_1/\alpha)^{1/2})$, $j=1,\dotsc,n$.  

Applying to (\ref{empiricalhankel}) the Strong Law of Large Numbers, we then can expect that, for large $n$, 
$\mathcal{H}_{n, \nu}(t) \simeq \mathcal{H}_{X_{1}, \nu}(\lambda t/\alpha)$, almost surely.  By Example \ref{hankeltransformgammadistn} and the Hankel Uniqueness Theorem \ref{uniqueness}, $\mathcal{H}_{X_{1}, \nu}(\lambda t/\alpha) = {_1}F_1(\alpha ; \nu+1 ; -t/\alpha)$, $t \ge 0$, if and only if $H_0$ is valid.  Therefore, small values of $T_{n,\nu}^2$ provide strong evidence in support of $H_0$, and we will reject $H_0$ for large values of $T^2_{n,\nu}$.

Henceforth, set $\nu = \alpha-1$; since $\nu  \ge -1/2$ then $\alpha \ge 1/2$.  We also denote $T^2_{n, \alpha-1}$ and $\mathcal{H}_{n,\alpha-1}$ by $T^2_n$ and $\mathcal{H}_n$, respectively.  By Kummer's formula (\ref{kummer}), the statistic (\ref{statisticc}) becomes
\begin{eqnarray}
\label{statistic}
T^2_n = n \int^{\infty}_0 \big[\mathcal{H}_n(t) - e^{-t/\alpha} \big]^2 \, \dd P_0(t).
\end{eqnarray}
This integral represents $T_n^2$ as a weighted integral of the squared difference between the empirical Hankel transform $\mathcal{H}_n$ and its almost sure limit under the null hypothesis.  

\medskip

\begin{remark}
\label{remark_weight_function}
{\rm 
Although there may seem to be many possible choices for $\nu$ and the weight measure in (\ref{statistic}), the choices are limited by the practical need to calculate the test statistic explicitly and to determine its distributions and statistical properties under the null and the alternative hypotheses.  The particular choices $\nu = \alpha-1$ and the weight measure $\dd P_0(t)$ lead to satisfactory outcomes under these criteria; moreover, in the exponential case, $\alpha = 1$, the results of Baringhaus and Taherizadeh \cite{BT2010} show that the resulting statistic has superior statistical advantages.  Consequently, we believe that these choices of $\nu$ and the weight measure are likely to be optimal from combined theoretical and practical standpoints.  
}\end{remark}

We now evaluate the test statistic $T_n^2$ for a given random sample. 

\begin{proposition}
\label{helpfulrepr}
The test statistic (\ref{statistic}) is a $V$-statistic of order 2.  Specifically,
$$
T^2_{n}=\frac{1}{n}\sum_{i=1}^{n} \sum_{j=1}^{n} h(Y_{i},Y_{j})
$$
where, for $x,y \ge 0$,  
\begin{align*}
h(x,y) = \Gamma&(\alpha) \, (xy)^{(1-\alpha)/2} \, \exp(-x-y) \, I_{\alpha-1}\big(2(xy)^{1/2}\big) {\phantom{\bigg[}} \\
& - \left(\frac{\alpha}{\alpha+1}\right)^\alpha \Big[\exp\Big(-\frac{\alpha x}{\alpha+1}\Big) + \exp\Big(-\frac{\alpha y}{\alpha+1}\Big)\Big] + \left(\frac{\alpha}{\alpha+2}\right)^{\alpha}.
\end{align*}
\end{proposition}

\Pro
By squaring the integrand in (\ref{statistic}), there are three terms to be calculated.  First,
\begin{align*}
n \int^{\infty}_0 \mathcal{H}^2_{n}(t) \, \dd P_0(t) &= \frac{1}{n} \int^{\infty}_0 \left(\sum_{i=1}^{n} \Gamma(\alpha)(Y_{i}t)^{(1-\alpha)/2} \, J_{\alpha-1}(2\sqrt{tY_{i}})\right)^2 \, \dd P_0(t) \\
&= \frac{\Gamma(\alpha)}{n} \sum_{i=1}^{n}\sum_{j=1}^{n} (Y_{i}Y_{j})^{(1-\alpha)/2} \int^{\infty}_0 {J_{\alpha-1}(2\sqrt{tY_{i}})J_{\alpha-1}(2\sqrt{tY_{j}})e^{-t}} \, \dd t.
\end{align*}
Replacing $t$ by $t^2$ transforms the latter integrals into special cases of Weber's second exponential integral \cite[(10.22.67)]{ref1}: 
\begin{equation}
\label{Webersintegral}
\int^\infty_0 \, t \, \exp(-p^2 t^2) \, J_{\nu}(at) \, J_{\nu}(bt) \, \dd t = \frac{1}{2p^2} \, \exp\big(-(a^2+b^2)/4p^2\big) \, I_{\nu}(ab/2p^2),
\end{equation}
valid for $\nu > -1$ and $p, a, b \in \bR$.  Simplifying the resulting expressions, we obtain 
$$
n\int^{\infty}_0 \mathcal{H}^2_n(t) \, \dd P_0(t) = \frac{\Gamma(\alpha)}{n} \sum_{i=1}^n \sum_{j=1}^n (Y_i Y_j)^{(1-\alpha)/2} \exp(-Y_i - Y_j) \, I_{\alpha-1}\big(2(Y_i Y_j)^{1/2}\big).
$$

Second, by proceeding as in Example \ref{hankeltransformgammadistn}, we obtain 
\begin{align*}
2n \int^{\infty}_0 & \mathcal{H}_{n}(t) \, e^{-t/\alpha} \, \dd P_0(t) \\
&\equiv 2 \sum_{i=1}^{n} \int^{\infty}_0 (tY_{i})^{(1-\alpha)/2} J_{\alpha-1}\big(2(tY_{i})^{1/2}) \, \exp(-(1+\alpha^{-1})t) \, t^{\alpha-1} \, \dd t \\
&= 2\sum_{i=1}^{n} (1+\alpha^{-1})^{-\alpha} \, e^{-\alpha Y_i/(\alpha+1)}.
\end{align*}
As this latter expression is symmetric in $Y_1,\ldots,Y_n$, we write it as a double sum, obtaining  
$$
2n \int^{\infty}_0 \mathcal{H}_{n}(t) \, e^{-t/\alpha} \, \dd P_0(t) = \frac{1}{n} \sum_{i=1}^{n}\sum_{j=1}^n \left(\frac{\alpha}{\alpha+1}\right)^{\alpha} \left[e^{-\alpha Y_{i}/(\alpha+1)}+e^{-\alpha Y_{j}/(\alpha+1)} \right].
$$

Third, the last integral is a gamma integral: 
\begin{align*}
n \int^{\infty}_0 e^{-2t/\alpha} \, \dd P_0(t) &= n \left(\frac{\alpha}{\alpha+2}\right)^\alpha 
= \frac{1}{n} \sum_{i=1}^n \sum_{j=1}^n \left(\frac{\alpha}{\alpha+2}\right)^\alpha.
\end{align*}
Collecting together the three terms, we obtain the desired result.
$\qed$

\subsection{The limiting null distribution of the test statistic}
\label{sectionlimitingnull}

We denote by $L^2=L^2(P_0)$ the space of (equivalence classes of) Borel measurable functions $f:[0,\infty) \rightarrow \mathbb{C}$ that are square-integrable with respect to the probability measure $P_0$, i.e. $\int^{\infty}_0 |f(t)|^2 \, \dd P_0(t) < \infty$.  
The space $L^2$ is a separable Hilbert space when equipped with the inner product
$$
\langle f, g \rangle_{L^2} = \int^{\infty}_0 f(t) \, \overline{g(t)} \, \dd P_0(t),
$$
and the corresponding norm, 
$$
\|f\|_{L^2} = \sqrt{\langle f, f \rangle_{L^2}},
$$
$f, g \in L^2$.  Moreover, it is well-known that the normalized Laguerre polynomials $\{\mathcal{L}_n^{(\alpha-1)}: n \in \mathbb{N}_0\}$, defined in Section \ref{Bessel_functions}, form an orthonormal basis, i.e. a complete orthonormal system, for the space $L^2$; see Szeg\"o \cite[p. 108, Chapter 5.7]{ref16}.  

After these preliminaries, we define the stochastic process
\begin{equation}
\label{zn}
Z_{n}(t) = \frac{1}{\sqrt{n}} \sum_{j=1}^n \Big[\Gamma(\alpha)(tY_{j})^{(1-\alpha)/2}J_{\alpha-1}(2\sqrt{tY_{j}}) - e^{-t/\alpha} \Big],
\end{equation}
$t \ge 0$.  We will view the process $Z_{n} := \{Z_{n}(t), t \ge 0\}$ as a random element in $L^2$ since, as we will observe in Lemma \ref{lemmazn} below, its sample paths are in $L^2$.

\begin{lemma}
\label{lemmazn}
The test statistic (\ref{statistic}) can be written as 
$$
T^2_{n} = \int^{\infty}_0 \big(Z_{n}(t)\big)^2 \, \dd P_0(t) = \|Z_n\|^2_{L^2}.
$$
In particular, $\|Z_n\|^2_{L^2} < \infty$.  
\end{lemma}

\Pro
The proof follows from the definition (\ref{statistic}) of the statistic $T_n^2$ and the observation that 
$n^{1/2} \big[\mathcal{H}_n(t) - e^{-t/\alpha} \big] \equiv Z_n(t)$.
$\qed$

\begin{remark}{\rm 
An advantage of using the random variables $(Y_1, \dotsc, Y_n)$ is that the statistic $T_n^2$ is scale-invariant in $X_1,\ldots,X_n$.  It is well-known that the joint distribution of $(Y_1, \dotsc, Y_n)$ under $H_0$ is a Dirichlet distribution, which does not depend on $\lambda$.  Therefore, without loss of generality, we will set $\lambda=1$ in deriving the null distribution of $T^2_{n}$. 
}\end{remark}

\begin{theorem}
\label{limitingnulldistribution}
Let $X_1, X_2, \ldots$ be i.i.d. $Gamma(\alpha,1)$-distributed random variables, where $\alpha \ge 1/2$, and let $Z_{n} := \{Z_{n}(t), t \ge 0\}$ be the stochastic process defined in (\ref{zn}).  Then, there exists a centered Gaussian process $Z := \{Z(t), t \ge 0\}$, with sample paths in $L^2$ and with covariance function, 
\begin{equation}
\label{covariance}
K(s,t) = e^{-(s+t)/\alpha} \Big[\Gamma(\alpha)(st/\alpha^2)^{(1-\alpha)/2}I_{\alpha-1}\big(2\sqrt{st}/\alpha\big) - \alpha^{-3}st - 1\Big],
\end{equation}
$s, t \ge 0$, such that $Z_{n} \xrightarrow{\ d\ } Z$ in $L^2$ as $n \to \infty$.  Moreover,
$$
T^2_{n} \xrightarrow{\ d \ } \int^\infty_0 Z^2(t) \, \dd P_0(t).
$$
\end{theorem}


\begin{remark}{\rm 
The method of proof of Theorem \ref{limitingnulldistribution} follows an approach similar to that used by Baringhaus and Taherizadeh \cite{BT2010}.  However, the technical details are more involved here because of the general parameter $\alpha$, so we now provide a synopsis of the proof for readers who wish to postpone reading the details of the proof.  

As the random variables $Y_1,\ldots,Y_n$ are not independent, we cannot directly apply a Central Limit Theorem to deduce convergence of $Z_n$ to the Gaussian process $Z$.  Instead, we apply a standard method of constructing auxiliary processes, $Z_{n,1}$, $Z_{n,2}$, and $Z_{n,3}$, and then decomposing the process $Z_n-Z$ into a sum of four parts, viz., 
$$
Z_n - Z = (Z_n - Z_{n,1}) + (Z_{n,1} - Z_{n,2}) + (Z_{n,2} - Z_{n,3}) + (Z_{n,3} - Z).
$$
Next, we show that $Z_n - Z_{n,1}$, $Z_{n,1} - Z_{n,2}$, and $Z_{n,2} - Z_{n,3}$ each converge to $0$ in probability, in $L^2$; then we apply a Central Limit Theorem to deduce that $Z_{n,3} \xrightarrow{\ d\ } Z$ in $L^2$, and so we obtain $Z_n \xrightarrow{\ d\ } Z$ in $L^2$.  Finally, we apply the Continuous Mapping Theorem (\cite[p. 67]{ref2}, \cite[p. 31]{ref21}) to conclude that $\|Z_n\|^2 \xrightarrow{\ d\ } \|Z\|^2$.  
}\end{remark}

\smallskip

\noindent{\it Proof of Theorem \ref{limitingnulldistribution}}. \ 
Recall from (\ref{besselderiv}) that $(s^{1-\alpha} J_{\alpha-1}(s))'= -s^{1-\alpha} J_{\alpha}(s)$.  Therefore, the Taylor expansion of order $1$ of the function $s^{1-\alpha} J_{\alpha-1}(s)$, at a point $s_0$, is
\begin{equation}
\label{taylor1}
s^{1-\alpha} J_{\alpha-1}(s) = s_0^{1-\alpha} J_{\alpha-1}(s_0) + (s_0 -s) u^{1-\alpha}J_{\alpha}(u),
\end{equation}
where $u$ lies between $t$ and $t_0$.  Setting $s = 2(tY_j)^{1/2}$ and $s_0 = 2(tX_j/\alpha)^{1/2}$, we obtain 
\begin{equation}
\label{taylor}
\begin{aligned}
2^{1-\alpha}(tY_j)^{(1-\alpha)/2} \, J_{\alpha-1}\big(2(tY_{j})^{1/2}\big) &= 2^{1-\alpha} (tX_{j}/\alpha)^{(1-\alpha)/2} J_{\alpha-1}\big(2(tX_{j}/\alpha)^{1/2}\big) \\
& \qquad + 2 \big[(tX_{j}/\alpha)^{1/2} - (tY_{j})^{1/2}\big] \, u_j^{1-\alpha} \, J_{\alpha}(u_j),
\end{aligned}
\end{equation}
where $u_j$ lies between $2(tY_{j})^{1/2}$and $2(tX_{j}/\alpha)^{1/2}$.  Define 
\begin{equation}
\label{Wn}
W_n = \alpha^{-1/2} - \overline{X}_n^{-1/2} = \frac{\overline{X}_n-\alpha}{(\alpha \overline{X}_n)^{1/2} (\alpha^{1/2}+\overline{X}_n^{1/2})};
\end{equation}
then 
$$
(tX_{j}/\alpha)^{1/2} - (tY_{j})^{1/2} = (tX_{j}/\alpha)^{1/2} - (tX_{j}/\overline{X}_n)^{1/2} = W_n \, (tX_{j})^{1/2},
$$
and (\ref{taylor}) reduces to 
\begin{multline}
2^{1-\alpha}(tY_j)^{(1-\alpha)/2} \, J_{\alpha-1}\big(2(tY_{j})^{1/2}\big) \\
= 2^{1-\alpha}(tX_{j}/\alpha)^{(1-\alpha)/2} \, J_{\alpha-1}\big(2(tX_{j}/\alpha)^{1/2}\big) 
\label{taylornew}
+ 2 W_n \, (tX_{j})^{1/2} \, u_j^{1-\alpha} \, J_{\alpha}(u_j).
\end{multline}
Multiplying both sides of (\ref{taylornew}) by $2^{\alpha-1}$, and then adding and subtracting the term 
$$
2 (tX_{j})^{1/2} W_n \, (tX_{j}/\alpha)^{(1-\alpha)/2} \, J_{\alpha}\big(2(tX_{j}/\alpha)^{1/2}\big)
$$
on the right-hand side, we obtain
\begin{equation}
\label{taylor2}
\begin{aligned}
(tY_j)&^{(1-\alpha)/2} \, J_{\alpha-1}\big(2(tY_{j})^{1/2}\big) \\
= \ & (tX_{j}/\alpha)^{(1-\alpha)/2} \, J_{\alpha-1}\big(2(tX_{j}/\alpha)^{1/2}\big) \\
& + 2^{\alpha} (tX_{j})^{1/2} \, W_n \, \big(2(tX_{j}/\alpha)^{1/2}\big)^{1-\alpha} \, J_{\alpha}\big(2(tX_{j}/\alpha)^{1/2}\big) \\
& + 2^{\alpha} (tX_{j})^{1/2} \, W_n \, \Big(u_j^{1-\alpha}J_{\alpha}(u_j)-\big(2(tX_{j}/\alpha)^{1/2}\big)^{1-\alpha} \, J_{\alpha}\big(2(tX_{j}/\alpha)^{1/2}\big)\Big).
\end{aligned}
\end{equation}
Simplifying the second term of the right-hand side of equation (\ref{taylor2}), we obtain 
\begin{equation}
\label{taylorexp}
\begin{aligned}
(tY_j)&^{(1-\alpha)/2} \, J_{\alpha-1}\big(2(tY_{j})^{1/2}\big) \\
= \ & (tX_{j}/\alpha)^{(1-\alpha)/2} \, J_{\alpha-1}\big(2(tX_{j}/\alpha)^{1/2}\big) \\
& + 2\alpha^{1/2} \, W_n \, (tX_{j}/\alpha)^{1-(\alpha/2)} \, J_{\alpha}\big(2(tX_{j}/\alpha)^{1/2}\big) \\
& + 2^{\alpha} \, W_n \, (tX_{j})^{1/2} \, \Big(u_j^{1-\alpha} \, J_{\alpha}(u_j)-\big(2(tX_{j}/\alpha)^{1/2}\big)^{1-\alpha} \, J_{\alpha}\big(2(tX_{j}/\alpha)^{1/2}\big)\Big).
\end{aligned}
\end{equation}

We now define the processes $Z_{n,1}(t)$, $Z_{n,2}(t)$, and  $Z_{n,3}(t)$, $t \ge 0$, by 
\begin{align*}
Z_{n,1}(t) &= \frac{1}{\sqrt{n}} \sum_{j=1}^{n} \Big[\Gamma(\alpha) \, (tX_{j}/\alpha)^{(1-\alpha)/2} \, J_{\alpha-1}\big(2(tX_{j}/\alpha)^{1/2}\big)\\
& \quad\quad\quad\quad\quad  + 2\Gamma(\alpha) \alpha^{1/2} \, W_n \, (tX_{j}/\alpha)^{1-(\alpha/2)} \, J_{\alpha}\big(2(tX_{j}/\alpha)^{1/2}\big)-e^{-t/\alpha} \Big], \\
Z_{n,2}(t) &= \frac{1}{\sqrt{n}}\sum_{j=1}^{n} \Big[\Gamma(\alpha) \, (tX_{j}/\alpha)^{(1-\alpha)/2} \, J_{\alpha-1}\big(2(tX_{j}/\alpha)^{1/2}\big) + 2 \alpha^{-1/2} \, W_n \, t e^{-t/\alpha} - e^{-t/\alpha} \Big], \\
Z_{n,3}(t) &= \frac{1}{\sqrt{n}}\sum_{j=1}^{n}\Big[\Gamma(\alpha) \, (tX_{j}/\alpha)^{(1-\alpha)/2} \, J_{\alpha-1}\big(2(tX_{j}/\alpha)^{1/2}\big) +\frac{(X_j -\alpha)}{\alpha^2} \, te^{-t/\alpha} - e^{-t/\alpha}\Big].
\end{align*}
The processes $Z_{n,k}(t)$, $k=1,2,3$ arise as follows.  
To define $Z_{n,1}(t)$, we use the first two terms in (\ref{taylorexp}), multiplied by $\Gamma(\alpha)$. To define $Z_{n,2}(t)$, we use the same expression from $Z_{n,1}(t)$ except that the term $(tX_{j}/\alpha)^{1-(\alpha/2)} \, J_{\alpha}(2(tX_{j}/\alpha)^{1/2})$ is replaced by its expected value, which is given in Example \ref{lemma4}.  To define $Z_{n,3}(t)$, we replace the term $W_n$ in $Z_{n,2}(t)$ by a constant multiple of $X_j-\alpha$, the constant being obtained by applying the Law of Large Numbers to the denominator of (\ref{Wn}).  We will show that
\begin{eqnarray}
\label{5.3.1}&Z_{n,3} \xrightarrow{d} Z \ \ \text{in $L^2$},&\\
\label{5.3.2}&\lVert Z_{n}-Z_{n,1} \rVert_{L^2} \xrightarrow{p} 0,&\\
\label{5.3.3}&\lVert Z_{n,1}-Z_{n,2} \rVert_{L^2} \xrightarrow{p} 0,&\\
\label{5.3.4}&\lVert Z_{n,2}-Z_{n,3} \rVert_{L^2} \xrightarrow{p} 0.&
\end{eqnarray}

To establish (\ref{5.3.1}), let 
\begin{equation}
\label{Zn3j}
Z_{n,3,j}(t) := \Gamma(\alpha) \, (tX_{j}/\alpha)^{(1-\alpha)/2} \, J_{\alpha-1}\big(2(tX_{j}/\alpha)^{1/2}\big) + \frac{(X_{j}-\alpha)}{\alpha^2} te^{-t/\alpha} - e^{-t/\alpha},
\end{equation}
for $t \ge 0$ and $j=1,\dotsc,n$. 
Since $X_j \sim Gamma(\alpha,1)$ then $E(X_j-\alpha)=0$; also, by Example \ref{hankeltransformgammadistn} and (\ref{confluenthgf}), we deduce that 
\begin{align*}
E \big[\Gamma(\alpha)(tX_{j}/\alpha)^{(1-\alpha)/2}J_{\alpha-1} \big(2(tX_{j}/\alpha)^{1/2}\big) \big] & = {_1}F_1(\alpha ; \alpha ; -t/\alpha) = e^{-t/\alpha}.
\end{align*}
Therefore, $E(Z_{n,3,j}(t))=0$, for all $t \ge 0$ and $j=1,\dotsc,n$, and it is also clear that $Z_{n,3,1},\dotsc,Z_{n,3,n}$ are i.i.d. random elements in $L^2$.  

We now show that $E(\lVert Z_{n,3,j} \rVert^2_{L^2}) < \infty$ for $j=1,\dotsc,n$.  
We have 
\begin{align*}
E(\lVert Z_{n,3,j} \rVert^2_{L^2}) &= E \int^\infty_0 {Z_{n,3,j}^2(t)}\, \dd P_0(t) \\
&= E \int^\infty_0 {\Big[\Gamma(\alpha)(tX_{j}/\alpha)^{(1-\alpha)/2}J_{\alpha-1}\big(2(tX_{j}/\alpha)^{1/2}\big)} \\
& \qquad\qquad\qquad\qquad\qquad + { \frac{(X_{j}-\alpha)}{\alpha^2} te^{-t/\alpha}-e^{-t/\alpha} \Big]^2 }\, \dd P_0(t).
\end{align*}
By the Cauchy-Schwarz inequality, $(a+b+c)^2 \le 3(a^2+b^2+c^2)$, for $a, b, c$ $\in \mathbb{R}$; so, to prove that $E(\lVert Z_{n,3,j} \rVert^2_{L^2}) < \infty$, it suffices to prove that
\begin{equation}
\label{Znlimitnormineq1}
E \int^\infty_0 {\big[ \Gamma(\alpha)(tX_{j}/\alpha)^{(1-\alpha)/2}J_{\alpha-1} \big(2(tX_{j}/\alpha)^{1/2} \big) \big]^2}\, \dd P_0(t)  < \infty,
\end{equation}
\begin{equation}
\label{Znlimitnormineq2}
E\Big[ (X_j-\alpha)^2 \int^\infty_0 \Big(\frac{t e^{- t/\alpha}}{\alpha^2}\Big)^2\, \dd P_0(t) \Big]< \infty,
\end{equation}
and 
\begin{equation}
\label{Znlimitnormineq3}
E\Big[ \int^{\infty}_0 (e^{-t/\alpha})^2 \, \dd P_0(t)  \Big] < \infty.
\end{equation}
To establish (\ref{Znlimitnormineq1}), we apply the inequality (\ref{besselineq2}) to obtain 
$$
E \int^\infty_0 { \big[ \Gamma(\alpha)(tX_{j}/\alpha)^{(1-\alpha)/2}J_{\alpha-1} \big( 2(tX_{j}/\alpha)^{1/2} \big) \big]^2 }\, \dd P_0(t)  \ \le \ E\int^\infty_0 \, 1 \cdot \dd P_0(t) =1. 
$$
As for (\ref{Znlimitnormineq2}) and (\ref{Znlimitnormineq3}), substituting for $\dd P_0(t)$ shows that both expectations are constant multiples of convergent gamma integrals.  

In summary, $E(Z_{n,3,j}(t)) = 0$ for $t \ge 0$ and $j=1,\dotsc,n$, and $Z_{n,3,1}(t),\dotsc,Z_{n,3,n}(t)$ are i.i.d. random elements in $L^2$ with $E(\lVert Z_{n,3,1} \rVert^2_{L^2}) < \infty$.  Therefore, by the Central Limit Theorem in $L^2$ (cf., Ledoux and Talagrand \cite[p. 281, Theorem 10.5]{ref19}), 
$$
\frac{1}{\sqrt{n}}\sum_{j=1}^{n} Z_{n,3,j} \xrightarrow{d} Z,
$$
where $Z:=(Z(t),t \ge 0)$ is a centered Gaussian random element in $L^2$. Moreover, $Z$ has the same covariance operator as $Z_{n,3,1}$.  

It is well-known that the covariance operator of the random element $Z_{n,3,1}$ is uniquely determined by the covariance function of the stochastic process $Z_{n,3,1}(t)$; see G{\=\i}khman and Skorohod \cite[pp.~218-219]{ref20}. We now show that the function $K(s,t)$ in (\ref{covariance}) is the covariance function of $Z_{n,3,1}$.  Noting that $E[Z_{n,3,1}(t)] = 0$ for all $t$, we obtain 
\begin{align*}
K(s,t) &= \Cov \big[Z_{n,3,1}(s),Z_{n,3,1}(t)\big] \\
&= \Cov \big[Z_{n,3,1}(s)+e^{-s/\alpha},Z_{n,3,1}(t)+e^{-t/\alpha}\big] \\
&= E \big[\big(Z_{n,3,1}(s)+e^{-s/\alpha}\big) \big(Z_{n,3,1}(t) + e^{-t/\alpha}\big)\big] - e^{-(s+t)/\alpha}.
\end{align*}
By (\ref{Zn3j}), 
\begin{equation}
\begin{aligned}
\label{fourterms}
E \big(Z_{n,3,1}(s) + & e^{-s/\alpha}\big) \big(Z_{n,3,1}(t) + e^{-t/\alpha}\big) \\
&= E\Big[\Gamma(\alpha) \, (sX_1/\alpha)^{(1-\alpha)/2} \, J_{\alpha-1}\big(2(sX_1/\alpha)^{1/2}\big) + \frac{(X_1-\alpha)}{\alpha^2} se^{-s/\alpha}\Big] \\
& \quad\times \Big[\Gamma(\alpha) \, (tX_1/\alpha)^{(1-\alpha)/2} \, J_{\alpha-1}\big(2(tX_1/\alpha)^{1/2}\big) + \frac{(X_1-\alpha)}{\alpha^2} te^{-t/\alpha}\Big],
\end{aligned}
\end{equation}
so the calculation of $K(s,t)$ reduces to evaluating the four terms obtained by expanding the product on the right-hand side of (\ref{fourterms}).  

The first term in the product in (\ref{fourterms}) is 
\begin{align}
\label{covarianceterm1}
E \, [\Gamma(\alpha)]^2 & (sX_{1}/\alpha)^{(1-\alpha)/2}(tX_{1}/\alpha)^{(1-\alpha)/2}J_{\alpha-1}(2(sX_{1}/\alpha)^{1/2})J_{\alpha-1}(2(tX_{1}/\alpha)^{1/2}) \nonumber\\
&= \Gamma(\alpha)(st/\alpha^2)^{(1-\alpha)/2}\int^\infty_0 {J_{\alpha-1}(2(sx/\alpha)^{1/2}) J_{\alpha-1}(2(tx/\alpha)^{1/2}) e^{-x} }\, \dd x.
\end{align}
Replacing $x$ by $x^2$, and applying Weber's second exponential integral (\ref{Webersintegral}), we find that the first term equals 
\begin{equation}
\label{K0kernelcalculation}
\Gamma(\alpha)(st/\alpha^2)^{(1-\alpha)/2}e^{-(s+t)/\alpha}I_{\alpha-1}(2\sqrt{st}/\alpha).
\end{equation}
The second term in the product in (\ref{fourterms}) is 
\begin{align*}
E\Big[ \Gamma(\alpha)&(sX_{1}/\alpha)^{(1-\alpha)/2}J_{\alpha-1}(2(sX_{1}/\alpha)^{1/2}) \frac{(X_{1}-\alpha)}{\alpha^2} te^{-t/\alpha} \Big] \\
&=\frac{t e^{-t /\alpha} s^{(1-\alpha)/2}}{\alpha^{(5-\alpha)/2}}\Big [ \int^\infty_0 {x^{(\alpha+1)/2}J_{\alpha-1}(2(sx/\alpha)^{1/2})e^{-x}}\, \dd x \\
& \quad\quad\quad\quad\quad\quad\quad\quad\quad\quad -\alpha \int^\infty_0 {x^{(\alpha-1)/2}J_{\alpha-1}(2(sx/\alpha)^{1/2})e^{-x}}\, \dd x \Big].
\end{align*}
The latter integrals are Hankel transforms of the type given in Example \ref{hankeltransformgammadistn}.  Applying that result to calculate each integral and simplifying the outcome, we deduce that the second term equals $- \alpha^{-3}s t \exp\big(-(s+t)/\alpha\big)$.  

The third term in the product in (\ref{fourterms}) is 
$$
E\Big[\Gamma(\alpha) \, (tX_{1}/\alpha)^{(1-\alpha)/2} \, J_{\alpha-1}(2(tX_{1}/\alpha)^{1/2}) \frac{(X_{1}-\alpha)}{\alpha^2} se^{-s/\alpha} \Big],
$$
which is the same as the second term but with $s$ and $t$ interchanged.  Therefore, the third term also equals $- \alpha^{-3}s t \exp\big(-(s+t)/\alpha\big)$.  

Finally, the fourth term in the product in (\ref{fourterms}) is 
$$
E\Big[\frac{(X_{1}-\alpha)^2}{\alpha^4}ste^{-(s+t)/\alpha}\Big] = \frac{s t e^{-(s+t)/\alpha}}{\alpha^4} \Var(X_1) = \frac{ste^{-(s+t)/\alpha}}{\alpha^3}.
$$
Combining all four terms, we obtain (\ref{covariance}).   

To establish (\ref{5.3.2}), we begin by showing that
$$
(\sqrt{n} W_n)^2 = \bigg( \frac{\sqrt{n} (\overline{X}_n-\alpha)}{(\alpha\overline{X}_n)^{1/2} (\alpha^{1/2}+\overline{X}_n^{1/2})} \bigg)^2
$$
converges in distribution to a random variable with finite variance. By the Central Limit Theorem, 
$
\sqrt{n}(\overline{X}_{n}-\alpha) \xrightarrow{d} \mathcal{N}(0,\alpha).
$  
Also, by the Law of Large Numbers and the Continuous Mapping Theorem, 
$$
(\alpha\overline{X}_n)^{1/2} (\alpha^{1/2}+\overline{X}_n^{1/2}) \xrightarrow{p} \alpha(\alpha^{1/2}+\alpha^{1/2}) = 2\alpha^{3/2}.
$$
By Slutsky's theorem (\cite[p. 249, Corollary 2]{ref2}), $\sqrt{n} W_n \xrightarrow{d} \mathcal{N}(0,\tfrac14 \alpha^{-2})$, so it follows from the Continuous Mapping Theorem that $(\sqrt{n} W_n)^2 \xrightarrow{d} \chi^2_1/4\alpha^2$, where $\chi_1^2$ denotes a random variable which has a chi-square distribution with one degree of freedom.  In particular, the limiting distribution of $\sqrt{n} W_n$ has finite variance.  

By the Taylor expansion (\ref{taylorexp}), 
\begin{align*}
Z_{n}&-Z_{n,1}\\
& =\frac{\Gamma(\alpha)}{\sqrt{n}}\sum_{j=1}^{n}\Big[ (tY_{j})^{(1-\alpha)/2}J_{\alpha-1}(2(tY_{j})^{1/2})-(tX_{j}/\alpha)^{(1-\alpha)/2}J_{\alpha-1}(2(tX_{j}/\alpha)^{1/2})\\
& {\hskip 2.6truein} -2 \alpha^{1/2} \, W_n \, (tX_{j}/\alpha)^{1-(\alpha/2)} \, J_{\alpha}\big(2(tX_{j}/\alpha)^{1/2}\big)\Big]\\
&=\frac{2^{\alpha} \Gamma(\alpha)}{n} (\sqrt{n} W_n) \ \sum_{j=1}^{n}  (tX_{j})^{1/2} \, \Big[u_j^{1-\alpha} \, J_{\alpha}(u_j)-\big(2(tX_{j}/\alpha)^{1/2}\big)^{1-\alpha} \, J_{\alpha}\big(2(tX_{j}/\alpha)^{1/2}\big)\Big].
\end{align*}
Define 
$$
V_n:=\frac{1}{n^2} \int^\infty_0 \Big[\sum_{j=1}^{n} (tX_{j})^{1/2} \, \Big(u_j^{1-\alpha} \, J_{\alpha}(u_j)-\big(2(tX_{j}/\alpha)^{1/2}\big)^{1-\alpha} \, J_{\alpha}\big(2(tX_{j}/\alpha)^{1/2}\big)\Big) \Big]^2 \, \dd P_0(t).
$$
Then, 
$$
\lVert Z_{n}-Z_{n,1} \rVert^2_{L^2} = 4^{\alpha} [\Gamma(\alpha)]^2 (\sqrt{n} W_n)^2 \, V_n.
$$
By the Cauchy-Schwarz inequality, 
$$
V_n \le \frac{1}{n} \int^\infty_0 { t \sum_{j=1}^{n} X_{j} \ \big| u_j^{1-\alpha} \, J_{\alpha}(u_j)-\big(2(tX_{j}/\alpha)^{1/2}\big)^{1-\alpha} \, J_{\alpha}\big(2(tX_{j}/\alpha)^{1/2}\big) \big| ^2}\, \dd P_0(t).
$$
Recall that $u_j$ lies between $2(tY_{j})^{1/2}$ and $2(tX_{j}/\alpha)^{1/2}$, so we can write 
\begin{align*}
u_j &= 2 (1-\theta_{n,j,t}) (tX_{j}/\alpha)^{1/2} + 2 \theta_{n,j,t} (tY_{j})^{1/2} \\
&= 2(tX_{j})^{1/2} \big( \alpha^{-1/2}+\theta_{n,j,t}(\overline{X}_n^{-1/2}-\alpha^{-1/2})\big),
\end{align*}
where $\theta_{n,j,t} \in [0,1]$.  By Lemma \ref{lemma_lipschitz}, the Lipschitz property of the Bessel functions, 
\begin{align*}
\big| u_j^{1-\alpha} \, J_{\alpha}(u_j) - \big(2(tX_j/\alpha&)^{1/2}\big)^{1-\alpha} \, J_{\alpha}\big(2(tX_{j}/\alpha)^{1/2}\big) \big|^2 \\
&\le \big|u_j- 2 (tX_j/\alpha)^{1/2} \big|^2 \big/ 4^{\alpha} [\Gamma(\alpha+1)]^2 \\
&= \big| 2(tX_{j})^{1/2} \theta_{n,j,t}(\overline{X}_n^{-1/2}-\alpha^{-1/2})\big|^2 \big/ 4^{\alpha} [\Gamma(\alpha+1)]^2 \\
&\le tX_{j} \, (\overline{X}_n^{-1/2}-\alpha^{-1/2})^2 \big/ 4^{\alpha-1} [\Gamma(\alpha+1)]^2,
\end{align*}
since $\theta_{n,j,t} \in [0,1]$. Therefore, 
$$
V_n \le \frac{1}{4^{\alpha-1} [\Gamma(\alpha+1)]^2} \Big(\frac{1}{n} \sum_{j=1}^n X^2_j \Big) (\overline{X}_n^{-1/2}-\alpha^{-1/2})^2 \int^{\infty}_0 {t^2}\, \dd P_0(t). 
$$
By the Law of Large Numbers, $(\overline{X}_n^{-1/2}-\alpha^{-1/2})^2 \xrightarrow{p} 0$ and 
$
n^{-1} \sum_{j=1}^n X^2_j \xrightarrow{p} E(X^2_1)=\alpha(\alpha+1),
$ 
so it follows that $V_{n} \xrightarrow{p} 0$. By Slutsky's theorem,
$$
\lVert Z_{n}-Z_{n,1} \rVert^2_{L^2} = 4^{\alpha} [\Gamma(\alpha)]^2 (\sqrt{n} W_n)^2 \cdot V_n \xrightarrow{d} 0,
$$
therefore $\lVert Z_{n}-Z_{n,1} \rVert_{L^2} \xrightarrow{p} 0$, as asserted in (\ref{5.3.2}).    

To establish (\ref{5.3.3}), define 
$$
\Delta_j(t):= \Gamma(\alpha)(tX_{j}/\alpha)^{1-(\alpha/2)} J_{\alpha}\big( 2(tX_{j}/\alpha)^{1/2} \big) - \alpha^{-1} te^{-t/\alpha},
$$
$t \ge 0$, $j=1,\dotsc,n$.  Then it is straightforward to verify that 
$$
Z_{n,1}-Z_{n,2} = \frac{2\alpha^{1/2}}{\sqrt{n}} \, W_n \,\sum_{j=1}^{n} \Delta_j(t)
$$
and therefore 
\begin{equation}
\label{ZN1minusZn2}
\lVert Z_{n,1} - Z_{n,2} \rVert^2_{L^2} = (2\alpha^{1/2} W_n)^2 \, \int^\infty_0 \Big[ \frac{1}{\sqrt{n}} \sum_{j=1}^n \Delta_j(t) \Big]^2 \, \dd P_0(t).
\end{equation}

By the Law of Large Numbers, $W_n \xrightarrow{p} 0$.  Also, as shown in Example \ref{lemma4}, 
$$
E \big[ \Gamma(\alpha)(tX_{j}/\alpha)^{1-(\alpha/2)} J_{\alpha}\big( 2(tX_{j}/\alpha)^{1/2} \big] = \alpha^{-1} t e^{-t/\alpha};
$$
therefore $E(\Delta_j(t)) = 0$, $t \ge 0$, $j=1,\dotsc,n$.  Also, $\Delta_1(t),\dotsc,\Delta_n(t)$ are i.i.d. random elements in $L^2$.  

We now show that $E(\lVert\Delta_1 \rVert^2_{L^2}) < \infty$.  We have 
\begin{align*}
E(\lVert \Delta_1 \rVert^2_{L^2}) 
&=E\int^\infty_0 {\Delta_1^2(t)}\, \dd P_0(t) \\
&=E\int^\infty_0 {\Big[  \Gamma(\alpha)(tX_{1}/\alpha)^{1-(\alpha/2)} J_{\alpha}\big( 2(tX_{1}/\alpha)^{1/2} \big) - \alpha^{-1} te^{-t/\alpha} \Big]^2}\, \dd P_0(t).
\end{align*}
To show that $E(\lVert\Delta_1 \rVert^2_{L^2}) < \infty$ it suffices, by the Cauchy-Schwarz inequality, to prove that
\begin{equation}
\label{Zn1Zn2inequality1}
E \int^\infty_0 {\Big[ \Gamma(\alpha) (tX_1/\alpha)^{1-(\alpha/2)} \, J_{\alpha} \big( 2(tX_1/\alpha)^{1/2} \big) \Big]^2}\, \dd P_0(t)  \ < \ \infty
\end{equation}
and 
\begin{equation}
\label{Zn1Zn2inequality2}
E \int^\infty_0 (\alpha^{-1} te^{-t/\alpha})^2 \, \dd P_0(t) \ < \ \infty.
\end{equation}
To establish (\ref{Zn1Zn2inequality1}), we apply the inequality (\ref{inequalitylemma5}) to obtain 
$$
|J_{\alpha} \big( 2(tX_1/\alpha)^{1/2} \big)| \le \frac{(tX_1/\alpha)^{-(1-\alpha)/2}}{\pi^{1/2} \Gamma(\alpha+\frac{1}{2})},
$$
for $t \ge 0$.  Therefore, 
\begin{align*}
E \int^\infty_0 \Big[ \Gamma(\alpha)(tX_1/\alpha)^{1-(\alpha/2)} & \, J_{\alpha} \big( 2(tX_1/\alpha)^{1/2} \big) \Big]^2 \, \dd P_0(t)\\
& \le \Big( \frac{\Gamma(\alpha)}{\pi^{1/2} \Gamma(\alpha+\frac{1}{2})} \Big)^2 \ E(X_1/\alpha) \ \int^\infty_0 {t}\, \dd P_0(t) < \infty.
\end{align*}
As for (\ref{Zn1Zn2inequality2}), that expectation is a convergent gamma integral.  Hence, $E(\lVert\Delta_1 \rVert^2_{L^2}) < \infty$.  

By the Central Limit Theorem in $L^2$, 
$n^{-1/2} \sum_{j=1}^{n} \Delta_j(t)$ converges to a centered Gaussian random element in $L^2$.  Thus, by the Continuous Mapping Theorem, 
$$
\Big\lVert \frac{1}{\sqrt{n}} \sum_{j=1}^{n} \Delta_j(t)\Big\rVert ^2_{L^2} 
:= \int^\infty_0 \Big[ \frac{1}{\sqrt{n}} \sum_{j=1}^n \Delta_j(t) \Big]^2 \, \dd P_0(t)
$$
converges in distribution to a random variable which has finite variance.  Since $W_n \xrightarrow{p} 0$ then  by (\ref{ZN1minusZn2}) and Slutsky's Theorem, we obtain $\lVert Z_{n,1}-Z_{n,2} \rVert^2_{L^2} \xrightarrow{d} 0$; therefore, $\lVert Z_{n,1}-Z_{n,2} \rVert_{L^2} \xrightarrow{p} 0$.  

To prove (\ref{5.3.4}), we observe that 
\begin{align*}
Z_{n,2} - Z_{n,3} &= \frac{1}{\sqrt{n}} \sum_{j=1}^{n} \Big( 2\alpha^{-1/2} W_n te^{-t/\alpha}-\frac{(X_{j}-\alpha)}{\alpha^2} te^{-t/\alpha} \Big)\\
&= \frac{te^{-t/\alpha}}{\sqrt{n}} \Big( \frac{2n(\overline{X}_n-\alpha)}{\alpha\overline{X}_n^{1/2}(\alpha^{1/2}+\overline{X}_n^{1/2})}-\frac{n(\overline{X}_n-\alpha)}{\alpha^2} \Big)\\
&= te^{-t/\alpha} \, \sqrt{n}(\overline{X}_n-\alpha) R_n,
\end{align*}
where 
$$
R_n = \frac{2}{\alpha \overline{X}_n^{1/2}(\alpha^{1/2}+\overline{X}_n^{1/2})}-\frac{1}{\alpha^2}.
$$
Therefore, 
\begin{align*}
\lVert Z_{n,2}-Z_{n,3} \rVert^2_{L^2} &= \big[ \sqrt{n}(\overline{X}_n-\alpha) R_n \big]^2 \ \int^\infty_0 {(te^{-t/\alpha})^2}\, \dd P_0(t).
\end{align*}
As noted earlier, $\int^\infty_0 (te^{-t/\alpha})^2\, \dd P_0(t) < \infty$.  Also, by the Central Limit Theorem, $\sqrt{n}(\overline{X}_n-\alpha) \xrightarrow{d} \mathcal{N}(0,\alpha)$; and by the Law of Large Numbers, $R_n \xrightarrow{p} 0$.  By Slutsky's theorem, $\big[ \sqrt{n}(\overline{X}_n-\alpha) R_n \big]^2 \xrightarrow{d} 0$; hence $\big[ \sqrt{n}(\overline{X}_n-\alpha) R_n \big]^2 \xrightarrow{p} 0$, and therefore $\lVert Z_{n,2}-Z_{n,3} \rVert_{L^2} \xrightarrow{p} 0$.

Finally, by the Continuous Mapping Theorem in $L^2$, 
$\lVert Z_{n} \rVert^2_{L^2} \xrightarrow{d} \lVert Z \rVert^2_{L^2}$, 
i.e.
$$
T^2_{n}= \int^\infty_0 {Z^2_{n}(t)}\, \dd P_0(t)  \xrightarrow{d} \int^\infty_0 {Z^2(t)}\, \dd P_0(t) .
$$
The proof now is complete.  
$\qed$

\subsection{Eigenvalues and eigenfunctions of the covariance operator}
\label{eigen_values_functions}

The covariance operator $\mathcal{S}:L^2 \rightarrow L^2$ of the random element $Z$ is defined for $s \ge 0$ and $f \in L^2$ by 
$$
\mathcal{S} f(s) = \int^\infty_0 K(s,t) f(t) \, \dd P_0(t),
$$
where $K(s,t)$ is the covariance function defined in equation (\ref{covariance}).  Let $\{\delta_k: k \ge 1\}$ be the positive eigenvalues, listed in non-increasing order, of $\mathcal{S}$; also, let $\{\chi^2_{1k}: k \ge 1\}$ be i.i.d. $\chi^2_{1}$-distributed random variables.  It is well-known that the integrated squared process, $\int^\infty_0 Z^2(t) \, \dd P_0(t)$, has the same distribution as $\sum_{k=1}^\infty \delta_k\chi^2_{1k}$.  This result follows from the Karhunen-Lo\'eve expansion of the Gaussian process $Z(t)$; see Le Ma{\^\i}tre and Knio \cite[Chapter 2]{lemaitreknio} or Vakhania \cite[p. 58]{ref22}.  Therefore, the limiting null distribution of $T^2_n$ is the same as $\sum_{k=1}^\infty \delta_k \chi^2_{1k}$.

For $s, t \ge 0$, define 
\begin{equation}
\label{K0kernel}
K_0(s,t) = e^{-(s+t)/\alpha} \Gamma(\alpha)(st/\alpha^2)^{(1-\alpha)/2} I_{\alpha-1}\big(2\sqrt{st}/\alpha\big),
\end{equation}
the first term in the covariance function defined in equation (\ref{covariance}); by (\ref{covarianceterm1}) and (\ref{K0kernelcalculation}), 
\begin{align}
\label{K0integral}
K_0(s,t) &= \int^\infty_0 \Gamma(\alpha) (sx/\alpha)^{(1-\alpha)/2} J_{\alpha-1}(2(sx/\alpha)^{1/2}) \nonumber \\
& \quad\quad\quad\quad \times \Gamma(\alpha)(tx/\alpha)^{(1-\alpha)/2} J_{\alpha-1}(2(tx/\alpha)^{1/2}) \, \dd P_0(x).
\end{align}

We will first find the eigenvalues and eigenfunctions of the integral operator $\mathcal{S}_0: L^2 \rightarrow L^2$, defined for $s \ge 0$ and $f$ in $L^2$ by 
\begin{eqnarray}
\label{defeigen}
\mathcal{S}_0 f(s)=\int^\infty_0 K_0(s,t) f(t) \, \dd P_0(t).
\end{eqnarray}

Before presenting the results on the eigenvalues and eigenfunctions of $\mathcal{S}_0$, we state for the sake of completeness some preliminary definitions pertaining to (linear) operators on $L^2$. Note that these definitions are provided by Sunder \cite{ref23} or Young \cite{ref17}.

An operator $\mathcal{T}: L^2 \rightarrow L^2$ is called \textit{symmetric (self-adjoint)} if, for all $f, g \in L^2$, $\langle \mathcal{T} f, g \rangle_{L^2} = \langle f,\mathcal{T} g \rangle_{L^2}$.   A symmetric operator $\mathcal{T}$ is called \textit{positive} if $\langle \mathcal{T} f, f \rangle_{L^2} \ge 0$ for all $f \in L^2$. 

An operator $\mathcal{T}$ is called \textit{compact} if for every bounded sequence $\{f_k : k \in \bN \}$ in $L^2$, the sequence $\{ \mathcal{T} f_k : k \in \bN \}$ has a convergent subsequence in $L^2$. The set of eigenvalues of a compact operator is countable. 

An operator $\mathcal{T}$ is \textit{Hilbert-Schmidt} if for every orthonormal basis $\{f_k : k \in \bN \}$ in $L^2$, the series $\sum_{k=1}^{\infty} \|\mathcal{T} f_k \|^2_{L^2}$ converges.  It is well-known that Hilbert-Schmidt operators are compact \cite[ p. 93, Theorem 8.8]{ref17}.

An operator $\mathcal{T}$ is \textit{of trace class} if for every orthonormal basis $\{f_k : k \in \bN \}$ in $L^2$, the series $\sum_{k=1}^{\infty} \| \mathcal{T} f_k \|_{L^2}$ converges.  An operator $\mathcal{T}$ is trace-class if and only if it can be expressed as a product of two Hilbert-Schmidt operators \cite[p. 74, Proposition 3.3.7]{ref23}.  Further, trace-class operators are Hilbert-Schmidt. 

\smallskip

Recall that $\alpha \ge 1/2$.  Throughout the remainder of the paper, we use the notation 
\begin{equation}
\label{betaandbalpha}
\beta = \Big(\frac{\alpha+4}{\alpha}\Big)^{1/2} \quad \hbox{and} \quad b_{\alpha} = \big(1+\tfrac12\alpha(1-\beta)\big)^{1/2}.
\end{equation}
We also set 
\begin{equation}
\label{defrho}
\rho_k = \alpha^{\alpha} b_{\alpha}^{4k+2\alpha},
\end{equation} 
$k \in \mathbb{N}_0$, and 
\begin{equation}
\label{deflen}
\textgoth{L}_k^{(\alpha-1)}(s) = \beta^{\alpha/2}\exp((1-\beta)s/2) \mathcal{L}_k^{(\alpha-1)}(\beta s),
\end{equation}
$s \ge 0$.  

\begin{theorem}
\label{theoremreducedcov}
The set $\{(\rho_k,\textgoth{L}_k^{(\alpha-1)}): k \in \mathbb{N}_0\}$ is a complete enumeration of the eigenvalues and eigenfunctions, respectively, of $\mathcal{S}_0$, and the eigenfunctions $\{\textgoth{L}_k^{(\alpha-1)}: k \in \mathbb{N}_0\}$ form an orthonormal basis in $L^2$. Moreover, $\mathcal{S}_0$ is positive and of trace-class.
\end{theorem}

\Pro
Recall from \cite[Eq. (18.18.27)]{ref1} the \textit{Poisson kernel}: For $r \in (0,1)$ and $x,y \ge 0$, 
\begin{align}
\label{poissonkernel}
I_{\alpha-1}\Big(\frac{2\sqrt{xyr}}{1-r}\Big) = (1-r) \, \exp&\Big(\frac{r(x+y)}{1-r}\Big) \, (xyr)^{(\alpha-1)/2} \nonumber \\
& \times \frac{1}{\Gamma(\alpha)}\sum_{k=0}^{\infty} \mathcal{L}_k^{(\alpha-1)}(x)  \mathcal{L}_k^{(\alpha-1)}(y) r^k.
\end{align}
In this expansion, set 
\begin{equation}
\label{r}
r = b_{\alpha}^4 = \big(1+\tfrac12\alpha(1-\beta)\big)^2,
\end{equation}
so that $r \in (0,1)$.  Note that $r^{1/2} = 1+\tfrac12\alpha(1-\beta)$ satisfies the quadratic equation 
$$
r - (\alpha+2)r^{1/2} + 1 = 0
$$
and also that this equation is equivalent to the identity 
\begin{align}
\frac{1-r}{\alpha r^{1/2}} &= 1 + \frac{2}{\alpha}(1 - r^{1/2}). \nonumber \\
\intertext{On the right-hand side of this identity, substitute for $r^{1/2}$ in terms of $\alpha$ and $\beta$ to obtain }
\label{defbeta}
\frac{1-r}{\alpha r^{1/2}} &= 1 + \frac{2}{\alpha}\big[1 - \big(1+\tfrac12\alpha(1-\beta)\big)\big] = \beta.
\end{align}
In (\ref{poissonkernel}), also set 
\begin{equation}
\label{xandy}
x = \frac{1-r}{\alpha r^{1/2}} s \equiv \beta s \quad \hbox{and} \quad y = \frac{1-r}{\alpha r^{1/2}} t \equiv \beta t.
\end{equation}
Then, 
\begin{equation}
\label{xyridentities}
\frac{\sqrt{xyr}}{1-r} = \frac{\sqrt{st}}{\alpha} \quad \hbox{and} \quad \frac{r(x+y)}{1-r} = \frac{(r^{1/2}-1)(s+t)}{\alpha} + \frac{(s+t)}{\alpha}
\end{equation}
Applying (\ref{defrho}),(\ref{deflen}) and (\ref{r})-(\ref{xyridentities}) to (\ref{poissonkernel}), and substituting the result in (\ref{K0kernel}), we obtain for $s, t \ge 0$, the pointwise convergent series expansion, 
\begin{equation}
\label{reducedcovariancepoisson2}
K_0(s,t) = \sum_{k=0}^{\infty} \rho_k \, \textgoth{L}_k^{(\alpha-1)}(s) \textgoth{L}_k^{(\alpha-1)}(t).
\end{equation}
Note that the system $\{\textgoth{L}_k^{(\alpha-1)}: k \in \mathbb{N}_0\}$ is orthonormal in $L^2$, i.e., 
\begin{equation}
\label{orthonormal}
 \int^\infty_0 {\textgoth{L}_k^{(\alpha-1)} (s) \textgoth{L}_l^{(\alpha-1)} (s)}\,\dd P_0(s) = \begin{cases}
    1, & k=l \\
    0, & k \neq l
  \end{cases}
\end{equation}
This result follows from the orthonormality property (\ref{laguerreorthog}) of the generalized Laguerre polynomials $\{\mathcal{L}_k^{(\alpha-1)}: k \in \mathbb{N}_0\}$.  

We verify that the series (\ref{reducedcovariancepoisson2}) converges in the separable tensor product Hilbert space $L^2 \otimes L^2 := L^2(P_0 \otimes P_0)$.  By the Cauchy criterion, it suffices to prove that for each $\epsilon > 0$, there exists $N \in \mathbb{N}$ such that 
\begin{eqnarray*}
\int_0^\infty\!\! \int_0^\infty {\Big[\sum_{k=m_1}^{m_2} \rho_k \, \textgoth{L}_k^{(\alpha-1)}(s) \textgoth{L}_k^{(\alpha-1)}(t) \Big]^2} \, \dd (P_0 \otimes P_0) (s,t) < \epsilon,
\end{eqnarray*}
for all $m_1, m_2 \in \mathbb{N}$ such that $m_2 \ge m_1 \ge N$.  By squaring the integrand, it suffices by Fubini's theorem to consider 
$$
\sum_{k_1=m_1}^{m_2} \sum_{k_2=m_1}^{m_2} \rho_{k_1} \rho_{k_2} \int_0^\infty \!\! \int_0^\infty \prod_{j=1}^2 {\textgoth{L}_{k_j}^{(\alpha-1)}(s) \textgoth{L}_{k_j}^{(\alpha-1)}(t)} \, \dd P_0(s) \,\dd P_0(t).
$$
Since the system $\{\textgoth{L}_k^{(\alpha-1)}: k \in \mathbb{N}_0\}$ is orthonormal, the latter sum reduces to 
$$
\sum_{k=m_1}^{m_2} \rho_k^2 = \rho_0^2 \sum_{k=m_1}^{m_2} b_{\alpha}^{8k}.
$$
By (\ref{r}), $0 < b_\alpha < 1$, so the geometric series $\sum_{k=0}^{\infty} b_{\alpha}^{8k}$ converges.  Since every convergent series is Cauchy, it follows that for each $\epsilon > 0$, there exists $N \in \mathbb{N}$ such that 
$\sum_{k=m_1}^{m_2} b_{\alpha}^{8k} < \epsilon$ 
for all $m_1, m_2 \in \mathbb{N}$ such that $m_2 \ge m_1 \ge N$.  Therefore, the series (\ref{reducedcovariancepoisson2}) is Cauchy in $L^2 \otimes L^2$ and hence, 
\begin{eqnarray*}
\lim_{m \rightarrow \infty} \int_0^\infty\!\! \int_0^\infty {\Big[ K_0(s,t) - \sum_{k=0}^m \rho_k \, \textgoth{L}_k^{(\alpha-1)}(s) \textgoth{L}_k^{(\alpha-1)}(t) \Big]^2}\, \dd (P_0 \otimes P_0) (s,t)=0.
\end{eqnarray*}
By Fubini's theorem, the latter expression equals 
\begin{eqnarray}
\label{seriesoconvarianceconv}
\lim_{m \rightarrow \infty} \int^{\infty}_0 \!\! {\int^{\infty}_0 {\Big[ K_0(s,t) - \sum_{k=0}^m \rho_k \, \textgoth{L}_k^{(\alpha-1)}(s) \textgoth{L}_k^{(\alpha-1)}(t) \Big]^2}}\,\dd P_0(s) \dd P_0(t) = 0.
\end{eqnarray}
It follows from the orthonormality property, (\ref{orthonormal}), of the system $\{\textgoth{L}_k^{(\alpha-1)}: k \in \mathbb{N}_0\}$, that for $l, m \in \mathbb{N}$ with $l \le m$, 
\begin{eqnarray}
\label{propertyeigen}
\int^{\infty}_0 { \Big[ \sum_{k=0}^m \rho_k \, \textgoth{L}_k^{(\alpha-1)}(s) \textgoth{L}_k^{(\alpha-1)}(t)\Big] \textgoth{L}_l^{(\alpha-1)}(t)}\,\dd P_0(t) = \rho_l \, \textgoth{L}_l^{(\alpha-1)}(s).
\end{eqnarray}
By (\ref{defeigen}) and (\ref{propertyeigen}), 
\begin{align*}
\int^{\infty}_0 & \Big| \mathcal{S}_0\textgoth{L}_l^{(\alpha-1)}(s) - \rho_l \, \textgoth{L}_l^{(\alpha-1)}(s) \Big| \,\dd P_0(s) \\
&= \int^{\infty}_0 {\Big| \int^{\infty}_0 { \Big[ K_0(s,t) - \sum_{k=0}^{m} \rho_k \textgoth{L}_k^{(\alpha-1)}(s) \textgoth{L}_k^{(\alpha-1)}(t)\Big] \ \textgoth{L}_l^{(\alpha-1)}(t)} \,\dd P_0(t) \Big|}\,\dd P_0(s).
\end{align*}
By the Cauchy-Schwarz inequality, this latter expression is bounded by 
\begin{align}
\Big( \int^{\infty}_0 \!\! \int^{\infty}_0  \Big| K_0(s,t) - \sum_{k=0}^{m} \rho_k & \, \textgoth{L}_k^{(\alpha-1)}(s) \textgoth{L}_k^{(\alpha-1)}(t) \Big|^2 \,\dd P_0 (s) \dd P_0 (t) \Big)^{1/2} \nonumber \\
\label{limiteigen}
& \times  \Big( \int^{\infty}_0 \!\! \int^{\infty}_0 {| \textgoth{L}_l^{(\alpha-1)}(t)| ^2} \,\dd P_0 (s) \dd P_0 (t) \Big)^{1/2}.
\end{align}
By the orthonormality property (\ref{orthonormal}) and the fact that $P_0$ is a probability distribution, the second term in (\ref{limiteigen}) equals $1$; therefore, 
\begin{multline}
\label{limiteigen2}
\int^{\infty}_0 { \Big| \mathcal{S}_0\textgoth{L}_l^{(\alpha-1)}(s) - \rho_l \, \textgoth{L}_l^{(\alpha-1)}(s) \Big| }\,\dd P_0(s) \\
\le \Big( \int^{\infty}_0 \!\! \int^{\infty}_0 { \Big| K_0(s,t) - \sum_{k=0}^{m} \rho_k \, \textgoth{L}_k^{(\alpha-1)}(s) \textgoth{L}_k^{(\alpha-1)}(t) \Big|^2 } \,\dd P_0 (s) \dd P_0 (t) \Big)^{1/2}.
\end{multline}
Since $m$ is arbitrary, we now let $m \rightarrow \infty$.  By (\ref{seriesoconvarianceconv}), the right-hand side of (\ref{limiteigen2}) converges to $0$, so we obtain 
$$
\int^{\infty}_0 { \Big| \mathcal{S}_0\textgoth{L}_l^{(\alpha-1)}(s) - \rho_l \, \textgoth{L}_l^{(\alpha-1)}(s) \Big| }\,\dd P_0(s) = 0,
$$
which proves that $\mathcal{S}_0\textgoth{L}_l^{(\alpha-1)}(s) = \rho_l \, \textgoth{L}_l^{(\alpha-1)}(s)$ for $P_0$-almost every $s$. Therefore, $\rho_k$ is an eigenvalue of $\mathcal{S}_0$ with corresponding eigenfunction $\textgoth{L}_k^{(\alpha-1)}$.  

Since the kernel $K_0(s,t)$ is symmetric in $(s,t)$, it follows that $\mathcal{S}_0$ is symmetric.  To show that $\mathcal{S}_0$ is positive, we observe that for $f \in L^2$, 
\begin{align*}
\langle \mathcal{S}_0 f, f \rangle_{L^2} &= \int_0^\infty S_0f(s) \, \overline{f(s)} \, \dd P_0(s) \\
&= \int^{\infty}_0 \Big[ \int^{\infty}_0 {K_0(s,t)f(t)}\, \dd P_0(t) \Big] \, \overline{f(s)} \,\dd P_0(s).
\end{align*}
Substituting for $K_0(s,t)$ from (\ref{K0integral}), we obtain 
\begin{align*}
\langle \mathcal{S}_0 f, f \rangle_{L^2} &= \int^{\infty}_0 { \Big[ \int^{\infty}_0 { \Big( \int^{\infty}_0 {\big( \Gamma(\alpha)(sx/\alpha)^{(1-\alpha)/2} J_{\alpha-1}(2(sx/\alpha)^{1/2}) \big)}}}\\
& \quad\quad\times  \Gamma(\alpha) (tx/\alpha)^{(1-\alpha)/2} J_{\alpha-1}(2(tx/\alpha)^{1/2}) \,\dd P_0(x) \Big) f(t) \,\dd P_0(t) \Big] \overline{f(s)} \,\dd P_0(s).
\end{align*}
Applying Fubini's theorem to reverse the order of integration, we find that the inner integrals with respect to $s$ and $t$ are complex conjugates of each other; therefore, 
\begin{multline}
\label{innerproduct}
\langle \mathcal{S}_0 f, f \rangle_{L^2} \\
= \int^{\infty}_0 \Big| \int^{\infty}_0 {\Gamma(\alpha) (sx/\alpha)^{(1-\alpha)/2} J_{\alpha-1}(2(sx/\alpha)^{1/2}) f(s) }\,\dd P_0(s)  \Big|^2 \,\dd P_0(x),
\end{multline}
which is nonnegative. Therefore, $\mathcal{S}_0$ is positive. 

Next, we prove that $\mathcal{S}_0$ is of trace-class. For $f \in L^2$, $s \ge 0$, it again follows by (\ref{K0integral}) and Fubini's theorem that 
\begin{align}
\mathcal{S}_0f(s)&=\int^\infty_0 {K_0(s,t)f(t)}\, \dd P_0(t)\nonumber\\
&=\int^\infty_0 \!\! \int^\infty_0 {\Gamma(\alpha)(tx/\alpha)^{(1-\alpha)/2} J_{\alpha-1}(2(tx/\alpha)^{1/2}) f(t)}\,\dd P_0(t) \nonumber\\
\label{redcovariancetrace}
& \quad\quad\quad\quad\quad\quad\quad\quad\quad \times {\Gamma(\alpha) (sx/\alpha)^{(1-\alpha)/2} J_{\alpha-1}(2(sx/\alpha)^{1/2}) }\,\dd P_0(x).
\end{align}
Denote by $\mathcal{T}_0:L^2 \rightarrow L^2$ the integral operator, 
$$
\mathcal{T}_0f(t) = \int_0^\infty \Gamma(\alpha)(tx/\alpha)^{(1-\alpha)/2} J_{\alpha-1}(2(tx/\alpha)^{1/2}) \, f(x) \, \dd x,
$$
$t \ge 0$.  By (\ref{besselineq2}), we have that
$$
|\Gamma(\alpha)(tx/\alpha)^{(1-\alpha)/2} J_{\alpha-1}(2(tx/\alpha)^{1/2})| \le 1, 
$$
for $t,x \ge 0$; therefore, 
$$
\left\| \Gamma(\alpha)(tx/\alpha)^{(1-\alpha)/2} J_{\alpha-1}(2(tx/\alpha)^{1/2})\right\|_{L^2 \otimes L^2}^2 < \infty, 
$$
for $t,x \ge 0$. By Young \cite[p. 93, Theorem 8.8]{ref17}, it follows that $\mathcal{T}_0$ is a Hilbert-Schmidt operator.  Now, we can write (\ref{redcovariancetrace}) as
\begin{align*}
\mathcal{S}_0 f(s)&=\int^\infty_0 {\mathcal{T}_0 f(x) \Big[ \Gamma(\alpha) (sx/\alpha)^{(1-\alpha)/2} J_{\alpha-1}(2(sx/\alpha)^{1/2}) \Big] }\,\dd P_0(x)\\
&=\mathcal{T}_0(\mathcal{T}_0 f)(s),
\end{align*}
$s \ge 0$, which proves that $\mathcal{S}_0$ is of trace-class.  

To complete the proof, we now show that the set $\{\textgoth{L}_k^{(\alpha-1)}: k \in \mathbb{N}_0\}$ is complete.  Here, it suffices to show that if $f \in L^2$ with $\langle f, \textgoth{L}_k^{(\alpha-1)} \rangle_{L^2} = 0$ for all $k$ then $f=0$, $P_0$-almost everywhere.  First, we note that 
\begin{align}
\int^{\infty}_0 &\Big| S_0f(s) \, \overline{f(s)}  - \sum_{k=0}^{m} \rho_k \langle f, \textgoth{L}_k^{(\alpha-1)} \rangle _{L^2} \textgoth{L}_k^{(\alpha-1)}(s) \, \overline{f(s)} \Big| \,\dd P_0(s) \nonumber\\
&= \int^{\infty}_0 {\Big| \int^{\infty}_0 { \Big[ K_0(s,t) - \sum_{k=0}^{m} \rho_k \textgoth{L}_k^{(\alpha-1)}(s) \textgoth{L}_k^{(\alpha-1)}(t) \Big] f(t)\, \overline{f(s)}} \,\dd P_0(t) \Big| }\,\dd P_0(s) \qquad\qquad \nonumber \\
\label{completenesseigen}
&\le \Big( \int^{\infty}_0 \!\! { \int^{\infty}_0 { \Big| K_0(s,t) - \sum_{k=0}^{m} \rho_k \textgoth{L}_k^{(\alpha-1)}(s) \textgoth{L}_k^{(\alpha-1)}(t) \Big|^2 }}\,\dd P_0(s) \dd P_0 (t) \Big) ^{1/2} \nonumber\\
& {\hskip 1.75truein} \times  \Big( \int^{\infty}_0 \!\! { \int^{\infty}_0 {|f(s)|^2 |f(t)|^2}}\,\dd P_0(s) \dd P_0 (t) \Big)^{1/2},
\end{align}
by the Cauchy-Schwarz inequality.  Since $f \in L^2$, the second term on the right-hand side of (\ref{completenesseigen}) is finite. Taking the limit on both sides of (\ref{completenesseigen}) as $m \rightarrow \infty$ and applying (\ref{seriesoconvarianceconv}), we obtain 
\begin{align}
\label{completenesseigen4}
\lim_{m \rightarrow \infty} \int^{\infty}_0 { \Big| S_0f(s) \, \overline{f(s)}-\sum_{k=0}^{m} \rho_k \langle f, \textgoth{L}_k^{(\alpha-1)} \rangle _{L^2} \textgoth{L}_k^{(\alpha-1)}(s) \, \overline{f(s)} \Big| }\,\dd P_0(s)=0.
\end{align}
Since $\langle f, \textgoth{L}_k^{(\alpha-1)} \rangle _{L^2}=0$ for all $k$ then (\ref{completenesseigen4}) reduces to 
$$
\langle \mathcal{S}_0 f, f \rangle_{L^2} \equiv \int^{\infty}_0 {\mathcal{S}_0f(s)\, \overline{f(s)} }\,\dd P_0(s) = 0.
$$
Therefore, by (\ref{innerproduct}), we obtain for $P_0$-almost every $x$, 
\begin{align}
\label{completenesseigen2}
\int^{\infty}_0 {\Gamma(\alpha) (sx/\alpha)^{(1-\alpha)/2} J_{\alpha-1}(2(sx/\alpha)^{1/2}) f(s) }\,\dd P_0(s) = 0.
\end{align}
Since the function $\Gamma(\alpha)(sx/\alpha)^{(1-\alpha)/2} J_{\alpha-1}(2(sx/\alpha)^{1/2})$ is continuous in $x \ge  0$ for each fixed $s \ge 0$ and, by (\ref{besselineq2}), 
$$
|\Gamma(\alpha) (sx/\alpha)^{(1-\alpha)/2} J_{\alpha-1}(2(sx/\alpha)^{1/2}) | \le 1,
$$
$x,s \ge 0$, then by the Dominated Convergence Theorem, the integral on the left-hand side of (\ref{completenesseigen2}) is a continuous function of $x$. If two continuous functions are equal $P_0$-almost everywhere then they are equal everywhere; this result follows from the fact that a continuous function is uniquely determined by its values on any dense set. Therefore, (\ref{completenesseigen2}) holds for all $x \ge 0$.

Henceforth, without loss of generality, we assume that $f$ is real-valued. Let $f^{+}$ and $f^{-}$ denote the positive and negative parts of $f$, respectively.  Then, $f = f^{+} - f^{-}$, $f^+$ and $f^-$ are nonnegative, and since $f \in L^2$ then by the Cauchy-Schwarz inequality, $f^{+}$ and $f^{-}$ are $P_0$-integrable.  Also, by (\ref{completenesseigen2}), 
\begin{align*}
\int^{\infty}_0 & {\Gamma(\alpha) (sx/\alpha)^{(1-\alpha)/2} J_{\alpha-1}(2(sx/\alpha)^{1/2}) f^{+}(s) }\,\dd P_0(s) \\
& =\int^{\infty}_0 {\Gamma(\alpha) (sx/\alpha)^{(1-\alpha)/2} J_{\alpha-1}(2(sx/\alpha)^{1/2}) f^{-}(s) }\,\dd P_0(s),
\end{align*}
$x \ge 0$.  By the Uniqueness Theorem for Hankel transforms, we notice that there are only two possible cases. Either
$$
\int^{\infty}_0 {f^{+} (s)}\,\dd P_0(s)=\int^{\infty}_0 {f^{-} (s)}\,\dd P_0(s)=0,
$$
or 
$$
\int^{\infty}_0 {f^{+} (s)}\,\dd P_0(s)=\int^{\infty}_0 {f^{-} (s)}\,\dd P_0(s)=c > 0.
$$
For the first case, we have that $f^{+}=f^{-}=0$ and so $f=0$ $P_0$-almost everywhere.  As for the second case, we have
\begin{align*}
\int^{\infty}_0 \Gamma(\alpha) & (sx/\alpha)^{(1-\alpha)/2} J_{\alpha-1}(2(sx/\alpha)^{1/2}) \frac{1}{c} f^{+}(s) \,\dd P_0(s) \\
& =\int^{\infty}_0 {\Gamma(\alpha) (sx/\alpha)^{(1-\alpha)/2} J_{\alpha-1}(2(sx/\alpha)^{1/2}) \frac{1}{c} f^{-}(s) }\,\dd P_0(s),
\end{align*}
$x \ge 0$.  By the Uniqueness Theorem for Hankel transforms, we obtain $f^{+}=f^{-}$ and hence $f = 0$ $P_0$-almost everywhere.  This proves that the orthonormal set $\{\textgoth{L}_k^{(\alpha-1)}: k \in \mathbb{N}_0\}$ is complete, and therefore it forms a basis in the separable Hilbert space $L^2$.  
$\qed$

\medskip

We remark that an alternative proof of completeness can be obtained by following the argument given by Szeg\"o \cite[p. 108, Chapter 5.7]{ref16}.  

\bigskip

The proof of the following theorem is similar to the proof of Theorem \ref{theoremreducedcov}, and the complete details are provided by Hadjicosta \cite{hadjicosta19}.  

\begin{theorem} 
\label{thmcovarianceS}
Let $\mathcal{S}:L^2 \rightarrow L^2$ be the covariance operator of the random element $Z$ defined as
$$
\mathcal{S} f(s)=\int^\infty_0 {K(s,t) f(t)}\, \dd P_0(t),
$$
for all $s \ge 0$ and for all functions $f$ in $L^2$, where $K(s,t)$ is the covariance function defined in equation (\ref{covariance}). Then, $\mathcal{S}$ is positive and of trace-class. 
\end{theorem}

Recall that a non-trivial function $f \in L^2$ is an \textit{eigenfunction} of $\mathcal{S}$ if there exists an \textit{eigenvalue} $\delta \in \bC$ such that $\mathcal{S}f = \delta f$.  As $\mathcal{S}$ is self-adjoint and positive, its eigenvalues are real and nonnegative.  In the next result, we find the positive eigenvalues and corresponding eigenfunctions of the operator $\mathcal{S}$, and we will show in Subsection \ref{decayrate} that $0$ is not an eigenvalue of $\mathcal{S}$.  

\begin{theorem}
\label{thmeigenvaluesofS}
For  $\delta \in \bR$, $\delta \neq \rho_k$ for any $k \in \bN$, define the functions 
\begin{align*}
A(\delta) &= 1-\beta^{\alpha}\sum_{k=0}^{\infty} \frac{(\alpha)_k}{k! (\rho_k-\delta)}\rho_k^2,\\
B(\delta) &= 1-\alpha \beta^{\alpha} \sum_{k=0}^{\infty} \frac{(\alpha)_k}{k! (\rho_k-\delta)}\rho_k^2 (b_{\alpha}^2-k\beta)^2,\\
\intertext{and}
D(\delta) &= \alpha^2 \beta^{\alpha}\sum_{k=0}^{\infty} \frac{(\alpha)_k}{k! (\rho_k-\delta)}\rho_k^2(b_{\alpha}^2-k\beta).
\end{align*}
Then the positive eigenvalues of the operator $\mathcal{S}$ are the positive roots of the function $G(\delta) := \alpha^3 A(\delta)B(\delta)-D^2(\delta)$.  Moreover, the eigenfunction corresponding to an eigenvalue $\delta$ has the Fourier-Laguerre expansion
$$
\beta^{\alpha/2} \sum_{k=0}^{\infty} \frac{\rho_k}{\rho_k-\delta}\Big( \frac{(\alpha)_k}{k!} \Big)^{1/2} \big( c_1+c_2\alpha^{-1}(b_{\alpha}^2-k\beta) \big) \textgoth{L}_k^{(\alpha-1)},
$$
where $c_1, c_2$ are not both equal to 0, $\alpha^3 c_1A(\delta)=c_2 D(\delta)$, and  $c_2B(\delta)=c_1D(\delta)$.
\end{theorem}

\Pro
Since the set $\{\textgoth{L}_k^{(\alpha-1)}: k \in \mathbb{N}_0\}$ is an orthonormal basis for $L^2$, the eigenfunction $\phi \in L^2$ corresponding to an eigenvalue $\delta$ can be written as
$$
\phi=\sum_{k=0}^{\infty} \langle \phi, \textgoth{L}_k^{(\alpha-1)} \rangle_{L^2} \textgoth{L}_k^{(\alpha-1)}.
$$
We restrict ourselves temporarily to eigenfunctions for which this series is pointwise convergent.  Substituting this series into the equation $\mathcal{S} \phi =\delta \phi $, we obtain 
\begin{align}
\label{eigen}
\int^\infty_0 {K(s,t) \sum_{k=0}^{\infty} \langle \phi, \textgoth{L}_k^{(\alpha-1)} \rangle_{L^2} \textgoth{L}_k^{(\alpha-1)}(t)}\, \dd P_0(t)&= \delta \sum_{k=0}^{\infty} \langle \phi, \textgoth{L}_k^{(\alpha-1)} \rangle_{L^2} \textgoth{L}_k^{(\alpha-1)}(s).
\end{align}
Substituting the covariance function $K(s,t)$ in the left-hand side of (\ref{eigen}), writing $K$ in terms of $K_0$, and assuming that we can interchange the order of integration and summation, we obtain 
\begin{multline}
\label{equationeigen}
\delta \sum_{k=0}^{\infty} \langle \phi, \textgoth{L}_k^{(\alpha-1)} \rangle_{L^2} \textgoth{L}_k^{(\alpha-1)}(s) \\
= \sum_{k=0}^{\infty} \langle \phi, \textgoth{L}_k^{(\alpha-1)}\rangle_{L^2} \int^\infty_0 \Big[K_0(s,t) - e^{-(s+t)/\alpha} (\alpha^{-3}st + 1)\Big] \textgoth{L}_k^{(\alpha-1)}(t) \, \dd P_0(t).
\end{multline}
By Theorem \ref{theoremreducedcov}, 
$$
\int^\infty_0 K_0(s,t) \textgoth{L}_k^{(\alpha-1)}(t) \, \dd P_0(t) = \rho_k \textgoth{L}_k^{(\alpha-1)}(s).
$$
On writing $\textgoth{L}_k^{(\alpha-1)}$ in terms of $L_k^{(\alpha-1)}$, the generalized Laguerre polynomial, applying the well-known formula \cite[(18.17.34)]{ref1} for the Laplace transform of $L_k^{(\alpha-1)}$, and making use of (\ref{r}) and (\ref{defbeta}), we obtain 
\begin{equation}
\label{18.17.34second}
 \langle e^{-t/\alpha},\textgoth{L}_k^{(\alpha-1)} \rangle_{L^2}:=\int^\infty_0 {e^{-t/\alpha}\textgoth{L}_k^{(\alpha-1)}(t) }\, \dd P_0(t) = \Big( \frac{(\alpha)_k}{k!} \Big)^{1/2} \beta^{\alpha/2} \rho_k.
\end{equation}
Again writing $\textgoth{L}_k^{(\alpha-1)}$ in terms of $L_k^{(\alpha-1)}$, applying Lemma \ref{lemma6}, and (\ref{r}) and (\ref{defbeta}), we obtain 
\begin{equation}
\label{secondintegral2}
\langle te^{-t/\alpha}, \textgoth{L}_k^{(\alpha-1)} \rangle_{L^2} :=\int^\infty_0{te^{-t/\alpha}\textgoth{L}_k^{(\alpha-1)}(t) }\, \dd P_0(t) = \Big( \frac{(\alpha)_k}{k!} \Big)^{1/2} \alpha^{2} \beta^{\alpha/2}\rho_k (b_{\alpha}^2-k\beta).
\end{equation}
In summary, (\ref{equationeigen}) reduces to 
\begin{multline}
\label{equationeigensummary}
\delta \sum_{k=0}^{\infty} \langle \phi, \textgoth{L}_k^{(\alpha-1)} \rangle_{L^2} \textgoth{L}_k^{(\alpha-1)}(s) \\
= \sum_{k=0}^{\infty} \rho_k \langle \phi, \textgoth{L}_k^{(\alpha-1)}\rangle_{L^2} \Big[\textgoth{L}_k^{(\alpha-1)}(s) - e^{-s/\alpha} \Big( \frac{(\alpha)_k}{k!} \Big)^{1/2} \beta^{\alpha/2} \big(\alpha^{-1}s (b_{\alpha}^2-k\beta) + 1\big)\Big].
\end{multline}

By applying (\ref{18.17.34second}), we also obtain the Fourier-Laguerre expansion of $e^{-s/\alpha}$ with respect to the orthonormal basis $\{\textgoth{L}_k^{(\alpha-1)} : k \in \mathbb{N}_0 \}$; indeed, 
\begin{align}
e^{-s/\alpha} &= \sum_{k=0}^{\infty} \langle e^{-s/\alpha},\textgoth{L}_k^{(\alpha-1)} \rangle_{L^2} \textgoth{L}_k^{(\alpha-1)} (s)\nonumber\\
\label{fourier1}
&= \beta^{\alpha/2} \sum_{k=0}^{\infty}\Big( \frac{(\alpha)_k}{k!} \Big)^{1/2} \rho_k \textgoth{L}_k^{(\alpha-1)} (s).
\end{align}
Similarly, by applying (\ref{secondintegral2}), we have 
\begin{align}
se^{-s/\alpha} &= \sum_{k=0}^{\infty} \langle se^{-s/\alpha}, \textgoth{L}_k^{(\alpha-1)} \rangle_{L^2} \textgoth{L}_k^{(\alpha-1)}(s) \nonumber \\
\label{fourier2}
&= \alpha^{2} \beta^{\alpha/2} \sum_{k=0}^{\infty}\Big( \frac{(\alpha)_k}{k!} \Big)^{1/2} \rho_k (b_{\alpha}^2-k\beta) \textgoth{L}_k^{(\alpha-1)} (s),
\end{align}

Let
\begin{align}
\label{c1}
c_1 &= \int_0^\infty e^{-t/\alpha} \phi(t) \, \dd P_0(t) \nonumber \\
&= \beta^{\alpha/2} \sum_{k=0}^\infty \langle \phi, \textgoth{L}_k^{(\alpha-1)} \rangle_{L^2}\Big( \frac{(\alpha)_k}{k!} \Big)^{1/2} \rho_k,
\end{align}
and
\begin{align}
\label{c2}
c_2 &= \int^\infty_0{te^{-t/\alpha}\phi(t) }\, \dd P_0(t) \nonumber \\
&= \alpha^{2}\beta^{\alpha/2} \sum_{k=0}^{\infty} \langle \phi, \textgoth{L}_k^{(\alpha-1)} \rangle_{L^2}\Big( \frac{(\alpha)_k}{k!} \Big)^{1/2}  \rho_k (b_{\alpha}^2-k\beta).
\end{align}
Combining (\ref{equationeigensummary})-(\ref{c2}), we find that (\ref{eigen}) reduces to 
\begin{multline}
\label{equationeigen2}
\delta \sum_{k=0}^{\infty} \langle \phi, \textgoth{L}_k^{(\alpha-1)} \rangle_{L^2} \textgoth{L}_k^{(\alpha-1)}(s) \\
= \sum_{k=0}^{\infty} \rho_k \Big[\langle \phi, \textgoth{L}_k^{(\alpha-1)} \rangle_{L^2} -  \beta^{\alpha/2}\Big( \frac{(\alpha)_k}{k!} \Big)^{1/2} \big( c_1 + c_2 \alpha^{-1} (b_{\alpha}^2-k\beta)\big) \Big] \textgoth{L}_k^{(\alpha-1)}(s),
\end{multline}
and now comparing the coefficients of $\textgoth{L}_k^{(\alpha-1)} (s)$, we obtain 
\begin{equation}
\label{equationlemma}
\delta \ \langle \phi, \textgoth{L}_k^{(\alpha-1)} \rangle_{L^2} = \rho_k \Big[\langle \phi, \textgoth{L}_k^{(\alpha-1)} \rangle_{L^2} - \beta^{\alpha/2}\Big( \frac{(\alpha)_k}{k!} \Big)^{1/2} \big( c_1+ c_2 \alpha^{-1} (b_{\alpha}^2-k\beta)\big) \Big],
\end{equation}
for all $k \in \bN_0$.  Since we have assumed that $\delta \neq \rho_k$ for any $k$ then we can solve this equation for $\langle \phi, \textgoth{L}_k^{(\alpha-1)} \rangle_{L^2}$ to obtain
\begin{equation}
\label{innerproduct_eigen}
\langle \phi, \textgoth{L}_k^{(\alpha-1)} \rangle_{L^2} = \beta^{\alpha/2} \frac{\rho_k }{\rho_k-\delta} \Big( \frac{(\alpha)_k}{k!} \Big)^{1/2} \big(c_1+c_2 \alpha^{-1}(b_\alpha^2 - k\beta) \big).
\end{equation}
Substituting (\ref{innerproduct_eigen}) into (\ref{c1}), we get
\begin{align*}
 c_1 &= c_1 \beta^{\alpha}\sum_{k=0}^{\infty} \frac{(\alpha)_k}{k! (\rho_k-\delta)}\rho_k^2 + c_2 \alpha^{-1}  \beta^{\alpha}\sum_{k=0}^{\infty} \frac{(\alpha)_k}{k! (\rho_k-\delta)}\rho_k^2(b_{\alpha}^2-k\beta) \\
&= c_1 \big(1-A(\delta)\big) + c_2 \alpha^{-3} D(\delta);
\end{align*}
therefore, 
\begin{equation}
\label{firstc1c2eq}
\alpha^{3} c_1 A(\delta) = c_2 D(\delta).
\end{equation}
Similarly, by substituting (\ref{innerproduct_eigen}) into (\ref{c2}), we obtain
\begin{align*}
c_2 &= c_1 \alpha^2 \beta^{\alpha}\sum_{k=0}^{\infty} \frac{(\alpha)_k}{k! (\rho_k-\delta)}\rho_k^2(b_{\alpha}^2-k\beta) + c_2 \alpha \beta^{\alpha} \sum_{k=0}^{\infty} \frac{(\alpha)_k}{k! (\rho_k-\delta)}\rho_k^2 (b_{\alpha}^2-k\beta)^2 \\
&= c_1 D(\delta) + c_2 \big(1-B(\delta)\big);
\end{align*}
hence, 
\begin{equation}
\label{secondc1c2eq}
c_2B(\delta) = c_1D(\delta).
\end{equation}

Suppose $c_1=c_2=0$; then it follows from (\ref{innerproduct_eigen}) that $\langle \phi, \textgoth{L}_k^{(\alpha-1)} \rangle_{L^2}=0$ for all $k$, which implies $\phi=0$, which is a contradiction since $\phi$ is a non-trivial eigenfunction.  Hence, $c_1$ and $c_2$ cannot be both equal to $0$.  

Combining (\ref{firstc1c2eq}) and (\ref{secondc1c2eq}), and using the fact that $c_1, c_2$ are not both $0$, it is straightforward to deduce that $\alpha^3 A(\delta) B(\delta)= D^2(\delta)$.  Therefore, if $\delta$ is a positive eigenvalue of $\mathcal{S}$ then it is a positive root of the function 
$G(\delta) = \alpha^3 A(\delta) B(\delta)-D^2(\delta)$.  

Conversely, suppose that $\delta$ is a positive root of $G(\delta)$ with $\delta \neq \rho_k$ for any $k \in \mathbb{N}_0$.  Define
\begin{equation}
\label{gammak}
\gamma_k := \beta^{\alpha/2} \Big(\frac{(\alpha)_k}{k!}\Big)^{1/2} \frac{\rho_k }{\rho_k-\delta} \big(c_1+c_2 \alpha^{-1}(b_{\alpha}^2-k\beta)\big),
\end{equation}
$k \in \mathbb{N}_0$, where $c_1$ and $c_2$ are real constants that are not both equal to $0$ and which satisfy (\ref{firstc1c2eq}) and (\ref{secondc1c2eq}).  That such constants exist can be shown by following a case-by-case argument similar to \cite[p. 48]{ref27}; for example, if $D(\delta) \neq 0$, $A(\delta) \neq 0$, and $B(\delta) \neq 0$, then we can choose $c_2$ to be any non-zero number and then set $c_1 = c_2 B(\delta)/D(\delta)$.  

Now define
\begin{align}
\label{eigenfunctionfull}
\widetilde\phi (s) :=\sum_{k=0}^{\infty} \gamma_k& \textgoth{L}_k^{(\alpha-1)}(s),
\end{align}
$s \ge 0$.  By applying the ratio test, we find that $\sum_{k=0}^{\infty} \gamma_k^2 < \infty$; therefore, $\widetilde\phi \in L^2$. 

We also verify that the series (\ref{eigenfunctionfull}) converges pointwise.  By (\ref{laguerre2}) and (\ref{deflen}), 
$$
\big|\textgoth{L}_k^{(\alpha-1)}(s)\big| = \beta^{\alpha/2}\exp((1-\beta)s/2) \Big( \frac{k!}{(\alpha)_k} \Big)^{1/2} \ |L_k^{(\alpha-1)}(\beta s)|, 
$$
$s \ge 0$.  By Erd\'elyi, {\it et al.} \cite[p. 207]{ref5}, 
$$
\big| L_k^{(\alpha-1)}(\beta s) \big| \le 
\begin{cases}
\displaystyle{\Big( 2-\frac{(\alpha)_k}{k!} \Big) e^{\beta s/2}}, & 1/2 \le \alpha < 1 \\
\displaystyle{\frac{(\alpha)_k}{k!} e^{\beta s/2}}, \phantom{\Bigg]} & \alpha \ge 1
\end{cases}
$$
$s \ge 0$.  Therefore, 
\begin{align}
\label{ineqeigen2}
\big|\textgoth{L}_k^{(\alpha-1)}(s)\big| \le 
\begin{cases}
\displaystyle{\beta^{\alpha/2} \Big[ 2 \Big( \frac{k!}{(\alpha)_k} \Big)^{1/2} -\Big( \frac{(\alpha)_k}{k!}\Big)^{1/2} \Big] e^{s/2}}, & 1/2 \le \alpha < 1 \\
\displaystyle{\beta^{\alpha/2} \Big( \frac{(\alpha)_k}{k!} \Big)^{1/2} e^{s/2}}, \phantom{\Bigg]} & \alpha \ge 1
\end{cases}
\end{align}
Thus, to establish the pointwise convergence of the series (\ref{eigenfunctionfull}), we need to show that 
\begin{equation}
\label{conditionsptwiseconv}
\sum_{k=0}^{\infty}\Big( \frac{(\alpha)_k}{k!} \Big)^{1/2} |\gamma_k| < \infty \quad\text{and}\quad \sum_{k=0}^{\infty} \Big( \frac{k!}{(\alpha)_k} \Big)^{1/2} |\gamma_k| < \infty.
\end{equation}
However, the convergence of each of these series follows from the ratio test. 

Next, we justify the interchange of summation and integration in our calculations.  By a corollary to Theorem 16.7 in Billingsley \cite[p. 224]{billi2}, we need to verify that 
\begin{equation}
\label{conditionsptwiseconv2}
\sum_{k=0}^{\infty} |\gamma_k| \int^\infty_0 { K(s,t) \ | \textgoth{L}_k^{(\alpha-1)}(t) |}\, \dd P_0(t) < \infty.
\end{equation}
First, we find a bound for $K_0(s,t)$. By the inequality (\ref{besselineq3}) for the modified Bessel function, we have 
$$
\Gamma(\alpha)(st/\alpha^2)^{(1-\alpha)/2}I_{\alpha-1}\big(2\sqrt{st}/\alpha\big) \le \exp(2\sqrt{st}/\alpha);
$$
$s, t \ge 0$.  Therefore, by (\ref{K0kernel}), 
\begin{align}
\label{inequalitycovreduced}
0 \le K_0(s,t) &\le \exp(-(s+t)/\alpha) \exp(2\sqrt{st}/\alpha) \nonumber \\
&= \exp(-(\sqrt{s} - \sqrt{t})^2/\alpha) \le 1.
\end{align}
By the triangle inequality and by (\ref{inequalitycovreduced}), we have 
\begin{align*}
0 \le K(s,t)  & \le K_0(s,t)+ e^{-(s+t)/\alpha}(\alpha^{-3}st+ 1)\\ 
& \le 2 + \alpha^{-3}st,
\end{align*}
$s, t \ge 0$.  Thus, to prove (\ref{conditionsptwiseconv2}), we need to establish that 
$$
\sum_{k=0}^{\infty} |\gamma_k| \ \int^\infty_0 { (2 + \alpha^{-3}st) \ | \textgoth{L}_k^{(\alpha-1)}(t) | }\,\dd P_0(t) < \infty.
$$
By applying the bound (\ref{ineqeigen2}), we see that it suffices to prove that
$$
\sum_{k=0}^{\infty}\Big( \frac{(\alpha)_k}{k!} \Big)^{1/2} |\gamma_k | \ \int^\infty_0 t^j \, \dd P_0(t) < \infty
$$
and
$$
\sum_{k=0}^{\infty} \Big( \frac{k!}{(\alpha)_k} \Big)^{1/2} |\gamma_k | \ \int^\infty_0 t^j \, \dd P_0(t) < \infty,
$$
$j=0,1$.  As these integrals are finite, the convergence of both series follows from (\ref{conditionsptwiseconv}).  

To calculate $\mathcal{S} \widetilde\phi (s)$ from (\ref{eigenfunctionfull}), we follow the same steps as before to obtain 
\begin{align*}
\mathcal{S} \widetilde\phi (s)&=\int^\infty_0 K(s,t) \sum_{k=0}^{\infty} \gamma_k \textgoth{L}_k^{(\alpha-1)}(t) \, \dd P_0(t) \nonumber\\
&= \sum_{k=0}^{\infty} \rho_k \gamma_k \textgoth{L}_k^{(\alpha-1)}(s) -c_1\beta^{\alpha/2} \sum_{k=0}^{\infty}\Big( \frac{(\alpha)_k}{k!} \Big)^{1/2}\rho_k \textgoth{L}_k^{(\alpha-1)}(s)\nonumber\\
& {\hskip 1.5truein} -c_2 \alpha^{-1} \beta^{\alpha/2}  \sum_{k=0}^{\infty}\Big( \frac{(\alpha)_k}{k!} \Big)^{1/2} \rho_k (b_{\alpha}^2-k\beta) \textgoth{L}_k^{(\alpha-1)} (s).
\end{align*}
By the definition (\ref{gammak}) of $\gamma_k$, and noting that 
$$
\frac{\rho_k}{\rho_k-\delta} -1=\frac{\delta}{\rho_k-\delta},
$$
we have  
\begin{align*}
\mathcal{S} \widetilde\phi(s) &= \beta^{\alpha/2} \sum_{k=0}^{\infty} \Big[\frac{\rho_k}{\rho_k-\delta} -1\Big]\Big( \frac{(\alpha)_k}{k!} \Big)^{1/2} \rho_k (c_1+c_2 \alpha^{-1}(b_{\alpha}^2-k\beta)) \textgoth{L}_k^{(\alpha-1)} (s) \\
&= \beta^{\alpha/2} \delta \sum_{k=0}^{\infty} \frac{\rho_k}{\rho_k-\delta} \Big( \frac{(\alpha)_k}{k!} \Big)^{1/2} (c_1+c_2 \alpha^{-1}(b_{\alpha}^2-k\beta)) \textgoth{L}_k^{(\alpha-1)}(s)\\
&= \delta \sum_{k=0}^{\infty} \gamma_k \textgoth{L}_k^{(\alpha-1)}(s)\\
&= \delta \widetilde\phi(s).
\end{align*}
Therefore, $\delta$ is an eigenvalue of $\mathcal{S}$ with corresponding eigenfunction $\widetilde\phi$. 
$\qed$

\medskip

In the previous result, we assumed that $\delta \notin \{\rho_k: k \in \bN_0\}$.  
As stated in the following conjecture, we believe that this assumption is valid for all $\alpha$.  

\begin{conjecture}
\label{lemmaeigen}
Let $\delta$ be an eigenvalue of the covariance operator $\mathcal{S}$.  Then, there is no $l \in \bN_0$ such that $\delta = \rho_l$.
\end{conjecture}

This conjecture is equivalent to: 

\begin{conjecture} 
\label{lambertconjecture}
There is no $l \in \bN_0$ such that  
\begin{equation}
\label{eq:conjecture}
\alpha \beta^{\alpha+2} \sum_{\substack{k=0 \\ k \neq l}}^{\infty} \frac{(\alpha)_k}{k!} \frac{\rho_k^2}{\rho_k-\rho_l} (l-k)^2 = 1+\alpha(b_{\alpha}^2-l\beta)^2.
\end{equation}
\end{conjecture} 

\noindent
{\it Proof of the equivalence of Conjectures \ref{lemmaeigen} and \ref{lambertconjecture}}: Suppose there exists $l \in \mathbb{N}_0$ such that $\delta = \rho_l$.  Substituting $k = l$ in (\ref{equationlemma}) and simplifying the outcome, we obtain 
\begin{equation}
\label{equationc1c2lemma}
c_1 = c_2 \alpha^{-1} (l\beta-b_{\alpha}^2).
\end{equation}
Substituting $\delta=\rho_l$ in (\ref{equationeigen2}), applying (\ref{equationc1c2lemma}), and then cancelling common terms in (\ref{equationeigen2}), we obtain 
\begin{equation}
\label{phiLk0}
\langle \phi, \textgoth{L}_k^{(\alpha-1)} \rangle_{L^2} = c_2 \, \alpha^{-1} \beta^{(2+\alpha)/2} \, \Big(\frac{(\alpha)_k}{k!}\Big)^{1/2} \frac{l - k}{\rho_k - \rho_l} \rho_k,
\end{equation}
for $k \neq l$.  Substituting this result for the inner product into (\ref{c1}), we obtain 
\begin{align*}
c_1 &= \beta^{\alpha/2} \Big[\sum_{\substack{k=0 \\ k \neq l}}^\infty \langle \phi,\textgoth{L}_k^{(\alpha-1)}\rangle_{L^2} \Big(\frac{(\alpha)_k}{k!}\Big)^{1/2} \rho_k + \langle \phi,\textgoth{L}_l^{(\alpha-1)}\rangle_{L^2} \Big(\frac{(\alpha)_l}{l!}\Big)^{1/2} \rho_l\Big] \nonumber \\
&= \beta^{\alpha/2} \Big[\sum_{\substack{k=0 \\ k \neq l}}^\infty c_2 \, \alpha^{-1} \beta^{(2+\alpha)/2} \, \frac{(\alpha)_k}{k!} \frac{l - k}{\rho_k - \rho_l} \rho_k^2 + \langle \phi,\textgoth{L}_l^{(\alpha-1)}\rangle_{L^2} \Big(\frac{(\alpha)_l}{l!}\Big)^{1/2} \rho_l\Big].
\end{align*}
Similarly, substituting (\ref{phiLk0}) into (\ref{c2}), we obtain 
\begin{align*}
c_2 &= \alpha^2 \beta^{\alpha/2} \Big[\sum_{\substack{k=0 \\ k \neq l}}^\infty \langle \phi,\textgoth{L}_k^{(\alpha-1)}\rangle_{L^2} \Big(\frac{(\alpha)_k}{k!}\Big)^{1/2} \rho_k (b_\alpha^2 - k\beta) \nonumber \\
& \qquad\qquad\qquad\qquad\qquad + \langle \phi,\textgoth{L}_l^{(\alpha-1)}\rangle_{L^2} \Big(\frac{(\alpha)_l}{l!}\Big)^{1/2} \rho_l (b_\alpha^2 - l\beta)\Big] \nonumber \\
&= \alpha^2 \beta^{\alpha/2} \Big[\sum_{\substack{k=0 \\ k \neq l}}^\infty c_2 \, \alpha^{-1} \beta^{(2+\alpha)/2} \, \frac{(\alpha)_k}{k!} \frac{l - k}{\rho_k - \rho_l} \rho_k^2 (b_\alpha^2 - k\beta) \nonumber \\
& \qquad\qquad\qquad\qquad\qquad +  \langle \phi,\textgoth{L}_l^{(\alpha-1)}\rangle_{L^2} \Big(\frac{(\alpha)_l}{l!}\Big)^{1/2} \rho_l (b_\alpha^2 - l\beta)\Big].
\end{align*}
On simplifying the above expressions and substituting $c_1$ from (\ref{equationc1c2lemma}), we obtain
\begin{align}
\label{conjecinner1}
\beta^{\alpha/2} \Big(\frac{(\alpha)_l}{l!}\Big)^{1/2} \rho_l & \langle \phi, \textgoth{L}_l^{(\alpha-1)} \rangle_{L^2} \nonumber \\
&= c_2 \Big[\alpha^{-1} (l\beta-b_{\alpha}^2)-\alpha^{-1} \beta^{\alpha+1} \sum_{\substack{k=0 \\ k \neq l}}^\infty \frac{(\alpha)_k}{k!} \frac{l - k}{\rho_k - \rho_l} \rho_k^2 \Big],
\end{align}
and
\begin{align}
\label{conjecinner2}
\alpha^2 \beta^{\alpha/2} \Big(\frac{(\alpha)_l}{l!}\Big)^{1/2} \rho_l & (b_{\alpha}^2-l\beta) \langle \phi, \textgoth{L}_l^{(\alpha-1)} \rangle_{L^2}  \nonumber\\
& = c_2 \Big[1-\alpha \beta^{\alpha+1} \sum_{\substack{k=0 \\ k \neq l}}^\infty  \frac{(\alpha)_k}{k!} \frac{l - k}{\rho_k - \rho_l} \rho_k^2 (b_\alpha^2 - k\beta) \Big].
\end{align}
Suppose that $c_2 = 0$ then it follows from (\ref{equationc1c2lemma}) that $c_1 = 0$, which contradicts the earlier observation that $c_1$ and $c_2$ are not both zero; therefore, $c_2 \neq 0$.  Also, by (\ref{betaandbalpha}), $b_\alpha^2 < 1 < \beta$, so $b_\alpha^2 - k\beta \neq 0$ for all $k \in \bN_0$.  Solving (\ref{conjecinner1}) and (\ref{conjecinner2}) for the inner product $\langle \phi, \textgoth{L}_l^{(\alpha-1)} \rangle_{L^2}$ and equating the two expressions, we obtain
\begin{align*}
1-\alpha \beta^{\alpha+1}  \sum_{\substack{k=0 \\ k \neq l}}^\infty & \frac{(\alpha)_k}{k!} \frac{l - k}{\rho_k - \rho_l} \rho_k^2 (b_\alpha^2 - k\beta)  \\
& = \alpha (b_{\alpha}^2-l\beta)\Big[ (l\beta-b_{\alpha}^2)-\beta^{\alpha+1} \sum_{\substack{k=0 \\ k \neq l}}^\infty  \frac{(\alpha)_k}{k!} \frac{l - k}{\rho_k - \rho_l} \rho_k^2 \Big].
\end{align*}
Simplifying the above equation, we obtain (\ref{eq:conjecture}).
$\qed$

\begin{remark}
{\rm
Since $b_\alpha < 1$ then $\rho_k < \rho_0$ for all $k \ge 1$.  Therefore, if $l = 0$ then each term in the sum on the left-hand side of (\ref{eq:conjecture}) is negative, hence the sum itself is negative.  On the other hand, the right-hand side clearly is positive.  Therefore, it follows that the conjecture is valid if $l = 0$.  

We remark that Taherizadeh \cite{ref27} proved Conjecture \ref{lambertconjecture} for the case in which $\alpha=1$.  Also, for $\alpha = 2$, the conjecture is proved by Hadjicosta \cite{hadjicosta19}.  A consequence of those proofs is that the left-hand side of (\ref{eq:conjecture}) is seen to be less than the right-hand side for $\alpha = 1,2$, so this leads us to conjecture further that, for all $l \in \bN$, 
$$
\alpha \beta^{\alpha+2} \sum_{\substack{k=0 \\ k \neq l}}^{\infty} \frac{(\alpha)_k}{k!} \frac{\rho_k^2}{\rho_k-\rho_l} (l-k)^2 > 1+\alpha(b_{\alpha}^2-l\beta)^2.
$$
}
\end{remark}

\smallskip

\subsection{An interlacing property of the eigenvalues}
\label{decayrate}

A difficulty of the eigenvalues $\delta_k$ is that they have no closed form expression; hence there is no simple formula for $m$, the number of terms in the truncated series $\sum_{k=1}^m \delta_k \chi^2_{1k}$ that should be used in practice to approximate the asymptotic distribution, $\sum_{k=1}^\infty \delta_k \chi^2_{1k}$, of the test statistic $T_n^2$.  

For $\alpha =1$, Baringhaus and Taherizadeh \cite{BT2010} calculated several $\delta_k$ numerically and found that the truncated sum $\sum_{k=1}^{10} \delta_k$ closely approximates the exact value of $Tr(S)$; hence, the distribution of the truncated sum, $\sum_{k=1}^{10} \delta_k \chi^2_{1k}$ is a good approximation to the asymptotic distribution, $\sum_{k=1}^\infty \delta_k \chi^2_{1k}$, of $T_n^2$.  This approach is feasible since, as $\mathcal{S}_0$ is of trace-class then by \cite[p. 237, Corollary 3.2]{ref6}, $Tr(\mathcal{S}_0)$ can be calculated either by integrating the kernel $K_0$ or by evaluating the sum of all eigenvalues $\rho_k$: 
\begin{equation}
\label{TraceS0}
\int_{0}^{\infty} { K_0(s,s)}\, \dd P_0(s) = Tr(\mathcal{S}_0) = \sum_{k=0}^{\infty} \rho_k=\alpha^{\alpha} b_{\alpha}^{2\alpha} (1-b_{\alpha}^4)^{-1}.
\end{equation}
Since $\mathcal{S}$ also is of trace-class then 
\begin{align}
\label{TraceS1}
\sum_{k=1}^{\infty} \delta_k = Tr(\mathcal{S}) &= \int_0^\infty K(s,s) \, \dd P_0(s) \nonumber \\
&= \int_0^\infty \big[K_0(s,s) - (\alpha^{-3} s^2 + 1) e^{-2s/\alpha}\big] \, \dd P_0(s) \nonumber \\
&= \alpha^{\alpha} \Big[\frac{b_{\alpha}^{2\alpha}}{1-b_{\alpha}^4} -\frac{1}{(\alpha+2)^{\alpha}} \Big(1 +\frac{(\alpha+1)}{(\alpha+2)^2} \Big) \Big].
\end{align}

To determine for general $\alpha$ the number of terms in the truncated series $\sum_{k=1}^m \delta_k \chi^2_{1k}$ that should be used in practice to approximate the asymptotic distribution of $T_n^2$, we will derive bounds for the eigenvalues $\delta_k$ in terms of the $\rho_k$ and then obtain a general formula for $m$ as a function of $\alpha$.  In this regard, we are reminded of the concept of a ``scree plot'' in principal component analysis \cite[p. 441]{Johnson}, and we refer to the ratio $(\sum_{k=1}^m \delta_k)/Tr(\mathcal{S})$ as the $m${\hskip1pt}{\it th scree ratio} for $T_n^2$.  

Since the operator $\mathcal{S}$ is compact and positive then the set of all its eigenvalues is countable and contains only nonnegative values \cite[Theorem 8.12, p. 98]{ref17}.   To prove that the eigenvalues are positive and also are simple, i.e., of multiplicity $1$, we will apply the theory of total positivity; see Karlin \cite{Karlin}.  In what follows, we denote by $\det(a_{ij})$ the $r \times r$ determinant with $(i,j)$th entry $a_{ij}$.  

Recall that a $C^\infty$ kernel $K:\bR^2 \to \bR$ is {\it extended totally positive} (ETP) if for all $r \ge 1$, all $s_1 \ge \cdots \ge s_r$, all $t_1 \ge \cdots \ge t_r$, there holds 
\begin{equation}
\label{ETP}
\frac{\det\big(K(s_i,t_j)\big)}{\prod_{1 \le i< j \le r} (s_i - s_j)(t_i - t_j)} > 0,
\end{equation}
where instances of equality for the variables $s_i$ and $t_j$ are to be understood as limiting cases, and then L'Hospital's rule is to be used to evaluate this ratio. 

\begin{proposition}
\label{eigenvaluessimple}
The eigenvalues $\{\delta_k: k \ge 1\}$ of $\mathcal{S}$ and the eigenvalues $\{\rho_k: k \ge 0\}$ of $\mathcal{S}_0$ are positive and simple.  In particular, $\mathcal{S}$ and $\mathcal{S}_0$ are injective.  Further, the corresponding eigenfunctions $\{\phi_k: k \ge 1\}$ of $\mathcal{S}$ satisfy the oscillation property, 
\begin{equation}
\label{oscillation}
(-1)^{r(r-1)/2} \det\big(\phi_i(s_j)\big) \ge 0
\end{equation}
for all $r \ge 1$ and $0 \le s_1 < \cdots < s_r < \infty$, and the same property holds for the eigenfunctions $\{\textgoth{L}_k: k \ge 0\}$ of $\mathcal{S}_0$.
\end{proposition}

\Pro
By (\ref{covariance}), the kernel $K(s,t)$ is of the form 
$$
K(s,t) = e^{-(s+t)/\alpha} s^2 t^2 \sum_{k=0}^\infty c_k s^k t^k,
$$
where the coefficients $c_k$ are positive for all $k = 0,1,2,\ldots$.  Therefore, 
\begin{align*}
\det\big(K(s_i,t_j)\big) &= \det\Big(e^{-(s_i+t_j)/\alpha} s_i^2 t_j^2 \sum_{k=0}^\infty c_k s_i^k t_j^k\Big) \\
&= \Big(\prod_{i=1}^r e^{-(s_i+t_i)/\alpha} s_i^2 t_i^2\Big) \cdot \det\Big(\sum_{k=0}^\infty c_k s_i^k t_j^k\Big).
\end{align*}
By Karlin \cite[p. 101]{Karlin} the series, $\sum_{k=0}^\infty c_k s^k t^k$, is ETP.  Therefore, it follows from (\ref{ETP}) that $K(s,t)$ is ETP.

In the case of $K_0$, we have 
$$
K_0(s,t) = e^{-(s+t)/\alpha} \sum_{k=0}^\infty c_k s^k t^k,
$$
where $c_k > 0$ for all $k$.  Then it follows by a similar argument that $K_0(s,t)$ is ETP.  

Finally, we apply some results of Karlin \cite{Karlin1964}.  According to those results, the eigenvalues of an integral operator are simple and positive if the kernel of the operator is ETP.  Therefore, the eigenvalues of $\mathcal{S}$  and  $\mathcal{S}_0$ are simple and positive.  In particular, $0$ is not an eigenvalue of $\mathcal{S}$ or $\mathcal{S}_0$; therefore, both operators are injective.  Also, the oscillation property (\ref{oscillation}) follows from \cite[Theorem 3]{Karlin1964}.  
$\qed$

\bigskip

We now derive an interlacing property of the eigenvalues $\delta_k$ and $\rho_k$.   

\begin{proposition}
\label{interlacing}
For $k \ge 1$, $\rho_{k-1} \ge \delta_k \ge \rho_{k+1}$.  In particular, $\delta_k = O(\rho_k)$ as $k \to \infty$.  
\end{proposition}

\Pro
Define the kernels $k_0(s,t) = - e^{-(s+t)/\alpha}$ and $k_1(s,t) = - e^{-(s+t)/\alpha} \alpha^{-3} st$, where $s, t \ge 0$.  Also, define on $L^2$ the corresponding integral operators, 
$$
\mathcal{U}_j f(s) = \int_0^\infty k_j(s,t) f(t) \dd P_0(t),
$$
$j=0,1$, $s \ge 0$.  Then it follows from (\ref{covariance}) that $\mathcal{S} = \mathcal{S}_0 + \mathcal{U}_0 + \mathcal{U}_1$.  

It is simple to verify that each $\mathcal{U}_j$ is self-adjoint and of rank one, i.e., the range of $\mathcal{U}_j$ is a one-dimensional subspace of $L^2$.  Also, $\mathcal{S}_0 + \mathcal{U}_0$ is compact and self-adjoint, and its corresponding kernel $K_0 + k_0$ is of the form 
$$
K_0(s,t) + k_0(s,t) = e^{-(s+t)/\alpha} st \sum_{j=0}^\infty c_j s^j t^j,
$$
where $c_j > 0$ for all $j$.  By the same argument as in the proof of Proposition \ref{eigenvaluessimple}, we find that the eigenvalues of $\mathcal{S}_0 + \mathcal{U}_0$ are simple and positive; hence, $\mathcal{S}_0 + \mathcal{U}_0$ is injective.  

Denote by $\omega_k$, $k  \ge 1$, the eigenvalues of $\mathcal{S}_0 + \mathcal{U}_0$, where $\omega_1 > \omega_2 > \cdots$.    Since $\mathcal{S}_0$ is compact, self-adjoint, and injective, and since $\mathcal{U}_0$ is self-adjoint and of rank one, it follows from Hochstadt \cite{hochstadt} or Dancis and Davis \cite{dancisdavis} that the eigenvalues of $\mathcal{S}_0$ interlace the eigenvalues of $\mathcal{S}_0 + \mathcal{U}_0$, i.e., $\rho_{k-1} \ge \omega_k \ge \rho_k$ for all $k \ge 1$.  Moreover, there is precisely one eigenvalue of $\mathcal{S}_0 + \mathcal{U}_0$ in one of the intervals $[\rho_k,\rho_{k-1})$, $(\rho_k,\rho_{k-1})$, or $(\rho_k,\rho_{k-1}]$.  

Since $\mathcal{U}_1$ is self-adjoint and of rank one then by a second application of Hochstadt's theorem, we find that the eigenvalues of $\mathcal{S}_0 + \mathcal{U}_0$ interlace the eigenvalues of $\mathcal{S}_0 + \mathcal{U}_0 + \mathcal{U}_1 \equiv \mathcal{S}$, i.e, $\omega_k \ge \delta_k \ge \omega_{k+1}$ for all $k \ge 1$.  Also, there is precisely one eigenvalue of $\mathcal{S}$ in one of the intervals $[\omega_{k+1},\omega_k)$, $(\omega_{k+1},\omega_k)$, or $(\omega_{k+1},\omega_k]$.  

Combining the conclusions of the preceding paragraphs, we deduce that $\rho_{k-1} \ge \delta_k \ge \rho_{k+1}$ for all $k \ge 1$.

Finally, since $\rho_k = \alpha^\alpha b_\alpha^{4k+2\alpha}$, it follows from the interlacing inequalities that $\delta_k = O(b_\alpha^{4k})$, hence $\delta_k = O(\rho_k)$.  
$\qed$

\bigskip

\begin{remark}{\rm 
\label{interlacingremark}
The preceding result yields the inequalities $\rho_0 \ge \delta_1 \ge \rho_2 \ge \delta_3 \ge \cdots$ and $\rho_1 \ge \delta_2 \ge \rho_3 \ge \delta_4 \ge \cdots$.  For the case in which $\alpha = 1$, we have observed from the tables of eigenvalues computed by Taherizadeh \cite[pp. 28, 54]{ref27} that the eigenvalues $\rho_k$ and $\delta_k$ satisfy the stronger, strict interlacing property: $\rho_k > \delta_k > \rho_{k+1}$ for all $k \ge 1$.  This leads us to conjecture that the strict interlacing property holds for general $\alpha$.  

There is also the issue of choosing the value of $m$ so that the $m$th scree ratio of $T_n^2$ exceeds $1-\epsilon$, where $0 < \epsilon< 1$.  Applying the interlacing inequalities for $\delta_k$, we obtain 
$\sum_{k=1}^m \delta_k \ge \sum_{k=2}^{m+1} \rho_k$.  Since $Tr(\mathcal{S}_0) > Tr(\mathcal{S})$, we advise that $m$ be chosen so that 
$$
\sum_{k=0}^{m+1} \rho_k \ge (1-\epsilon) Tr(\mathcal{S}_0).
$$
This criterion leads to a value for $m$ that is readily applicable in the analysis of data.  Substituting $\rho_k = \alpha^\alpha b_\alpha^{4k+2\alpha}$, evaluating the resulting geometric series in closed form, and substituting for $Tr(\mathcal{S}_0)$ from (\ref{TraceS0}), we obtain 
$$
\alpha^\alpha b_\alpha^{2\alpha} \frac{1 - b_\alpha^{4(m+2)}}{1 - b_\alpha^4} = \sum_{k=0}^{m+1} \rho_k \ge (1-\epsilon) Tr(\mathcal{S}_0) = (1-\epsilon) \alpha^\alpha b_\alpha^{2\alpha} \frac{1}{1 - b_\alpha^4}.
$$
Solving this inequality for $m$, we obtain 
\begin{equation}
\label{noeigenvalues}
m \ge \frac{\log\epsilon}{4\log b_\alpha} - 2.
\end{equation}
We illustrate this bound by calculating it for various values of $\alpha$.  For $\epsilon = 10^{-10}$, which represents accuracy to ten decimal places, this results in the following table:
}\end{remark}

\smallskip

\begin{table}[h!]
\caption{Values of the lower bound on $m$ for the scree ratio of $T_n^2$.}
\medskip
\label{screeratiotable}
\centering
\begin{tabular}{|r|rrrrrrrr|}
  \hline
  $\alpha$ & 0.5 & 0.75 &  1 & 3 & 5 & 10 & 20 & 50 \\
  $m$      &  15 &   12 & 10 & 6 & 4 &  3 &  2 &  1 \\
  \hline
  \end{tabular}
\end{table}

\subsection{Two applications to data}

We will now apply our test to two data sets.  

\begin{table}[h!]
\caption{Waiting times for a Geiger counter to observe 100 alpha particles.}
\label{dataset1}
\centering
   \begin{tabular}{cccccccccc}
    6.9 & 5.9 & 7.2 & 7.6 & 7.5 & 7.3 & 7.0 & 7.1 & 6.7 & 5.3\\
    6.7 & 7.1 & 6.1 & 6.3 & 5.4 & 6.4 & 6.5 & 7.3 & 5.7 & 7.4\\
    6.3 & 7.6 & 7.6 & 6.7 & 6.9
   \end{tabular}
\end{table}

The first data set, given above, is obtained from \cite[p.~155]{hogg} and provides $n=25$ waiting times (in seconds) for a Geiger counter to observe 100 alpha-particles emitted by barium-133.  
A Kolmogorov-Smirnov test of $H_0$, the null hypothesis that the data were drawn from a $Gamma(\alpha=100,\lambda=14.7)$ distribution, failed to reject $H_0$ at the 10\% level of significance \cite[p.~464]{hogg}.  

We will apply the statistic $T_n^2$ to test the null hypothesis that the data are drawn from a gamma distribution with $\alpha=100$ and unspecified $\lambda$.  The observed value of the test statistic $T^2_n$, calculated using the data in Table \ref{dataset1}, is $6.301 \times 10^{-10}$.  

We used the limiting null distribution of $T_n^2$ to estimate $T^2_{n\,;\,0.05}$.  For $\alpha=100$, it follows from Table \ref{screeratiotable} that only one eigenvalue is needed to approximate accurately the asymptotic distribution of $T_n^2$; therefore, $T_n^2 \approx \delta_1 \chi^2_1$.  By (\ref{TraceS1}), we obtain $\delta_1 \simeq Tr(S) = 6.721718 \times 10^{-6}$.  Therefore, $T^2_{n\,;\,0.05} \simeq \delta_1 \chi^2_{1\,;\,0.05}$, where $\chi^2_{1\,;\,0.05}$ is the 95th percentile of the $\chi^2_1$ distribution, so we obtain $T^2_{n\,;\,0.05} = 2.582 \times 10^{-5}$.  As this critical value exceeds the observed value of $T_n^2$, we fail to reject the null hypothesis that the waiting times are drawn from a $Gamma(\alpha = 100,\lambda)$ distribution. 

For an alternative approach, we conducted a simulation study to approximate $T^2_{n\,;\,0.05}$\,, the 95th-percentile of the null distribution of $T_n^2$.  We generated $10,000$ random samples of size $n=25$ from the $Gamma(100,1)$ distribution, calculated the value of $T^2_n$ for each sample, and recorded the 95th-percentile of all 10,000 simulated values of $T_n^2$.   We repeated this process a total of ten times, finally approximating $T^2_{n\,;\,0.05}$ as the 20\%-trimmed mean of all 10 simulated 95th-percentiles, viz., $T^2_{n\,;\,0.05} = 2.368 \times 10^{-5}$.  Since this critical value exceeds the observed value of $T_n^2$ then we fail to reject the null hypothesis at the 5\% level of significance.  Moreover, we derived from our simulation study an approximate P-value of $0.99$ for the test.  

\medskip

\begin{table}[h!]
  \caption{Failure times for right rear brakes on some tractors.}
   \medskip
 \label{dataset2}
   \centering
  \begin{tabular}{cccccccccc}
56 & 83 & 104 & 116 & 244 & 305 & 429 & 452 & 453 & 503 \\
552 & 614 & 661 & 673 & 683 & 685 & 753 & 763 & 806 & 834 \\
838 & 862 & 897 & 904 & 981 & 1007 & 1008 & 1049 & 1060 & 1107 \\
1125 & 1141 & 1153 & 1154 & 1193 & 1201 & 1253 & 1313 & 1329 & 1347\\
1454 & 1464 & 1490 & 1491 & 1532 & 1549 & 1568 & 1574 & 1586 & 1599 \\
1608 & 1723 & 1769 & 1795 & 1927 & 1957 & 2005 & 2010 & 2016 & 2022\\ 2037 & 2065 & 2096 & 2139 & 2150 & 2156 & 2160 & 2190 & 2210 & 2220 \\
2248 & 2285 & 2325 & 2337 & 2351 & 2437 & 2454 & 2546 & 2565 & 2584\\
2624 & 2675 & 2701 & 2755 & 2877 & 2879 & 2922 & 2986 & 3092 & 3160 \\
3185 & 3191 & 3439 & 3617 & 3685 & 3756 & 3826 & 3995 & 4007 & 4159\\ 4300 & 4487 & 5074 & 5579 & 5623 & 6869 & 7739
   \end{tabular}
\end{table}

This second data set provides $n=107$ failure times (in hours) for right rear brakes on a collection of construction tractors. The data, reproduced from \cite{barlow}, were analyzed recently in \cite{cuparic}, where the null hypothesis of exponentiality was rejected.  

To test the hypothesis that these data are drawn from a gamma-distributed population, we will assume for illustrative purposes that $\alpha = 2.3$.  This value was obtained by setting the maximum likelihood estimate of the mode of the $Gamma(\alpha,\lambda)$ density, viz., $(\alpha-1)\bar{X}_n/\alpha$, equal to a mode of the histogram, and then solving the resulting equation for $\alpha$.  Consequently, the observed value of the test statistic $T^2_n$ is $0.0053$. 

For $\alpha = 2.3$, it follows from (\ref{noeigenvalues}) that $T_n^2 \approx \sum_{k=1}^7 \delta_k \chi^2_{1k}$.  This requires that we first calculate the $\delta_k$ numerically as the positive roots of the function $G(\delta)$ in Theorem \ref{thmeigenvaluesofS}, and then we would apply the results of Kotz, et al. \cite{KJB} to derive the distribution of $\sum_{k=1}^7 \delta_k \chi^2_{1k}$ and carry out the test.  We recommend in practice the one-term approximation \cite[Eqs.~(71),~(79)]{KJB},
$$
P\Big(\sum_{k=1}^m \delta_k \chi^2_{1k} \ge t\Big) \simeq P\big(\chi^2_m \ge 2t/(\delta_1 + \delta_m)\big)
$$
which is well-known to be accurate for other problems \cite{GR1983} and which leads to the explicit expression, $T^2_{n\,;\,0.05} \simeq \frac12 (\delta_1+\delta_m)\chi^2_{m\,;\,0.05}$, for an approximate critical value of $T_n^2$.  

As an alternative to calculating $\delta_1,\ldots,\delta_m$, we can apply the interlacing inequalities in Proposition \ref{interlacing} to obtain a stochastic lower bound, $\sum_{k=1}^m \delta_k \chi^2_{1k} \ge \sum_{k=2}^{m+1} \rho_k \chi^2_{1k}$.  If we carry out the test by using the lower bound $\sum_{k=2}^{m+1} \rho_k \chi^2_{1k}$ with its exact distribution or a one-term approximation obtained from Kotz, et al. \cite[loc.~cit.]{KJB}, we will obtain a conservative test of the null hypothesis, i.e., with a level of significance at most 5\%.

By performing a simulation study as for the previous data set, we obtained the approximation, $T^2_{n\,;\,0.05} = 0.0356$, which exceeds the observed value of the test statistic.  Therefore, we fail to reject the null hypothesis at the 5\% level of significance.  Moreover, the simulation study provided an approximate P-value of $0.47$.   

\begin{table}[h!]
\caption{The outcome of testing the tractor brakes data with different values of $\alpha$.}
\smallskip
\label{caterpillartable}
\centering
\begin{tabular}{|r|l|l|l|l|l|l|}
  \hline
  $\alpha$                 &  1.0          & 1.8       & 2.3       &  3         & 5          & 8 \\[2pt]
  \hline 
  Observed $T_n^2$ & 0.6965             & 0.0162 & 0.0053 & 0.0534 & 0.1977 & 0.3180 \\[2pt]
  $T^2_{n\,;\,0.05}$ & 0.1406 & 0.0559 & 0.0356 & 0.0228 & 0.0088 & 0.0037 \\[2pt]
  P-value                    & 0.0000 & 0.3013 & 0.4702 & 0.0026 & 0.0000  & 0.0000 \\
  \hline
  \end{tabular}
\end{table}

Since was assumed for illustrative purposes that $\alpha = 2.3$, we repeated the test for several values of $\alpha$, obtaining the results in Table \ref{caterpillartable}.  We note that the null hypothesis is rejected for the case in which $\alpha = 1$ where, under the null hypothesis, the data are exponentially distributed; hence, we deduce in this case the same conclusion as in \cite{cuparic}.


\subsection{Consistency of the test}

\begin{theorem}
\label{consistency}
Let $X_1,X_2,\dotsc$ be a sequence of non-negative, i.i.d., random variables with finite mean $\mu$.  Let $\gamma \in (0,1)$ denote the level of significance of the test and $c_{n,\gamma}$ be the $(1-\gamma)$-quantile of the test statistic $T^2_n$ under $H_0$.  If $X_1,X_2,\dotsc$ are not Gamma-distributed then  
$$
\lim_{n \rightarrow \infty} P(T^2_n > c_{n,\gamma})=1.
$$
\end{theorem}

\Pro
Recall that
\begin{align*}
n^{-1} T^2_n=\int^\infty_0{\Big[ \frac{1}{n}\sum_{j=1}^{n} \Gamma(\alpha) (tY_{j})^{(1-\alpha)/2} J_{\alpha-1} (2(tY_{j})^{1/2})-e^{-t/\alpha}\Big]^2} \, \dd P_0(t) .
\end{align*}
By subtracting and adding the term
$$
\frac{1}{n}\sum_{j=1}^n \Gamma(\alpha) (tX_{j}/\mu)^{(1-\alpha)/2} J_{\alpha-1} (2(tX_{j}/\mu)^{1/2})
$$
inside the squared term, and then expanding the integrand, we obtain
\begin{align}
n^{-1} T^2_n
\label{consistency1}
&=\int^\infty_0{\Big[ \frac{1}{n}\sum_{j=1}^{n} \Gamma(\alpha) (tX_{j}/\mu)^{(1-\alpha)/2} J_{\alpha-1} (2(tX_{j}/\mu)^{1/2})-e^{-t/\alpha}\Big]^2}\, \dd P_0(t) \\
& \quad + \int^\infty_0 {\Big[ \frac{1}{n}\sum_{j=1}^{n} \Gamma(\alpha) \big( (tY_{j})^{(1-\alpha)/2} J_{\alpha-1} (2(tY_j)^{1/2}) } \nonumber \\
\label{consistency2}
& \quad\quad\quad\quad\quad\quad\quad\quad\quad\quad\quad {-(tX_{j}/\mu)^{(1-\alpha)/2} J_{\alpha-1} (2(tX_{j}/\mu)^{1/2}) \big) \Big]^2} \, \dd P_0(t) \\
& \quad + 2\int^\infty_0{ \Big[ \frac{1}{n}\sum_{j=1}^{n} \Gamma(\alpha) (tX_{j}/\mu)^{(1-\alpha)/2} J_{\alpha-1} (2(tX_{j}/\mu)^{1/2})-e^{-t/\alpha}\Big]} \nonumber \\
& \quad\quad\quad\quad \times \Big[ \frac{1}{n}\sum_{j=1}^{n} \Gamma(\alpha) \big( (tY_{j})^{(1-\alpha)/2} J_{\alpha-1} (2(tY_j)^{1/2}) \nonumber \\
\label{consistency3}
& \quad\quad\quad\quad\quad\quad\quad\quad\quad\quad\quad {-(tX_{j}/\mu)^{(1-\alpha)/2} J_{\alpha-1} (2(tX_{j}/\mu)^{1/2}) \big) \Big]}\, \dd P_0(t).
\end{align}

Next, we proceed as in the first part of the proof of Theorem \ref{limitingnulldistribution}.  Applying the Taylor expansion (\ref{taylor1}) with $s = 2(tY_j)^{1/2}$ and $s_0 = 2(tX_j/\mu)^{1/2}$, we obtain 
\begin{equation}
\label{taylorexpconsistency}
\begin{aligned}
2^{1-\alpha}(tY_j)^{(1-\alpha)/2} \, J_{\alpha-1}\big(2(tY_{j})^{1/2}\big) &= 2^{1-\alpha} (tX_{j}/\mu)^{(1-\alpha)/2} J_{\alpha-1}\big(2(tX_{j}/\mu)^{1/2}\big) \\
& \qquad + 2 \big[(tX_{j}/\mu)^{1/2} - (tY_{j})^{1/2}\big] \, u_j^{1-\alpha} \, J_{\alpha}(u_j),
\end{aligned}
\end{equation}
where $u_j$ lies between $2(tY_{j})^{1/2}$and $2(tX_{j}/\mu)^{1/2}$.  Since 
\begin{align*}
(tY_{j})^{1/2}-(tX_{j}/\mu)^{1/2} &= (tX_{j}/\bar{X}_n)^{1/2}-(tX_{j}/\mu)^{1/2}\\
&= (tX_j)^{1/2}({\bar{X}_{n}}^{-1/2}-\mu^{-1/2}),
\end{align*}
then equation (\ref{taylorexpconsistency}), multiplied by $2^{\alpha-1}$, becomes
\begin{align}
(tY_j)^{(1-\alpha)/2}J_{\alpha-1}(2(tY_{j})^{1/2})&=(tX_{j}/\mu)^{(1-\alpha)/2} J_{\alpha-1}(2(tX_{j}/\mu)^{1/2})\nonumber\\
\label{taylorexpconsistency2}
& \quad\quad -2^{\alpha}(tX_j)^{1/2}({\bar{X}_{n}}^{-1/2}-\mu^{-1/2})u_j^{1-\alpha}J_{\alpha}(u_j).
\end{align}

First, we will establish that (\ref{consistency2}) converges to 0, almost surely. By Lemma \ref{lemma5},
$$
| u_j^{1-\alpha}J_{\alpha}(u_j) | \le \frac{1}{2^{\alpha-1}\pi^{1/2} \Gamma(\alpha+\frac{1}{2})},
$$
and therefore, by (\ref{taylorexpconsistency2}), 
\begin{multline*}
\frac{1}{n}\sum_{j=1}^n \Gamma(\alpha) |(tY_j)^{(1-\alpha)/2} J_{\alpha-1}(2(tY_{j})^{1/2})-(tX_{j}/\mu)^{(1-\alpha)/2} J_{\alpha-1} (2(tX_{j}/\mu)^{1/2}) | \\
\le  \frac{2 \Gamma(\alpha)}{\pi^{1/2} \Gamma(\alpha+\frac{1}{2})} \ t^{1/2} \ | {\bar{X}_{n}}^{-1/2}-\mu^{-1/2} | \ \frac{1}{n} \sum_{j=1}^n X_j^{1/2}.
\end{multline*}
By the triangle inequality, we conclude that the integral (\ref{consistency2}) is less than or equal to
$$
\Big( \frac{2 \Gamma(\alpha)}{\pi^{1/2} \Gamma(\alpha+\frac{1}{2})} \Big)^2 ({\bar{X}_{n}}^{-1/2}-\mu^{-1/2})^2 \big( \frac{1}{n}\sum_{j=1}^n X_j^{1/2}\big)^2 \int^{\infty}_0 {t}\, \dd P_0(t).
$$
By the Cauchy-Schwarz inequality,
$
\big(n^{-1}\sum_{j=1}^n X_j^{1/2}\big)^2 \le \bar{X}_n;
$ 
also $\int^{\infty}_0 {t}\, \dd P_0(t) = \alpha$.  Moreover, by the Strong Law of Large Numbers, $\bar{X}_n \to \mu$, almost surely.  Thus, the integral (\ref{consistency2}) converges to $0$ almost surely. 

Second, we show that (\ref{consistency3}) tends to $0$, almost surely.  Due to inequality (\ref{besselineq2}), the fact that $e^{-t/\alpha} \le 1$ for $t \ge 0$, and the triangle inequality, we have 
$$
\bigg| \frac{1}{n}\sum_{j=1}^{n} \Gamma(\alpha) (tX_{j}/\mu)^{(1-\alpha)/2} J_{\alpha-1} (2(tX_{j}/\mu)^{1/2})-e^{-t/\alpha} \bigg| \le 2.
$$
By the triangle inequality, the absolute value of (\ref{consistency3}) is less than or equal to
\begin{align}
\label{consistency3_ineq}
& 2 \int^{\infty}_0{\Big| \frac{1}{n}\sum_{j=1}^{n} \Gamma(\alpha) \big( (tY_{j})^{(1-\alpha)/2} J_{\alpha-1} (2(tY_j)^{1/2})}\nonumber\\
& \quad\quad\quad\quad\quad\quad\quad {-(tX_{j}/\mu)^{(1-\alpha)/2} J_{\alpha-1} (2(tX_{j}/\mu)^{1/2}) \big) \Big|}\, \dd P_0(t).
\end{align}
By the Cauchy-Schwarz inequality and the fact that $\int^\infty_0 \, \dd P_0(t) =1$, (\ref{consistency3_ineq}) is seen to be less than or equal to
\begin{align*}
& 2 \Big(\int^{\infty}_0{\Big[ \frac{1}{n}\sum_{j=1}^{n} \Gamma(\alpha) \big( (tY_{j})^{(1-\alpha)/2} J_{\alpha-1} (2(tY_j)^{1/2})}\nonumber\\
& \quad\quad\quad\quad\quad\quad\quad {-(tX_{j}/\mu)^{(1-\alpha)/2} J_{\alpha-1} (2(tX_{j}/\mu)^{1/2}) \big) \Big]^2}\, \dd P_0(t) \Big)^{1/2}.
\end{align*}
Following the same argument as for integral (\ref{consistency2}), we conclude that the integral (\ref{consistency3}) tends to $0$, almost surely. 

As for the integral (\ref{consistency1}), we subtract and add inside the squared term the Hankel transform of the random variable $X_1/\mu$, i.e., the term 
$$
E[\Gamma(\alpha)(tX_1/\mu)^{(1-\alpha)/2}J_{\alpha-1}(2(tX_1/\mu)^{1/2})],
$$
and expand the integrand.  Then we find that (\ref{consistency1}) equals 
\begin{align}
\int^\infty_0 & \Big[\frac{1}{n}\sum_{j=1}^{n} \Gamma(\alpha) (tX_{j}/\mu)^{(1-\alpha)/2} J_{\alpha-1} (2(tX_{j}/\mu)^{1/2}) \nonumber\\
\label{consistency11}
& {\hskip 1.2truein} - E[\Gamma(\alpha)(tX_1/\mu)^{(1-\alpha)/2}J_{\alpha-1}(2(tX_1/\mu)^{1/2})]\Big]^2 \, \dd P_0(t) \\
& + \int^\infty_0 \Big[ E[\Gamma(\alpha)(tX_1/\mu)^{(1-\alpha)/2}J_{\alpha-1}(2(tX_1/\mu)^{1/2})]-e^{-t/\alpha}\Big]^2 \, \dd P_0(t) \nonumber\\
\label{consistency13}
& + 2\int^\infty_0{\Big[ \frac{1}{n}\sum_{j=1}^{n} \Gamma(\alpha) (tX_{j}/\mu)^{(1-\alpha)/2} J_{\alpha-1} (2(tX_{j}/\mu)^{1/2})}\nonumber\\
& {\hskip 1.2truein} - E[\Gamma(\alpha)(tX_1/\mu)^{(1-\alpha)/2}J_{\alpha-1}(2(tX_1/\mu)^{1/2})] \Big] \nonumber\\
& {\hskip 0.5truein} \times {\Big[ E[\Gamma(\alpha)(tX_1/\mu)^{(1-\alpha)/2}J_{\alpha-1}(2(tX_1/\mu)^{1/2})]-e^{-t/\alpha} \Big]}\, \dd P_0(t). 
\end{align}
By the Strong Law of Large Numbers in $L^2$ \cite[p. 189, Corollary 7.10]{ref19}, we conclude that the term (\ref{consistency11}) converges to $0$, almost surely. 

Next, we show that (\ref{consistency13}) converges to $0$, almost surely.  By inequality (\ref{besselineq2}), and the fact that $e^{-t/\alpha} \le 1$ for $t \ge 0$, we have 
$$
\Big| E[\Gamma(\alpha)(tX_1/\mu)^{(1-\alpha)/2}J_{\alpha-1}(2(tX_1/\mu)^{1/2})]-e^{-t/\alpha} \Big| \le 2.
$$
Therefore, the absolute value of the integral (\ref{consistency13}) is less than or equal to
\begin{align*}
2 &\int^\infty_0 \Big| \frac{1}{n} \sum_{j=1}^{n} \Gamma(\alpha) (tX_j/\mu)^{(1-\alpha)/2}J_{\alpha-1}(2(tX_{j}/\mu)^{1/2}) \\
& {\hskip 1 truein} -E[\Gamma(\alpha)(tX_1/\mu)^{(1-\alpha)/2}J_{\alpha-1}(2(tX_1/\mu)^{1/2})] \Big|  \, \dd P_0(t) \\
& \le  2\Big( \int^\infty_0 \Big[ \frac{1}{n} \sum_{j=1}^{n} \Gamma(\alpha) (tX_j/\mu)^{(1-\alpha)/2}J_{\alpha-1}(2(tX_{j}/\mu)^{1/2}) \\
& {\hskip 1 truein} -E[\Gamma(\alpha)(tX_1/\mu)^{(1-\alpha)/2}J_{\alpha-1}(2(tX_1/\mu)^{1/2})] \Big]^2 \, \dd P_0(t) \Big)^{1/2},
\end{align*}
where the latter bound follows from the Cauchy-Schwarz inequality. Again, by the Strong Law of Large Numbers in $L^2$, 
we conclude that (\ref{consistency13}) converges to $0$, almost surely.  

We have now shown that 
\begin{eqnarray}
\label{almostsurelyconv}
n^{-1} T^2_n \xrightarrow{a.s.} \int^\infty_0{\Big( E[\Gamma(\alpha)(tX_1/\mu)^{(1-\alpha)/2}J_{\alpha-1}(2(tX_1/\mu)^{1/2})]-e^{-t/\alpha}\Big)^2}\, \dd P_0(t).
\end{eqnarray}
Denote by $\Lambda$ the right-hand side of (\ref{almostsurelyconv}); then, $\Lambda \ge 0$.  Suppose that $\Lambda = 0$, then 
$$
E[\Gamma(\alpha)(tX_1/\mu)^{(1-\alpha)/2}J_{\alpha-1}(2(tX_1/\mu)^{1/2})]-e^{-t/\alpha} = 0,
$$
equivalently, $\mathcal{H}_{X_1,\alpha-1}(t/\mu) - e^{-t/\alpha} = 0$, $P_0$-almost everywhere.  By continuity, we obtain $\mathcal{H}_{X_1,\alpha-1}(t/\mu) = e^{-t/\alpha}$ for all $t \ge 0$.  By the Uniqueness Theorem for Hankel transforms, it follows that $X_1$ has a gamma distribution, which contradicts the assumption that $X_1$ does not have a gamma distribution.  Therefore, $\Lambda > 0$.  

Under $H_0$, $n^{-1} T^2_n \xrightarrow{a.s.} 0$, and therefore $n^{-1} T^2_n\xrightarrow{p} 0$, i.e., for any $\epsilon >0$,
$$
\lim_{n \rightarrow \infty} P_{H_0} \big( n^{-1} T^2_n \ge \epsilon \big)=0.
$$
Thus, for any $\epsilon > 0$ and $\gamma>0$, there exists $n_0(\epsilon,\gamma) \in \mathbb{N}$ such that
$$
P_{H_0} \big( n^{-1} T^2_n \ge \epsilon \big) \le \gamma,
$$
for all $n \ge n_0(\epsilon,\gamma)$. Now, let $c_{n,\gamma}$ be the $(1-\gamma)$-quantile of the test statistic $T^2_n$ under $H_0$, where $\gamma$ is fixed; then
$
0 \le c_{n,\gamma} \le n\epsilon, 
$ 
for all $n \ge n_0 (\epsilon)$, since $c_{n,\gamma}:=\inf\{ x \ge 0: P_{H_0}(T^2_n >x) \le \gamma\}$.  Therefore, $0 \le n^{-1} c_{n,\gamma} \le \epsilon$ for all $n \ge n_0 (\epsilon)$.  In summary, for any $\epsilon > 0$, there exists $n_0(\epsilon) \in \bN$ such that
$n^{-1} c_{n,\gamma} \le \epsilon$ 
for all $n \ge n_0(\epsilon)$, i.e. 
\begin{eqnarray}
\label{quantileconv}
\lim_{n \rightarrow \infty} n^{-1} c_{n,\gamma} =0.
\end{eqnarray}

By (\ref{almostsurelyconv}) and (\ref{quantileconv}), $n^{-1} T^2_n-n^{-1} c_{n,\gamma} \xrightarrow{a.s.} \Lambda$, and hence $n^{-1} T^2_n-n^{-1}c_{n,\gamma} \xrightarrow{p} \Lambda$. Thus, by Severini \cite[ p. 340, Corollary 11.3 (i)]{ref12}), $n^{-1} T^2_n-n^{-1}c_{n,\gamma} \xrightarrow{d} \Lambda$.
Further, 
\begin{align*}
\lim_{n \rightarrow \infty} P(T^2_n > c_{n,\gamma})&=\lim_{n \rightarrow \infty} P\big( n^{-1}T^2_n -n^{-1} c_{n,\gamma} > 0 \big)\\
&=1-\lim_{n \rightarrow \infty} P\big( n^{-1}T^2_n -n^{-1}c_{n,\gamma} \le 0 \big).
\end{align*}
Since the distribution function of the constant positive random variable $\Lambda$ is continuous at 0, we conclude that 
$$
\lim_{n \rightarrow \infty} P(T^2_n > c_{n,\gamma})=1-P(\Lambda \le 0)=1-0=1.
$$
This concludes the proof.
$\qed$

\medskip

\begin{remark}{\rm
By applying Theorem 1 of Baringhaus, Ebner, and Henze \cite{henze17}, we also find that under fixed alternatives to the null hypothesis, 
$$
n^{-1/2}(n^{-1} T_n^2 - \Lambda) \xrightarrow{d} N(0,\sigma^2)
$$
as $n \to \infty$, where $\sigma^2$ is a constant that is evaluated from an integral involving the alternative distribution.
}\end{remark}

\section{Contiguous Alternatives to the Null Hypothesis}
\label{contiguous}

In this section, we find the limiting distribution of the test statistic under a sequence of contiguous alternatives. 

\subsection{Assumptions}
\label{assumptionscontiguous}

For $n \in \bN$, let $X_{n1},\dotsc,X_{nn}$ be a triangular array of row-wise independent random variables. As usual, let $P_0=Gamma(\alpha,1)$, $\alpha \ge 1/2$, and let $Q_{n1}$ be a probability measure dominated by $P_0$. 
We wish to test the hypothesis
$$
H_0: \ \text{The marginal distribution of each} \ X_{ni}, \ i=1,\dotsc,n, \ \text{is} \ P_0
$$
against the alternative
$$
H_1: \ \text{The marginal distribution of each} \ X_{ni}, \ i=1,\dotsc,n, \ \text{is} \ Q_{n1}.
$$
We write the Radon-Nikodym derivative of $Q_{n1}$ with respect to $P_0$ in the form
\begin{align}
\label{radon_nikodym}
\frac{\dd Q_{n1}}{\dd P_0}=1+n^{-1/2}h_n
\end{align}
We will need two assumptions in the sequel.

\begin{assumption}
\label{2assumptions}
{\rm
We assume that: 
\begin{itemize}
\item[(A1)] The functions $\{h_n: n \in \bN \}$ form a sequence of $P_0$-integrable functions converging pointwise, $P_0$-almost everywhere, to a function $h$, and 
\item[(A2)] $\sup_{n \in \bN} E_{P_0} |h_n|^4 < \infty$. 
\end{itemize}
}\end{assumption}

Note that since $ \int (\dd Q_{n1}/ \dd P_0) \, \dd P_0 = 1$ then we also have $\int {h_n}\, \dd P_0 = 0$, for all $n \in \bN$.  Also, by applying (A2), we deduce the uniform integrability of $|h_{n}|^2$: 
\begin{align*}
\lim_{k \rightarrow \infty} \sup_{n \in \bN} E_{P_0} \big( |h_n|^2 I( |h_n|^2 > k) \big)&=\lim_{k \rightarrow \infty} \sup_{n \in \bN} \int {|h_n|^2 \ I( |h_n|^2 >k)} \, \dd P_0\\
&\le \lim_{k \rightarrow \infty} \sup_{n \in \bN} \int {k^{-1} |h_n|^4}\, \dd P_0\\ 
&=\lim_{k \rightarrow \infty} k^{-1} \sup_{n \in \bN} E_{P_0} |h_n|^4=0.
\end{align*}
By Bauer \cite[p. 95, Theorem 2.11.4]{bauer}, the $P_0$-almost everywhere convergence of $h_n$ to $h$ implies the $P_0$-stochastic convergence of $h_n$ to $h$.  Again by Bauer \cite[p. 104, Theorem 2.12.4]{bauer}, the uniform integrability of $h_n^2$ together with the $P_0$-stochastic convergence of $h_n$ to $h$ imply the convergence of $h_n$ in mean square: 
$$
\lim_{n \rightarrow \infty} \int {|h_n-h|^2}\, dP_0=0,
$$
and therefore 
$$
\lim_{n \rightarrow \infty} \int |h_n|^2 \, dP_0 = \int |h|^2 \, dP_0.
$$
Since convergence in mean square implies convergence in mean, we have
$$
\lim_{n \rightarrow \infty} \int {|h_n-h|}\, dP_0=0,
$$
and thus, 
$$
\lim_{n \rightarrow \infty} \int {h_n}\, dP_0=\int {h}\, dP_0.
$$
Now, due to the fact that $\int {h_n}\, dP_0=0$ for all $n \in \bN$, we obtain 
$$
\lim_{n \rightarrow \infty} \int {h_n}\, dP_0=\int {h}\, dP_0=0.
$$

\subsection{Examples}
\label{subsec_examples}

In this subsection, we verify that Assumptions \ref{2assumptions} are valid for a broad collection of sequences of contiguous alternatives.

\subsubsection{Gamma alternatives with rate parameter \texorpdfstring{$\lambda_n=1+n^{-1/2}$}{lambda}}

Let $Q_{n1}:=Gamma(\alpha, \lambda_n)$ with $\alpha \ge 1/2$ and $\lambda_n=1+n^{-1/2}$.  Then, 
$$
\frac{\dd Q_{n1}}{\dd P_0}=\big( 1+n^{-1/2} \big)^{\alpha} \ \exp (-n^{-1/2} x ),
$$
$x \ge 0$. We equate this Radon-Nikodym derivative to $1+n^{-1/2}h_n(x)$, obtaining 
$$
h_n(x) =  n^{1/2} \big[\big( 1+n^{-1/2} \big)^{\alpha} \exp (-n^{-1/2} x )-1\big] ,
$$
for $x \ge 0$.  By applying L'Hospital's rule, we obtain 
$$
h(x):= \lim_{n \rightarrow \infty} h_n(x)=x+\alpha,
$$
for $x \ge 0$.  Next, we find $E_{P_0} |h_n|^4$. Define 
\begin{align}
\label{remainderterm}
R_n(x) &=\exp(-n^{-1/2} x) - (1-n^{-1/2} x) \\
& =\sum_{k=2}^{\infty} \frac{1}{k!} \big ( -n^{-1/2}x\big)^{k}, \nonumber
\end{align}
the remainder term of the Taylor series expansion of $\exp (-n^{-1/2} x)$, for $x \ge 0$. Then, by elementary algebraic manipulations, we obtain 
\begin{align*}
h_n(x) &= n^{1/2} (1+n^{-1/2})^{\alpha} \exp (-n^{-1/2} x ) - n^{1/2} \\
&= n^{1/2} (1+n^{-1/2})^{\alpha} \big(R_n(x) + 1-n^{-1/2} x\big) - n^{1/2} \\
&= (1+n^{-1/2})^{\alpha-1} \big[1+(1+n^{1/2})R_n(x) - (1+n^{-1/2})x \big] \\
& \hskip 2.5truein + n^{1/2}[(1+n^{-1/2})^{\alpha-1}  - 1].
\end{align*}
By (\ref{remainderterm}), the triangle inequality, and the Lipschitz continuity of the exponential function, we have
\begin{align*}
|R_n(x)| \le n^{1/2} |R_n(x) | &= n^{1/2} |\exp(-n^{-1/2} x) - (1-n^{-1/2} x)| \nonumber \\
&\le n^{1/2} \big[|\exp(-n^{-1/2} x) - 1| + n^{-1/2} x\big] \nonumber \\
&\le n^{1/2} \big[n^{-1/2} x + n^{-1/2} x\big] \nonumber \\
&= 2x,
\end{align*}
$x \ge 0$.  Therefore, 
\begin{align*}
|h_n(x)| &\le (1+n^{-1/2})^{\alpha-1} \Big[1+(1+n^{1/2})|R_n(x)| + (1+n^{-1/2})x \Big] \\
& \hskip 2.5truein + \big| n^{1/2}\big((1+n^{-1/2})^{\alpha-1}  - 1\big)\big| \\
&\le (1+n^{-1/2})^{\alpha-1} (1+4x + 2x ) + \big| n^{1/2}\big((1+n^{-1/2})^{\alpha-1}  - 1\big)\big| \\
&= (1+n^{-1/2})^{\alpha-1} (1+6x) + \big| n^{1/2}\big((1+n^{-1/2})^{\alpha-1}  - 1\big)\big|.
\end{align*}
It is elementary that $(1+n^{-1/2})^{\alpha-1} \to 1$ and that $n^{1/2}\big((1+n^{-1/2})^{\alpha-1} - 1\big) \to \alpha-1$ as $n \to \infty$; therefore, there exists a positive constant $M$ such that $(1+ n^{-1/2})^{\alpha-1} \le M$ and $\big|n^{1/2}\big((1+n^{-1/2})^{\alpha-1}  - 1\big) \big| \le M$ for all $n$.  Therefore, 
$|h_n(x)| \le  M(1 + 6x) + M = M(2+6x)$, $x \ge 0$, so we obtain 
\begin{align*}
E_{P_0}|h_n|^4 \le M^4 \int^{\infty}_0 (2 + 6x)^4 \, \dd P_0(x) < \infty,
\end{align*}
and this bound does not depend on $n$.  Thus, 
$\sup_{n \in \bN} E_{P_0} |h_n|^4 < \infty$.

\subsubsection{Gamma alternatives with shape parameter \texorpdfstring{$\alpha_n=\alpha+n^{-1/2}$}{alpha}}

Let $Q_{n1}:=Gamma(\alpha_n, 1)$ with $\alpha_n=\alpha+n^{-1/2}$, $\alpha \ge  1/2$. We have
$$
\frac{\dd Q_{n1}}{\dd P_0}= \frac{\Gamma(\alpha)}{\Gamma(\alpha_n)} x^{1/\sqrt{n}},
$$
$x \ge 0$. Following (\ref{radon_nikodym}), we equate this Radon-Nikodym derivative to $1+n^{-1/2}h_n(x)$, obtaining
$$
h_n(x)=n^{1/2} \Big( \frac{\Gamma(\alpha)}{\Gamma(\alpha_n)} x^{1/\sqrt{n}}-1 \Big) ,
$$
for $x \ge 0$. Recall the \textit{digamma function} \cite[Section 5.2]{ref1}
$$
\psi(z):=\frac{\dd}{\dd z} \log \Gamma(z)=\frac{\Gamma'(z)}{\Gamma(z)}
$$
$z > 0$.  Applying L'Hospital's rule, we obtain 
\begin{align*}
h(x):= \lim_{n \rightarrow \infty} h_n(x)&=\lim_{n \rightarrow \infty} n^{1/2} \Big(  \frac{\Gamma(\alpha)}{\Gamma(\alpha + n^{-1/2})}x^{1/\sqrt{n}}-1 \Big)\\
&= \log x-\psi(\alpha),
\end{align*}
$x > 0$.  To calculate $E_{P_0} |h_n|^4$, we apply the binomial expansion, obtaining 
\begin{align*}
\Big| n^{1/2} \Big( \frac{\Gamma(\alpha)}{\Gamma(\alpha + n^{-1/2})} x^{1/\sqrt{n}}-1 \Big) \Big|^4 = n^{2} \sum_{j=0}^4 (-1)^j \binom{4}{j} \Big(\frac{\Gamma(\alpha)}{\Gamma(\alpha+n^{-1/2})} \Big)^j x^{j/\sqrt{n}} ,
\end{align*}
thus, 
\begin{align}
\label{fourthmoment}
E_{P_0} |h_n|^4 &=n^{2} \sum_{j=0}^4 (-1)^j \binom{4}{j} \Big(\frac{\Gamma(\alpha)}{\Gamma(\alpha+n^{-1/2})} \Big)^j \, \frac{\Gamma(\alpha + jn^{-1/2})}{\Gamma(\alpha)}
\end{align}
Next, the Taylor expansion of $\Gamma(\alpha)/\Gamma(\alpha+n^{-1/2})$ for sufficiently large values of $n$ is 
\begin{eqnarray}
\label{taylorexpgammaalter}
\frac{\Gamma(\alpha)}{\Gamma(\alpha+n^{-1/2})}= \sum_{j=0}^4 a_j n^{-j/2} + o(n^{-2}),
\end{eqnarray}
where $a_0 = 1$.  
Recall the \textit{Hurwitz zeta function} \cite[Section 25.11]{ref1},
$$
\zeta(s,z)=\sum_{k=0}^{\infty} (k+z)^{-s},
$$
for $s>1$ and $z>0$. It is known \cite[Equation (25.11.12)]{ref1} that for $n \in \bN$ and $z > 0$,
$$
\psi^{(n)}(z)=(-1)^{n+1} n! \zeta(n+1,z).
$$
After lengthy but straightforward calculations, we obtain 
\begin{align*}
a_1&=-\psi(\alpha),\ \quad a_2=\frac{1}{2} \psi^2(\alpha)-\frac{1}{2}\zeta(2,\alpha),\\
a_3&=-\frac{1}{6}\psi^3(\alpha) +\frac{1}{3}\zeta(3,\alpha)+\frac{1}{2}\psi(\alpha)\zeta(2,\alpha),\\
\intertext{and}
a_4&=-\frac{1}{4}\zeta(4,\alpha)+\frac{1}{8}\zeta^2(2,\alpha)-\frac{1}{3}\zeta(3,\alpha)\psi(\alpha)-\frac{1}{4}\psi^2(\alpha)\zeta(2,\alpha)+\frac{1}{24}\psi^4(\alpha).
\end{align*}
Next, we substitute the Taylor expansion (\ref{taylorexpgammaalter}) in (\ref{fourthmoment}) and then take the limit as $n \to \infty$. By applying L'Hospital's rule four times then, after some lengthy but straightforward calculations, we obtain 
$$
\lim_{n \to \infty} E_{P_0} |h_n|^4 = 9a_1^4 + 24a_2^2 + 24a_1a_3 - 36a_1^2a_2 - 24a_4.
$$
Thus, $E_{P_0} |h_n^4|$ is a bounded sequence, and therefore $\sup_{n \in \bN} E_{P_0} |h_n|^4 < \infty$.

\subsubsection{Contaminated gamma models}
\label{contaminatedgammamodels}

Consider the contamination model, 
\begin{equation}
\label{Qn1contamination}
Q_{n1} = (1 - n^{-1/2}) P_0 + n^{-1/2} Gamma(2\alpha,1),
\end{equation}
where, as usual, $\alpha \ge 1/2$.  We have
\begin{align*}
\frac{\dd Q_{n1}}{\dd P_0}&=n^{-1/2} \Big( \frac{\Gamma(\alpha)}{\Gamma(2\alpha)} x^{\alpha}-1 \Big) +1, 
\end{align*}
for $x \ge 0$.  Following (\ref{radon_nikodym}), we equate this Radon-Nikodym derivative to $1+n^{-1/2}h_n(x)$, obtaining
$$
h_n(x)=\frac{\Gamma(\alpha)}{\Gamma(2\alpha)} x^{\alpha}-1, 
$$
for $x \ge 0$. Thus, 
$$
h(x):=\lim_{n \rightarrow \infty} h_n(x)=\frac{\Gamma(\alpha)}{\Gamma(2\alpha)} x^{\alpha}-1,
$$
$x \ge 0$.  Since 
$$
E_{P_0} |h_n|^4 = \int^{\infty}_0 \Big( \frac{\Gamma(\alpha)}{\Gamma(2\alpha)} x^{\alpha}-1 \Big)^4 \, \dd P_0(x)
$$
clearly is finite and does not depend on $n$ then $\sup_{n \in \bN} E_{P_0} |h_n|^4 < \infty$.  

We note also that the model (\ref{Qn1contamination}) is a special case of the contamination model 
$$
Q_{n2} = (1-n^{-1/2}) P_0 + n^{-1/2} P_1,
$$
where $P_1$ is a probability measure dominated by $P_0$, and $\int {(\dd P_1/ \dd P_0)^4}\, \dd P_0 < \infty$.  The preceding calculations can also be done for many choices of $P_1$, so we deduce that the Assumptions \ref{2assumptions} also hold for broad classes of the model $Q_{n2}$.

\subsection{The distribution of the test statistic under contiguous alternatives}

Let $\boldsymbol{P}_n=P_0 \otimes \cdots \otimes P_0$ and $\boldsymbol{Q}_n=Q_{n1} \otimes \cdots \otimes Q_{n1}$, where $P_0=Gamma(\alpha,1)$, $\alpha \ge 1/2$, be the $n$-fold product probability measures of $P_0$ and $Q_{n1}$ respectively. 

\begin{theorem}
\label{theoremcontiguous}
Let $X_{n1},\dotsc, X_{nn}$, $n \in \bN$, be a triangular array of non-negative row-wise i.i.d. random variables, where $X_{nj} \equiv X_j$, $j=1,\ldots,n$. We assume that the distribution of $X_{nj}$ is $Q_{n1}$ for every $j=1,\dotsc,n$. Further, let $Z_n:=\{ Z_n(t),t \ge 0 \}$ be a stochastic process with 
$$
Z_{n}(t)=\frac{1}{\sqrt{n}} \sum_{j=1}^n \Big[ \Gamma(\alpha)(tX_{nj}/\bar{X}_n)^{(1-\alpha)/2}J_{\alpha-1}\big(2(tX_{nj}/\bar{X}_n)^{1/2}\big)-e^{-t/\alpha} \Big],
$$
$t \ge 0$.  Under the assumptions of Section \ref{assumptionscontiguous}, there exists a centered Gaussian process $Z:=\{Z(t), t \ge 0\}$ with sample paths in $L^2$ and the covariance function $K(s,t)$ in (\ref{covariance}), and a function
$$
c(t)=\int^{\infty}_0 { \Big[ \Gamma(\alpha)(tx/\alpha)^{(1-\alpha)/2}J_{\alpha-1}(2(tx/\alpha)^{1/2})+\frac{(x-\alpha)}{\alpha^2} te^{-t/\alpha}-e^{-t/\alpha} \Big] h(x)}\,\dd P_0(x), 
$$
$t \ge 0$, such that
$Z_n \xrightarrow{d} Z +c$ in $L^2$.  Moreover, as $n \to \infty$, 
$$ 
T^2_{n} \xrightarrow{d} \int^\infty_0 \big(Z(t)+c(t)\big)^2 \, \dd P_0(t).
$$ 
 \end{theorem}
 
 \medskip
 	
We note that the proof of this theorem and the subsequent results can be obtained by following the approach in \cite[pp. 79--91]{ref27}.  In order to maintain a relatively self-contained presentation, we provide some of the details here.

Before proceeding to those details, we will present some preliminary results.  Consider the log-likelihood ratio,
\begin{align*}
\Lambda_n(X_{n1},\dotsc,X_{nn}) := \log \frac{ \dd \boldsymbol{Q}_n(X_{n1},\dotsc,X_{nn})}{\dd\boldsymbol{P}_n(X_{n1},\dotsc,X_{nn})}.
\end{align*}
From the definition of $\boldsymbol{P}_n$ and $\boldsymbol{Q}_n$, we obtain 
\begin{align*}
\Lambda_n(X_{n1},\dotsc,X_{nn}) &= \log \prod_{j=1}^n (1+n^{-1/2}h_n(X_{nj})) \\
&=\sum_{j=1}^n \log(1+n^{-1/2}h_n(X_{nj})).
\end{align*}
Since $\Lambda_n=-\infty$ if and only if $1+n^{-1/2}h_n(X_{nj} = 0$ for some $j$, we obtain 
\begin{align*}
\boldsymbol{P}_n (\Lambda_n = -\infty) &= \boldsymbol{P}_n\Big(\bigcup_{j=1}^n \{1+n^{-1/2}h_n(X_{n1})=0\}\Big) \\
&\le n \boldsymbol{P}_n(1+n^{-1/2}h_n(X_{n1})=0) \\
&= n \boldsymbol{P}_n(h_n(X_{n1})=-n^{1/2}).
\end{align*}
Denote the indicator function of an event $A$ by $A$.  Since $h_n(X_{n1})=-n^{1/2}$ if and only if $|n^{-1/2}h_n(X_{n1})|^4=1$ then 
\begin{align*}
n \boldsymbol{P}_n(h_n(X_{n1}) = -n^{1/2}) &= n E_{P_0}  ( |n^{-1/2}h_n(X_{n1})|^4 I(h_n(X_{n1})=-n^{-1/2}))\\
&\le n^{-1} E_{P_0}  ( |h_n(X_{n1})|^4)\\
&\le n^{-1} \sup_{n \in \mathbb{N}} E_{P_0} ( |h_n|^4).
\end{align*}
Under the assumption that $\sup_{n \in \bN} E_{P_0} (|h_n|^4) < \infty$, we obtain $n \boldsymbol{P}_n(h_n(X_{n1}) \to 0$ as $n \to \infty$.  Therefore, without loss of generality, we shall assume that $\Lambda_n > -\infty$ and $1+n^{-1/2}h_n(X_{nj}) > 0$ for all $j=1,\dotsc,n$ and $n \ge 1$ (see \cite[p. 140, Appendix D.2]{ref27} or \cite[p. 303, Example 6.118]{wittingmuller}).

The Taylor expansion of order $2$ of the function $\log(1+n^{-1/2}h_n(X_{nj}))$, at $1$ is 
$$
\log(1+n^{-1/2}h_n(X_{nj}))=n^{-1/2}h_n(X_{nj})-2^{-1}n^{-1}h_n^2(X_{nj}) + R(h_n(X_{nj})),
$$ 
with remainder term 
$$
R(h_n(X_{nj})) = \frac13 n^{-3/2}\big(1 + n^{-1/2} t_{nj}(X_{nj}) h_n(X_{nj})\big)^{-3} h_n^3(X_{nj}),
$$
where $t_{nj} : [0, \infty) \to [0,1]$ is a measurable function.  Therefore, 
$$
\Lambda_n(X_{n1},\dotsc,X_{nn})=\sum_{j=1}^{n} \big[ n^{-1/2}h_n(X_{nj})-2^{-1}n^{-1}h_n^2(X_{nj}) + R(h_n(X_{nj})) \big].
$$
In the following result, we use the notation $\sigma^2:=\int {h^2}\,dP_0$. 

\smallskip

\begin{lemma}
\label{lemmacontiguous} As $n \to \infty$, 
\begin{itemize}
\item[(1)]
$n^{-1/2} \sum_{j=1}^{n} h_n(X_{nj}) \xrightarrow{d} \mathcal{N}(0,\sigma^2)$ in $\boldsymbol{P}_n$-distribution.
\item[(2)]
$n^{-1} \sum_{j=1}^{n} h_n^2(X_{nj}) \rightarrow \sigma^2$ in $\boldsymbol{P}_n$-probability.
\item[(3)]
$\sum_{j=1}^{n} R(h_n(X_{nj})) \rightarrow 0$ in $\boldsymbol{P}_n$-probability.
\end{itemize}
\end{lemma}

\medskip

\noindent
The proofs of these results are given in \cite[pp. 80-83]{ref27}.  Combining these three results, we conclude that under $\boldsymbol{P}_n$, 
\begin{equation}
\label{convln_contiguous}
\Lambda_n(X_{n1},\dotsc,X_{nn}) \xrightarrow{d} \mathcal{N}(-\tfrac12 \sigma^2,\sigma^2).
\end{equation}

We introduced in Section \ref{sectionlimitingnull} the stochastic process
$$
Z_{n}(t)=\frac{1}{\sqrt{n}} \sum_{j=1}^n \big[ \Gamma(\alpha)(tX_{j}/\bar{X}_n)^{(1-\alpha)/2}J_{\alpha-1}(2(tX_{j}/\bar{X}_n)^{1/2}\big) - e^{-t/\alpha} \big], 
$$
$t \ge 0$.  Also, we introduced in Theorem \ref{limitingnulldistribution} the centered stochastic process 
$$
Z_{n,3}(t)=\frac{1}{\sqrt{n}}\sum_{j=1}^{n} \big[ \Gamma(\alpha)(tX_{j}/\alpha)^{(1-\alpha)/2}J_{\alpha-1}( 2(tX_{j}/\alpha)^{1/2}) +\frac{(X_{j}-\alpha)}{\alpha^2}te^{-t/\alpha}-e^{-t/\alpha}\big],
$$
$t \ge 0$. We proved that there exists a centered Gaussian process $Z:=\{ Z(t), t \ge 0 \}$ with sample paths in $L^2$ and with the covariance function $K(s,t)$ given in (\ref{covariance}) such that, under $\boldsymbol{P}_n$, $\|Z_n -Z_{n,3}\|_{L^2} \xrightarrow{p} 0$ and $Z_{n,3} \xrightarrow{d} Z$ in $L^2$.  For $k \in \mathbb{N}$ and $t_1,\dotsc,t_k \ge 0$, it follows from the multivariate Central Limit Theorem that $\big(Z_{n,3}(t_1),\dotsc,Z_{n,3}(t_k)\big)' \xrightarrow{d} \mathcal{N}_k \big(\boldsymbol{0},\Sigma\big)$ under $\boldsymbol{P}_n$, where $\Sigma = \big(K(t_i,t_j)\big)_{1 \le i,j \le k}$ is the $k \times k$ positive definite matrix with $(i,j)$th entry $K(t_i,t_j)$.  

Let $\|\cdot\|_{\bR^{k+1}}$ denote the standard Euclidean norm on $\bR^{k+1}$.  Then, it follows from Lemma \ref{lemmacontiguous} that 
\begin{align}
\Big\|\big(&Z_{n,3}(t_1),\dotsc,Z_{n,3}(t_k),\Lambda_n\big)' \nonumber\\
& \quad - \Big(Z_{n,3}(t_1),\dotsc,Z_{n,3}(t_k), \sum_{j=1}^n \big[ n^{-1/2} h_n(X_{nj})-2^{-1} n^{-1}h_n^2(X_{nj})\big] \Big)' \Big\|_{\mathbb{R}^{k+1}}\nonumber\\
\label{convprobzn3ln}
&=\Big( \Lambda_n - \sum_{j=1}^n \big[ n^{-1/2} h_n(X_{nj})-2^{-1} n^{-1}h_n^2(X_{nj})\big] \Big)^2 \nonumber\\
& =\Big(\sum_{j=1}^{n} R(h_n(X_{nj})\Big)^2 \rightarrow 0, \quad\text{in $\boldsymbol{P}_n$-probability}.
\end{align}

\medskip

\begin{lemma}
\label{lemma4contiguous}
For $t \ge 0$, define 
\begin{equation}
\label{cofti}
c(t) :=\lim_{n \rightarrow \infty} \Cov \Big( Z_{n,3}(t), \sum_{j=1}^n \big[ n^{-1/2}  h_n(X_{nj})-2^{-1}n^{-1}h_n^2(X_{nj})\big] \Big),
\end{equation}
and set $\boldsymbol{c} = \big( c(t_1),\dotsc, c(t_k) \big)'$.  Then, under $\boldsymbol{P}_n$,
\begin{multline}
\label{lemma4contiguous_eq}
\Big( Z_{n,3}(t_1),\dotsc,Z_{n,3}(t_k), \sum_{j=1}^n \big[ n^{-1/2} h_n(X_{nj})-2^{-1} n^{-1}h_n^2(X_{nj})\big] \Big)' \\
\xrightarrow{d} \mathcal{N}_{k+1} \left( (0,\dotsc,0,-2^{-1}\sigma^2)', \begin{bmatrix} \Sigma & \boldsymbol{c} \\ \boldsymbol{c'} & \sigma^2 \end{bmatrix} \right)
\end{multline}
and 
\begin{equation}
\label{convz3ln_contiguous}
\big(Z_{n,3}(t_1),\dotsc,Z_{n,3}(t_k), \Lambda_n\big)' \xrightarrow{d} \mathcal{N}_{k+1} \left( (0,\dotsc,0,-2^{-1}\sigma^2)', \begin{bmatrix} \Sigma & \boldsymbol{c} \\ \boldsymbol{c'} & \sigma^2 \end{bmatrix} \right).
\end{equation}
\end{lemma}

\Pro
Substituting for $Z_{n,3}$ in (\ref{cofti}), applying the assumptions in Subsection \ref{assumptionscontiguous}, and carrying out some straightforward calculations, we obtain 
$$
c(t) = \int^{\infty}_0 { \Big[ \Gamma(\alpha)(tx/\alpha)^{(1-\alpha)/2}J_{\alpha-1}(2(tx/\alpha)^{1/2})+\frac{(x-\alpha)}{\alpha^2} te^{-t/\alpha}-e^{-t/\alpha} \Big] h(x)}\,\dd P_0(x), 
$$
$t \ge 0$.  Letting
$$
g(t,X_{nj}) := \Gamma(\alpha)(t X_{nj}/\alpha)^{(1-\alpha)/2}J_{\alpha-1}(2(t X_{nj}/\alpha)^{1/2}) +\frac{(X_{nj}-\alpha)}{\alpha^2} t e^{-t/\alpha} -e^{-t/\alpha};
$$
then 
$$
Z_{n,3}(t_i)=\frac{1}{\sqrt{n}} \sum_{j=1}^n g(t_i, X_{nj}).
$$
To establish (\ref{lemma4contiguous_eq}), we will apply the Cram\'er-Wold device.  Then it suffices to establish that for every $\boldsymbol{u}=(u_1,\dotsc, u_{k+1})' \in \mathbb{R}^{k+1}$, 
\begin{multline}
\label{cramerwoldconv}
\frac{1}{\sqrt{n}} \sum_{j=1}^n \big(g(t_1,X_{nj}),\dotsc, g(t_k,X_{nj}), h_n(X_{nj})-2^{-1} n^{-1/2}h_n^2(X_{nj}) \big) \boldsymbol{u} \\
\xrightarrow{d} \mathcal{N}_{k+1} \left( (0,\dotsc,0,-2^{-1}\sigma^2) \boldsymbol{u}, \boldsymbol{u}'\begin{bmatrix} \Sigma & \boldsymbol{c} \\ \boldsymbol{c'} & \sigma^2 \end{bmatrix} \boldsymbol{u} \right).
\end{multline}
Now, let $Y_1,Y_2,\dotsc$ be i.i.d. $P_0$-distributed random variables, and define
$$
k_n(Y_j) = \sum_{l=1}^k g(t_l,Y_j) u_l + \big( h_n(Y_j)-2^{-1} n^{-1/2}h_n^2(Y_j) \big) u_{k+1}.
$$
Under $\boldsymbol{P}_n$, $n^{-1/2} \sum_{j=1}^n k_n(Y_j)$ has the same distribution as
$$
\frac{1}{\sqrt{n}} \sum_{j=1}^n \big(g(t_1,X_{nj}),\dotsc, g(t_k,X_{nj}), h_n(X_{nj})-2^{-1} n^{-1/2}h_n^2(X_{nj}) \big)\boldsymbol{u},
$$
$\boldsymbol{u} \in \mathbb{R}^{k+1}$.  Since $E(g(t_i, Y_1))=0$, $i=1,\dotsc, k$, then 
$$
\mu_n := E(k_n(Y_1))=-2^{-1} n^{-1/2} u_{k+1} E(h_n^2(Y_1)).
$$
Denote by $\tau^2_n$ the variance of $k_n(Y_1)$.  Then, 
\begin{align*}
\tau^2_n &=\sum_{i=1}^k \sum_{j=1}^k \Cov \big(g(t_i,Y_1),g(t_j,Y_1)\big) u_i u_j + \Var (h_n(Y_1)) u^2_{k+1} \\
& \qquad + (4n)^{-1} \Var (h^2_n(Y_1)) u^2_{k+1} - n^{-1/2} \Cov (h_n(Y_1), h^2_n(Y_1)) u^2_{k+1} \\
& \qquad\qquad + 2 u_{k+1} \sum_{i=1}^k \Cov \Big(g(t_i, Y_1),\big( h_n(Y_1)-2^{-1} n^{-1/2}h_n^2(Y_1) \big)\Big) u_i.
\end{align*}
By the assumptions (A1) and (A2) in Subsection \ref{subsec_examples}, we obtain $\Var (h^2_n(Y_1)) < \infty$ and $\Cov \big( h_n(Y_1), h^2_n(Y_1) \big) < \infty$.  Thus, as $n \rightarrow \infty$, 
\begin{eqnarray}
\label{tao_contiguous}
\tau^2_n \to \sum_{i=1}^k \sum_{j=1}^k K(t_i, t_j) u_i u_j + \sigma^2 u^2_{k+1} + 2 u_{k+1}\sum_{i=1}^k u_i c(t_i):= \tau^2.
\end{eqnarray}
Similarly, it can be shown that, as $n \rightarrow \infty$, 
\begin{align}
\label{asconvcontiguous}
(k_n&(Y_1)- \mu_n)^2 \nonumber \\
& \to \sum_{i=1}^k \sum_{j=1}^k g(t_i, Y_1) g(t_j,Y_1) u_i u_j + h(Y_1) \Big(h(Y_1) u_{k+1} + \sum_{i=1}^k g(t_i,Y_1) u_i \Big) u_{k+1},
\end{align}
$P_0$-almost surely.  In addition, we notice that 
\begin{multline}
\label{expected_asconvcontiguous}
E \Big[\sum_{i=1}^k \sum_{j=1}^k g(t_i, Y_1) g(t_j,Y_1) u_i u_j + h(Y_1) \Big(h(Y_1) u_{k+1} + \sum_{i=1}^k g(t_i, Y_1) u_i \Big) u_{k+1} \Big] \\
\ = \ \sum_{i=1}^k \sum_{j=1}^k K(t_i, t_j) u_i u_j + \sigma^2 u^2_{k+1} + 2 u_{k+1} \sum_{i=1}^k c(t_i) u_i \ \equiv \ \tau^2.
\end{multline}
For every $\epsilon > 0$, 
\begin{eqnarray}
\label{result1contiguous}
0 \le (k_n(Y_1)- \mu_n)^2  \ I(|(k_n(Y_1)- \mu_n| > \epsilon \sqrt{n} \tau_n) \le (k_n(Y_1)- \mu_n)^2.
\end{eqnarray}
Also, for every $\epsilon > 0$, 
$$
0 \le (k_n(Y_1)- \mu_n)^2 \ I(|(k_n(Y_1)- \mu_n| > \epsilon \sqrt{n} \tau_n) \le  \frac{(k_n(Y_1)- \mu_n)^4}{ \epsilon^2 n \tau^2_n},
$$
from which we conclude that as $n \rightarrow \infty$, 
\begin{eqnarray}
\label{result2contiguous}
(k_n(Y_1)- \mu_n)^2 \  I(|(k_n(Y_1)- \mu_n| > \epsilon \sqrt{n} \tau_n) \rightarrow 0,
\end{eqnarray}
$P_0$-almost surely.  As the results (\ref{tao_contiguous}) -- (\ref{result2contiguous}) are the sufficient conditions in Pratt's version of the Dominated Convergence Theorem \cite[p. 221, Theorem 5.5]{gut}, we conclude that as $n \to \infty$, 
$$
E\big( (k_n(Y_1)- \mu_n)^2  \ I(|(k_n(Y_1)- \mu_n| > \epsilon \sqrt{n} \tau_n) \big) \to 0.
$$
This result is equivalent to the Lindeberg condition, i.e. for every $\epsilon > 0$, 
$$
\lim_{n \to \infty} \frac{1}{n \tau^2_n} \sum_{j=1}^n E\big( (k_n(Y_j)- \mu_n)^2 \ I(|(k_n(Y_j)- \mu_n| > \epsilon \sqrt{n} \tau_n) \big) = 0.
$$
Thus, we deduce from the Lindeberg-Feller Central Limit Theorem that 
$$
\frac{1}{\sqrt{n} \tau_n} \sum_{j=1}^n (k_n(Y_j)-\mu_n) \xrightarrow{d} \mathcal{N}(0,1),
$$
therefore,
$$
\frac{1}{\sqrt{n}} \sum_{j=1}^n k_n(Y_j) \xrightarrow{d} \mathcal{N} \big( -2^{-1} \sigma^2 u_{k+1}, \tau^2 \big).
$$
Note also that $(0,\dotsc,0,-2^{-1}\sigma^2)\boldsymbol{u} = -2^{-1} \sigma^2 u_{k+1}$ and that 
\begin{align*}
\tau^2 &= \sum_{i=1}^k \sum_{j=1}^k u_i u_j K(t_i, t_j) + u^2_{k+1}\sigma^2 +2 u_{k+1}\sum_{i=1}^k u_i c(t_i)\\
&= \boldsymbol{u}' \begin{bmatrix} \Sigma & \boldsymbol{c} \\ \boldsymbol{c'} & \sigma^2 \end{bmatrix} \boldsymbol{u}.
\end{align*}
Therefore, (\ref{cramerwoldconv}) is proved. 

Finally, (\ref{convz3ln_contiguous}) follows from (\ref{convprobzn3ln}), Lemma \ref{lemma4contiguous}, and \cite[p. 25, Theorem 4.1]{ref21}.
$\qed$

\medskip

Now, we proceed to the proof of Theorem \ref{theoremcontiguous}.

\medskip

\noindent{\it Proof of Theorem 4.1}.  
By (\ref{convln_contiguous}) and Le Cam's first lemma (see \cite[p. 140, Theorem D.5]{ref27} or \cite[p. 311, Corollary 6.124]{wittingmuller}), $\boldsymbol{P}_n$ and $\boldsymbol{Q}_n$ are mutually contiguous.  Also, by (\ref{convz3ln_contiguous}) and Le Cam's third lemma (see \cite[p. 141, Theorem D.6]{ref27} or  \cite[p. 329, Corollary 6.139]{wittingmuller}), under $\boldsymbol{Q}_n$, 
\begin{eqnarray}
\label{convz3_contiguous}
\big(Z_{n,3}(t_1),\dotsc,Z_{n,3}(t_k) \big)' \xrightarrow{d} \mathcal{N}_{k} (\boldsymbol{c},\Sigma).
\end{eqnarray}
By \cite[p. 138, Theorem D.2]{ref27} or \cite[p. 56, Theorem 5.51]{wittingmuller}, the convergence in distribution of $Z_{n,3}$ under $\boldsymbol{P}_n$ in $L^2$ implies that $Z_{n,3}$ is tight in $L^2$ under $\boldsymbol{P}_n$. Further, since $\boldsymbol{Q}_n$ is contiguous to $\boldsymbol{P}_n$, by \cite[p. 139, Theorem D.4]{ref27} or \cite[p. 295, Theorem 6.113 (a)]{wittingmuller}, $Z_{n,3}$ is tight in $L^2$ under $\boldsymbol{Q}_n$. 

By (\ref{convz3_contiguous}) and the tightness of $Z_{n,3}$ in $L^2$ under $\boldsymbol{Q}_n$, we obtain $Z_{n,3} \xrightarrow{d} Z + c$ under $\boldsymbol{Q}_n$ (see \cite[Theorem 2, Example 4]{cremerskadelka}).  Moreover, since $\| Z_n - Z_{n,3}\|_{L^2} \xrightarrow{p} 0$ under $\boldsymbol{P}_n$ and since $\boldsymbol{Q}_n$ is contiguous to $\boldsymbol{P}_n$, we have that under $\boldsymbol{Q}_n$, $\| Z_n - Z_{n,3}\|_{L^2} \xrightarrow{p} 0$.  Therefore, by Billingsley \cite[p. 25, Theorem 4.1]{ref21}, we obtain $Z_n \xrightarrow{d} Z + c$ under $\boldsymbol{Q}_n$.  

Finally, by the Continuous Mapping Theorem \cite[p. 31, Corollary 1]{ref21}, $\lVert Z_{n} \rVert^2_{L^2} \xrightarrow{d} \lVert Z +c \rVert^2_{L^2}$ under $\boldsymbol{Q}_n$, i.e., 
$$
T^2_{n} \xrightarrow{d} \int^\infty_0 {(Z(t)+c(t))^2}\,dP_0(t)
$$
under $\boldsymbol{Q}_n$.  The proof now is complete. 
$\qed$

\section{The Efficiency of the Test}
\label{sec:efficiency}

In this section, we investigate the approximate Bahadur slope of the test statistic $T_n^2$ under local alternatives.  We show further that the Wieand condition, under which the approaches through the approximate Bahadur efficiency and the Pitman efficiency are in accord, is valid.  By applying the results of this section, we are able to calculate the approximate asymptotic relative efficiency (ARE) of the proposed test relative to alternative tests. 

Let $X_1, X_2, \dotsc $ be i.i.d, non-negative random variables with some unknown distribution $\mathcal{P}$. Henceforth, we assume that $\mathcal{P}$ is indexed by a parameter $\theta \in \Theta :=(-\eta, \eta)$, for some $\eta >0$. We let $\theta \in \Theta_0=\{\theta_0 \}=\{ 0 \}$ to represent the null hypothesis and $\theta \in \Theta_1=\Theta \setminus \{0\}$ to represent the alternative hypothesis. In Chapter 3, we showed that $T^2_n$ is scale-invariant, i.e. it does not depend on the unknown rate parameter $\lambda$. Thus, under the null hypothesis $\theta_0=0$, we assume that $X_1, X_2,\dotsc $ are i.i.d, non-negative $P_0$-distributed random variables and under the local alternatives, represented by $\theta \in \Theta_1$, we assume that $X_1, X_2,\dotsc $ are i.i.d, non-negative $P_\theta$-distributed random variables. 

The Radon-Nikodym derivative of $P_{\theta}$ with respect to $P_0$ is $\dd P_{\theta}/{\dd P_0}=1+\theta h_{\theta}$. We assume that as $\theta \rightarrow 0$, the function $h_{\theta}$ converges to some function $h$ in $L^2$. Since $ \int {(\dd P_{\theta}/ \dd P_0)}\, \dd P_0=1$, we have that 
\begin{equation}
\label{meanzero} 
\int_0^{\infty} {h_{\theta} (x)}\, \dd P_0(x)=0,
\end{equation}
for $\theta \in \Theta_1$. Further, we shall assume that for $\theta \in \Theta_1$,
\begin{equation}
\label{assumption_efficiency}
\int_0^{\infty} {x h_{\theta}(x) }\, \dd P_0(x)=0.
\end{equation}

Before we proceed to the approximate Bahadur slope of our test, we recall the definition of the exact Bahadur slope as presented by Bahadur \cite[Section 7]{bahadur2} or Taherizadeh \cite[Chapter 5]{ref27}, for general parameter spaces $\Theta$ and $\Theta_j$, $j=0,1$, and for a general sequence of test statistics $\{ U_n : n \in \bN \}$. For $\theta \in \Theta_0$, let $F_n(t)=P_{\theta} (U_n < t)$, $t \in \mathbb{R}$, be the null distribution of $U_n$.  Then {\it the level attained by} $U_n$ is $L_n := 1-F_n(U_n)$. For $\theta \in \Theta_1$, the \textit{exact Bahadur slope} of the sequence of test statistics $\{U_n: n \in \bN\}$ is 
$$
c(\theta) = - 2 \lim_{n \rightarrow \infty} n^{-1} \log L_n,
$$
whenever the limit exists (almost surely).  For $\theta \in \Theta_0$, this limit exists with $c(\theta)=0$. 

Let $\{U_{1,n}: n \in \bN \}$ and $\{U_{2,n}: n \in \bN\}$ be sequences of test statistics with exact Bahadur slopes $c_1(\theta)$ and $c_2(\theta)$ respectively. Then the \textit{exact Bahadur asymptotic relative efficiency} (ARE) of $ \{ U_{1,n} : n \in \bN \}$ with respect to $ \{ U_{2,n} : n \in \bN \}$ is 
$$
e_{1,2}^B(\theta)=\frac{c_1(\theta)}{c_2(\theta)}, 
$$
$\theta \in \Theta_1$.  If $e_{1,2}^B(\theta) >1 $, then we prefer the sequence of test statistics $\{U_{1,n}: n \in \bN\}$. 

For further results on the exact Bahadur slope, we refer to Bahadur \cite[Section 7]{bahadur2}.  In particular, \cite[Theorem 7.2]{bahadur2} provides a method of detecting the exact Bahadur slope. In general, it is not easy to compute the exact Bahadur slope and therefore, we investigate the approximate Bahadur slope.  We note that Bahadur \cite{bahadur3} showed that for $\Theta_0=\{ \theta_0 \}$, the approximate Bahadur slope is very close to the exact Bahadur slope for $\theta$ in a neighborhood of $\theta_0$, i.e. under local alternatives.

In order to find the approximate Bahadur slope of a test under local alternatives, we need to establish that the sequence of test statistics $ \{ U_n : n \in \bN \}$ is a {\it standard sequence} in the sense defined by Bahadur \cite{bahadur1}, viz., 

\begin{itemize}
\item[(i)] There exists a continuous probability distribution function $F$ such that, for each $\theta \in \Theta_0$ and all $t \in \mathbb{R}$, 
$$
\lim_{n \rightarrow \infty} P_\theta (U_n < t)=F(t), 
$$
\item[(ii)] 
There exists $a > 0$ such that
$$
-2 \lim_{t \rightarrow \infty} t^{-2} \log[1-F(t)]=a. 
$$ 
\item[(iii)] 
There exists a function $b: \Theta_1 \rightarrow (0,\infty)$, such that, for all $\theta \in \Theta_1$, $n^{-1/2} U_n \rightarrow b(\theta)$, in $P_\theta$-probability. 
\end{itemize}
The \textit{approximate Bahadur slope} of $ \{ U_n : n \in \bN \}$ is defined by
\begin{align*}
c^{(a)}(\theta):=
\begin{cases}
0, & \theta \in \Theta_0 \\
a b^2(\theta), & \theta \in \Theta_1.
\end{cases}
\end{align*}

\subsection{Approximate Bahadur slope of the test}

Now we have the following result for the test statistic $T_n^2$ arising from the Hankel transform.

\begin{theorem}
The sequence of test statistics $\{T_n : n \in \bN \}$ is a standard sequence.  Further, $a = \delta_1^{-1}$, the inverse of the largest eigenvalue of the covariance operator $\mathcal{S}$, 
\begin{equation}
\label{bsquaredfunction}
b^2(\theta) = \theta^2 \int_0^{\infty} {\Big[\int_0^{\infty} { \Gamma(\alpha) (tx/\alpha)^{(1-\alpha)/2} J_{\alpha-1}(2 (tx/\alpha)^{1/2}) h_{\theta} (x)}\, \dd P_0(x) \Big]^2}\, \dd P_0(t)
\end{equation}
and 
$$
\lim_{\theta \rightarrow 0} \frac{c^{(a)}(\theta)}{\theta^2} = \delta_1^{-1} \int_0^{\infty} {\Big[ \int_0^{\infty} \Gamma(\alpha) (tx/\alpha)^{(1-\alpha)/2} J_{\alpha-1}(2 (tx/\alpha)^{1/2} ) h(x) \, \dd P_0(x) \Big] ^2}\, \dd P_0(t).
$$
\end{theorem}

\Pro
The proof of this theorem follows along the lines of the proof of Theorem 5.4 in \cite[p. 98]{ref27}. For completeness, we provide the details here. 

First, we will establish that $\{T_n : n \in \bN \}$ is a standard sequence. In Chapter 3, we showed that the limiting null distribution of the test statistic $T^2_n$ is the same as that of $\sum_{k=1}^{\infty}\delta_k \chi^2_{1k}$, where the coefficients $\{\delta_k : k \ge 1 \}$ are the positive eigenvalues, listed in decreasing order, of the covariance operator $\mathcal{S}$ with kernel function defined in equation (\ref{covariance}), and $\{\chi^2_{1k} : k \ge 1 \}$ are i.i.d. $\chi^2_{1}$-distributed random variables. From the Monotone Convergence Theorem, we have 
$$
\lim_{N \rightarrow \infty} E \Big( \sum_{k=1}^{N}\delta_k \chi^2_{1k} \Big)=E \Big(\sum_{k=1}^{\infty}\delta_k \chi^2_{1k} \Big) = E \Big(\sum_{k=1}^{\infty}\delta_k \chi^2_{1k} \Big)=\sum_{k=1}^{\infty}\delta_k,
$$
which is finite since $\mathcal{S}$ is of trace-class.  Thus, $\sum_{k=1}^{\infty}\delta_k \chi^2_{1k}$ is almost surely a positive real random variable with continuous probability distribution function. 

By Zolotarev \cite{zolotarev}, 
\begin{eqnarray*}
1-F(t) &=&P \Big( \sum_{k=1}^{\infty} \delta_k \chi^2_{1k} > t^2 \Big)\\
&=& \prod_{k \ge 2} \Big( 1-\frac{\delta_k}{\delta_1} \Big)^{-1} \Big( \frac{\pi}{2\delta_1} \Big)^{-1/2} \ t^{-1} \exp(-t^2/2\delta_1) \ [1+ o_{p}(1)],
\end{eqnarray*}
where $o_{p}(1) \xrightarrow{t \rightarrow \infty} 0$. Therefore, 
\begin{align*}
-2 & t^{-2} \log[1-F(t)] \\
& = 2 t^{-2} \Big[ \sum_{k \ge 2} \log \Big( 1-\frac{\delta_k}{\delta_1} \Big) + 2^{-1}\log \Big( \frac{\pi}{2\delta_1} \Big) + \log t - \log(1+o_{p}(1)) \Big] + \delta_1^{-1},
\end{align*}
which converges to $\delta_1^{-1}$ as $t \rightarrow \infty$.  

By assumption (\ref{assumption_efficiency}), for $\theta \in \Theta_1$, 
\begin{align*}
E_\theta (X_1):=\int^{\infty}_0 { x }\, \dd P_{\theta}(x)=\int^{\infty}_0 { x (1+\theta h_{\theta}(x))}\, \dd P_0(x)=\alpha.
\end{align*}
From the proof of Theorem \ref{consistency}, we have that
\begin{align*}
n^{-1} T^2_n \xrightarrow{p} \ \int^\infty_0{\Big( E_{\theta}[\Gamma(\alpha)(tX_1/\alpha)^{(1-\alpha)/2}J_{\alpha-1}(2(tX_1/\alpha)^{1/2})]-e^{-t/\alpha}\Big)^2}\, \dd P_0(t).
\end{align*}
Since $\dd P_\theta/\dd P_0 = 1 + \theta h_\theta$ then, by (\ref{meanzero}), 
\begin{multline*}
E_{\theta}[\Gamma(\alpha)(tX_1/\alpha)^{(1-\alpha)/2}J_{\alpha-1}(2(tX_1/\alpha)^{1/2})]\\
=e^{-t/\alpha} + \theta \int^{\infty}_0 { \Gamma(\alpha)(tx/\alpha)^{(1-\alpha)/2}J_{\alpha-1}(2(tx/\alpha)^{1/2}) h_{\theta}(x)}\, \dd P_{0}(x),
\end{multline*}
and then it follows that $n^{-1} T^2_n \xrightarrow{p} b^2(\theta)$, the function defined in (\ref{bsquaredfunction}).  Therefore, $n^{-1/2} T_n \xrightarrow{p} b(\theta)$ in $P_\theta$-probability, so the sequence of test statistics $\{T_n: n \in \bN\}$ is a standard sequence.  

Finally, we find the limiting approximate Bahadur slope, as $\theta \rightarrow 0$. By applying the Cauchy-Schwarz inequality, the inequality  (\ref{besselineq2}), and the $L^2$-convergence of $h_\theta$ to $h$, it is straightforward to establish that 
$$
\limsup_{\theta \to 0} \Big|\frac{b^2(\theta)}{\theta^2}-\int^{\infty}_0 \Big[\int^\infty_0 \Gamma(\alpha)(tx/\alpha)^{(1-\alpha)/2}J_{\alpha-1}(2(tx/\alpha)^{1/2}) h(x) \dd P_0(x)\Big]^2 \dd P_{0}(t)\Big| = 0.
$$
Therefore, 
$$
\lim_{\theta \rightarrow 0} \Big|\frac{b^2(\theta)}{\theta^2}-\int^{\infty}_0 \Big[\int^\infty_0 \Gamma(\alpha)(tx/\alpha)^{(1-\alpha)/2}J_{\alpha-1}(2(tx/\alpha)^{1/2}) h(x) \, \dd P_0(x)\Big]^2 \, \dd P_{0}(t)\Big|=0.
$$
The proof now is complete. $\qed$

\subsection{Wieand's condition}

Wieand \cite{wieand} showed that if two standard sequences of test statistics $ \{ U_{1,n}: n \in \bN \}$ and $ \{ U_{2,n}: n \in \bN \}$ satisfy an additional condition, now called the \textit{Wieand condition}, then the limiting approximate Bahadur efficiency is in accord with the limiting Pitman efficiency, as the level of significance decreases to $0$.  In this section, we show that Wieand's condition is valid for our sequence of test statistics $ \{ T_n : n \in \bN \}$. 

Before we proceed to establish the Wieand condition, we will briefly describe Pitman's asymptotic relative efficiency.  Here, we use the same notation as in \cite[Chapter 5, Section 5.3]{ref27}. Let $\Theta \subset \bR$ and for $\theta_0 \in \Theta$, set $\Theta_0=\{\theta_0 \}$  and $\Theta_1=\Theta \setminus \{\theta_0\}$. Suppose that $\{\theta_j : j \in \bN \}$, $\theta_j \in \Theta_1$, is a sequence of alternatives with $\lim_{j \rightarrow \infty} \theta_j=\theta_0$, that $\{ \beta_{1,j} : j  \in \bN \}$ and $\{ \beta_{2,j} : j  \in \bN \}$ are real sequences with $\lim_{j \to \infty} \beta_{1,j}=\beta$ and $\lim_{j \to \infty} \beta_{2,j}=\beta$, where $\beta \in (0,1)$. Further, for two standard sequences of test statistics $ \{ U_{1,n} : n \in \bN \}$ and $ \{ U_{2,n} : n \in \bN \}$, and for $i=1,2$ and $j \in \bN$, define 
$$
N(i,j) := \min \big\{ m \in \mathbb{N}: P_{\theta_j} (U_{i, n} > t_{i, n}) \ge \beta_{i,j} \hbox{ for all } n \ge m \big\}, 
$$
where $t_{i,n}$ is the $(1-\gamma)$-quantile of the null distribution of $U_{i,n}$.  The interpretation is that $N(i,j)$ is the smallest sample size necessary to ensure that the test based on $U_{i,n}$ has power $\beta_{i,j}$.  

The \textit{Pitman asymptotic relative efficiency} (ARE) of $\{ U_{1, n} : n \in \bN \}$ with respect to $\{ U_{2, n} : n \in \bN \}$ is 
$$
e_{1,2}^P(\gamma, \beta)=\lim_{j \rightarrow \infty} \frac{N(2,j)}{N(1,j)},
$$
assuming that the limit exists. If $e_{1,2}^P(\gamma, \beta) > 1 $, then we would prefer the sequence of test statistics $\{ U_{n,1} : n \in \bN \}$. 

\begin{theorem}
\label{wieandTn}
The sequence of test statistics $\{ T_n : n \in \bN \}$ satisfies the Wieand condition, i.e., there exists a constant $\theta^{*} > 0$ such that for any $\epsilon > 0$ and $\gamma \in (0,1)$, there exists a constant $C > 0$ such that 
$$
P(|n^{-1/2}T_n - b(\theta)| < \epsilon b(\theta)) > 1-\gamma.
$$
for any $\theta \in \Theta_1 \cap (-\theta^{*}, \theta^{*})$ and $n > C/b^2(\theta)$.  
\end{theorem}

\Pro
For $t \ge 0$ and $\theta \in \Theta$, consider the Hankel transform, 
$$
\mathcal{H}_{X_1, \theta} (t) = E_{\theta} [\Gamma(\alpha)(tX_1/\alpha)^{(1-\alpha)/2}J_{\alpha-1}(2(tX_1/\alpha)^{1/2}) ].
$$
We have that
$$
 n^{-1/2} T_n =\Big[ \int^{\infty}_0 {  \Big( \frac{1}{n} \sum_{j=1}^n \Gamma(\alpha)(tY_j)^{(1-\alpha)/2}J_{\alpha-1}\big( 2(tY_j)^{1/2}) -e^{-t/\alpha}  \Big)^2}\, \dd P_0(t)  \Big]^{1/2}.
$$
By adding and subtracting the term $\mathcal{H}_{X, \theta} (t)$ inside the squared term, and then applying Minkowski's inequality, we obtain 
\begin{align}
\label{wieand_ineq1}
 n^{-1/2} T_n & \le \ \Big[ \int^{\infty}_0 {  \Big( \frac{1}{n} \sum_{j=1}^n \Gamma(\alpha)(tY_j)^{(1-\alpha)/2}J_{\alpha-1}(2(tY_j)^{1/2}\big)-\mathcal{H}_{X_1, \theta}(t)\Big)^2}\, \dd P_0(t)\Big]^{1/2}\nonumber\\
& \quad\quad\quad + \Big[ \int^{\infty}_0 { \big( \mathcal{H}_{X_1, \theta}(t) -e^{-t/\alpha} \big)^2}\, \dd P_0(t)  \Big]^{1/2}.
\end{align}
Now set 
\begin{align*}
b(\theta) := \Big[ \int^{\infty}_0 { \big( \mathcal{H}_{X_1, \theta}(t) -e^{-t/\alpha} \big)^2}\, \dd P_0(t)  \Big]^{1/2}.
\end{align*}
By adding and subtracting the term $$ \frac{1}{n} \sum_{j=1}^n \Gamma(\alpha)(tY_j)^{(1-\alpha)/2}J_{\alpha-1}\big(2(tY_j)^{1/2}\big)$$ inside the squared term, and then applying Minkowski's inequality again, we get
\begin{align}
\label{wieand_ineq2}
b(\theta) & \le n^{-1/2} T_n \nonumber\\
& \quad + \Big[ \int^{\infty}_0 {  \Big( \frac{1}{n} \sum_{j=1}^n \Gamma(\alpha)(tY_j)^{(1-\alpha)/2}J_{\alpha-1}\big(2(tY_j)^{1/2}\big)-\mathcal{H}_{X_1, \theta}(t) \Big)^2}\, \dd P_0(t)\Big]^{1/2}.
\end{align}
Combining (\ref{wieand_ineq1}) and (\ref{wieand_ineq2}), we conclude that
\begin{multline}
\label{wieand_ineq3}
| n^{-1/2}T_n - b(\theta)| \\
\le \ \Big[ \int^{\infty}_0 {  \Big( \frac{1}{n} \sum_{j=1}^n \Gamma(\alpha)(tY_j)^{(1-\alpha)/2}J_{\alpha-1}\big(2(tY_j)^{1/2}\big)-\mathcal{H}_{X_1, \theta}(t) \Big)^2}\, \dd P_0(t) \Big]^{1/2}.
\end{multline}
Further, by subtracting and adding the term $$\frac{1}{n} \sum_{j=1}^n \Gamma(\alpha)(tX_{j}/\alpha)^{(1-\alpha)/2}J_{\alpha-1}\big(2(tX_{j}/\alpha)^{1/2}\big)$$ inside the squared term
$$
\Big( \frac{1}{n} \sum_{j=1}^n \Gamma(\alpha)(tY_j)^{(1-\alpha)/2}J_{\alpha-1}\big( 2(tY_j)^{1/2}\big) -\mathcal{H}_{X_1, \theta}(t)\Big)^2,
$$
and then applying the Cauchy-Schwarz inequality, we obtain 
\begin{align}
\Big( \frac{1}{n} \sum_{j=1}^n & \Gamma(\alpha) (tY_j)^{(1-\alpha)/2}J_{\alpha-1}\big( 2(tY_j)^{1/2}\big) -\mathcal{H}_{X_1, \theta}(t)\Big)^2 \nonumber\\
\label{wieand_ineq4}
 & \le 2 \ \Big[ \frac{1}{n} \sum_{j=1}^n \Gamma(\alpha) \Big( (tY_j)^{(1-\alpha)/2}J_{\alpha-1}\big( 2(tY_j)^{1/2}\big)\nonumber\\
 & \quad\quad\quad\quad\quad\quad\quad\quad\quad\quad\quad -(tX_{j}/\alpha)^{(1-\alpha)/2}J_{\alpha-1}\big( 2(tX_{j}/\alpha)^{1/2}\big) \Big) \Big]^2 \nonumber\\
 & \quad + 2 \Big[ \frac{1}{n} \sum_{j=1}^n \Gamma(\alpha) (tX_{j}/\alpha)^{(1-\alpha)/2}J_{\alpha-1}\big( 2(tX_{j}/\alpha)^{1/2}\big)  - \mathcal{H}_{X_1, \theta}(t)\Big]^2.
\end{align}
Next, by proceeding as in the first part of the proof of Theorem \ref{consistency}, we get 
\begin{multline*}
\frac{1}{n}\sum_{j=1}^n \Gamma(\alpha) |(tY_j)^{(1-\alpha)/2} J_{\alpha-1}(2(tY_j)^{1/2})-(tX_{j}/\alpha)^{(1-\alpha)/2} J_{\alpha-1} (2(tX_j/\alpha)^{1/2}) | \\
\le  \frac{2 \Gamma(\alpha)}{\pi^{1/2} \Gamma(\alpha+\frac{1}{2})} \ t^{1/2} \ | {\bar{X}_{n}}^{-1/2}-\alpha^{-1/2} | \ \frac{1}{n} \sum_{j=1}^n X_j^{1/2}.
\end{multline*}
By the Cauchy-Schwarz inequality, $n^{-1} \sum_{j=1}^n X_j^{1/2} \le \bar{X}_n^{1/2}$, and thus, 
\begin{align}
\label{inequality_w}
| \bar{X}_n^{-1/2}-\alpha^{-1/2} | \ \frac{1}{n} \sum_{j=1}^n X_j^{1/2} \le  | 1- \alpha^{-1/2} \bar{X}_n^{1/2} | =\alpha^{-1/2} \  \frac{| \bar{X}_n-\alpha |}{ |\bar{X}_n^{1/2} +\alpha^{1/2}| }.
\end{align}
Since $ |\bar{X}_n^{1/2} +\alpha^{1/2}| \ge \alpha^{1/2}$, we see that (\ref{inequality_w}) is bounded above by $\alpha^{-1/4} \, | \bar{X}_n-\alpha |$. Therefore,
\begin{multline}
\label{wieand_ineq5}
\frac{1}{n}\sum_{j=1}^n \Gamma(\alpha) |(tY_j)^{(1-\alpha)/2} J_{\alpha-1}(2(tY_j)^{1/2})-(tX_{j}/\alpha)^{(1-\alpha)/2} J_{\alpha-1} (2(tX_j/\alpha)^{1/2}) | \\
\le \frac{2\Gamma(\alpha)}{\alpha^{1/4} \pi^{1/2} \Gamma(\alpha+\tfrac{1}{2})} \ t^{1/2} \ | \bar{X}_n-\alpha|.
\end{multline}
By (\ref{wieand_ineq3}), Markov's inequality, and Fubini's theorem, 
\begin{multline}
\label{inequality_w2}
P\big(| n^{-1/2}T_n - b(\theta)| \le \epsilon b(\theta) \big) \\
\ge 1 - \frac{1}{\epsilon^2 b^2(\theta)} \int^{\infty}_0 E_{\theta} \Big( \frac{1}{n} \sum_{j=1}^n \Gamma(\alpha)(tY_j)^{(1-\alpha)/2}J_{\alpha-1}\big( 2(tY_j)^{1/2}\big)-\mathcal{H}_{X_1, \theta}(t) \Big)^2 \, \dd P_0(t).
\end{multline}
By (\ref{wieand_ineq4}), (\ref{wieand_ineq5}), and by noting that $\int_0^{\infty} t \dd P_0(t) =\alpha$, we conclude that (\ref{inequality_w2}) is greater than or equal to 
\begin{multline*}
1 - \frac{1}{\epsilon^2 b^2(\theta)} \Big[ 8 \alpha^{1/2} \Big( \frac{ \Gamma(\alpha)}{\pi^{1/2} \Gamma(\alpha+\tfrac{1}{2})} \Big)^2 \ E_{\theta} (\bar{X}_n-\alpha)^2 \\
+ 2 \int^{\infty}_0 { E_{\theta} \Big( \frac{1}{n} \sum_{j=1}^n \Gamma(\alpha)(tX_{j}/\alpha)^{(1-\alpha)/2}J_{\alpha-1}\big(2(tX_{j}/\alpha)^{1/2}\big)-\mathcal{H}_{X_1, \theta}(t)\Big)^2 }\, \dd P_0(t) \Big].
\end{multline*}
By inequality (\ref{besselineq2}), 
\begin{multline*}
E_{\theta} \Big( \frac{1}{n} \sum_{j=1}^n  \Gamma(\alpha)(tX_{j}/\alpha)^{(1-\alpha)/2}J_{\alpha-1}\big(2(tX_{j}/\alpha)^{1/2}\big)-\mathcal{H}_{X_1, \theta}(t)\Big)^2 \\
= n^{-1} \Var_{\theta}  \big( \Gamma(\alpha)(tX_1/\alpha)^{(1-\alpha)/2}J_{\alpha-1}\big( 2(tX_1/\alpha)^{1/2}\big)  \big)  \ \le \  n^{-1};
\end{multline*}
therefore 
\begin{align}
\label{wieand_ineq6}
P\big(| n^{-1/2}T_n - b(\theta)| & \le \epsilon b(\theta) \big) \nonumber \\
&\ge 1- \frac{1}{\epsilon^2 b^2(\theta)} \Big[ 8 \alpha^{1/2} \Big( \frac{ \Gamma(\alpha)}{\pi^{1/2} \Gamma(\alpha+\tfrac{1}{2})} \Big)^2 \ E_{\theta} (\bar{X}_n-\alpha)^2 + \frac{2}{n}  \Big].
\end{align}

Next, we write 
$$
(\bar{X}_n-\alpha)^2 = \frac{1}{n^2} \Big(\sum_{j=1}^n (X_j -\alpha) \Big)^2.
$$
and expand the sum.  Using the i.i.d. property of $X_1,\ldots,X_n$, we obtain 
\begin{align}
E_{\theta}(\bar{X}_n-\alpha)^2 &= E_{0} \Big[ \Big( \frac{1}{n} \sum_{j=1}^n (X_j -\alpha) \Big)^2 \prod_{j=1}^n (1+\theta h_{\theta}(X_j)) \Big] \nonumber\\
\label{expected_are}
&= \frac{1}{n}  E_{0} \Big[ (X_1-\alpha)^2 \prod_{j=1}^n (1+\theta h_{\theta}(X_j))\Big] \nonumber\\
& \quad\quad\quad + \frac{(n-1)}{n}  E_{0} \Big[ (X_1-\alpha) (X_2-\alpha) \prod_{j=1}^n (1+\theta h_{\theta}(X_j)) \Big].
\end{align}
Since $\int^{\infty}_0 { h_{\theta}(x)}\, \dd P_0(x)=0$, and by the assumption (\ref{assumption_efficiency}), $\int^{\infty}_0 {xh_{\theta}(x)}\,dP_0(x)=0$, for $\theta \in \Theta_1$, we conclude that $E_0 (1+ \theta h_{\theta} (X_1))=1$, and $E_{0} [(X_1-\alpha) (1+\theta h_{\theta}(X_1))]=0$.  Thus, the first term in the right-hand side of (\ref{expected_are}) equals $n^{-1} E_{0} [(X_1-\alpha)^2 (1+\theta h_{\theta}(X_1))]$, and the second term equals $0$. 

Further, by applying the Cauchy-Schwarz inequality, we obtain 
\begin{align*}
E_{0} [(X_1-\alpha)^2 (1+\theta h_{\theta}(X_1))]&=E_{0} (X_1-\alpha)^2 +\theta E_{0} \big( (X_1-\alpha)^2 h_{\theta}(X_1) \big)\\
& \le  \alpha + | \theta | \big( E_0 (X_1-\alpha)^4 \big )^{1/2} \big( E_0 (h^2_{\theta}(X_1)) \big)^{1/2}.
\end{align*}
Since $X_1 \sim Gamma(\alpha,1)$ under $H_0$ then $E_0(X_1-\alpha)^4 = 3\alpha^2 + 6\alpha$, so we have 
$$
E_{0} [(X_1-\alpha)^2 (1+\theta h_{\theta}(X_1))] \le  \alpha + | \theta | (3\alpha^2+6\alpha)^{1/2} \big( E_0 (h^2_{\theta}(X_1)) \big)^{1/2}.
$$
By the $L^2$-convergence of $h_{\theta}$ to $h$, we conclude that there exists $\theta^{*} \in (0, \eta)$ such that
$$
{\bar\sigma}^2 := \sup_{\theta \in (-\theta^{*}, \theta^{*})} E_0 [(X_1-\alpha)^2 (1+\theta h_{\theta}(X_1))] < \infty.
$$
Therefore, (\ref{wieand_ineq6}) can be written as
\begin{align*}
P\big(| n^{-1/2}T_n - b(\theta)| \le \epsilon b(\theta) \big) & \ge 1- \frac{1}{n \epsilon^2 b^2(\theta)}  \Big[ 8 \alpha^{1/2} \Big( \frac{ \Gamma(\alpha)}{\pi^{1/2} \Gamma(\alpha+\tfrac{1}{2})} \Big)^2 \, {\bar\sigma}^2 + 2 \Big],
\end{align*}
for all $\theta \in (-\theta^{*}, \theta^{*})$.  Set
$$
C := \frac{1}{\epsilon^2 \gamma} \Big[ 8 \alpha^{1/2} \Big( \frac{ \Gamma(\alpha)}{\pi^{1/2} \Gamma(\alpha+\tfrac{1}{2})} \Big)^2 \, {\bar\sigma}^2 + 2 \Big].
$$
Then for all $\theta \in (-\theta^{*}, \theta^{*})$ and $n > C/b^2(\theta)$, 
$$
P\big(| n^{-1/2}T_n - b(\theta)| \le \epsilon b(\theta) \big) \ge 1-\frac{\gamma C}{n b^2(\theta)} > 1-\gamma. 
$$
The proof now is complete. 
$\qed$

\addcontentsline{toc}{section}{References}
\bibliographystyle{ims}

\begin{thebibliography}{00}


\bibitem{allen} 
Allen, A. O. \, \textsl{Probability, Statistics, and Queueing Theory}, second edition. Academic Press, San Diego, CA, 1990.


\bibitem{bahadur3}
Bahadur, R. R. \, Rates of convergence of estimates and test statistics. \textit{The Annals of Mathematical Statistics}, {\bf 38} (1967), 303--324. 

\bibitem{bahadur2}
Bahadur, R. R. \, \textsl{Some Limit Theorems in Statistics}. SIAM, Philadelphia, PA, 1971. 

\bibitem{bahadur1}
Bahadur, R. R. \, Stochastic comparison of tests. \textit{Annals of Mathematical Statistics}, {\bf 31} (1960), 276--295.

\bibitem{henze17}
Baringhaus, L., Ebner, B., and Henze, N. \, The limit distribution of weighted $L^2$-goodness-of-fit statistics under fixed alternatives, with applications.  \textit{Annals of the Institute of Statistical Mathematics}, {\bf 69} (2017), 969--995.

\bibitem{henze92}
Baringhaus, L., and Henze, N. \, A goodness-of-fit test for the Poisson distribution based on the empirical generating function.  \textit{Statistics \& Probability Letters}, {\bf 13} (1992), 269--274.

\bibitem{henze08}
Baringhaus, L., and Henze, N. \, A new weighted integral goodness-of-fit statistic for exponentiality.  \textit{Statistics \& Probability Letters}, {\bf 78} (2008), 1006--1016.

\bibitem{ks_exp}
Baringhaus, L., and Taherizadeh, F. \, A K-S type test for exponentiality based on empirical Hankel transforms. \textit{Communications in Statistics - Theory and Methods}, {\bf 42} (2013), 3781--3792.

\bibitem{BT2010}
Baringhaus, L., and Taherizadeh, F. \, Empirical Hankel transforms and their applications to goodness-of-fit tests. \textit{Journal of Multivariate Analysis}, {\bf 101} (2010), 1445--1467.

\bibitem{barlow}
Barlow R. E., and Campo, R. \, Total time on test processes and applications to failure data analysis. In: \textsl{Reliability and Fault Tree Analysis},  pp. 451--481, SIAM, Philadelphia, PA, 1975.

\bibitem{BarLev}
Bar-Lev, S. K., Bshouty, D., Enis, P., Letac, G., Lu, I., and Richards, D. \, The diagonal multivariate natural exponential families and their classification.  \textit{Journal of Theoretical Probability}, {\bf 7} (1994), 883--929.

\bibitem{bauer}
Bauer, H. \, \textsl{Probability Theory and Elements of Measure Theory}, second English edition. Academic Press, New York, NY, 1981.

 
\bibitem{ref21}
Billingsley, P. \, \textsl{Convergence of Probability Measures}. Wiley, New York, NY, 1968. 

\bibitem{billi2}
Billingsley, P. \, \textsl{Probability and Measure}.  Wiley, New York, NY, 1979. 

\bibitem{ref6}
Brislawn, C. \, Traceable integral kernels on countably generated measure spaces. \textit{Pacific Journal of Mathematics}, {\bf 150} (1991), 229--240.

\bibitem{ref2}
Chow, Y. S., and Teicher, H. \, \textsl{Probability Theory: Independence, Interchangeability, Martingales}, second edition.  Springer, New York, NY, 1988.

\bibitem{cremerskadelka}
Cremers, H., and Kadelka, D. \, On weak convergence of stochastic processes with Lusian path spaces. \textit{Manuscripta Mathematica}, {\bf 45} (1984), 115--125.

\bibitem{czaplicki}
Czaplicki, J. M. \, \textsl{Statistics for Mining Engineering}.  CRC Press, Boca Raton, FL, 2014.

\bibitem{cuparic}
Cupari\'c, M., Milo\v{s}evi\'c, B., and Obradovi\'c M. \, New $L^2$-type exponentiality tests.  Preprint, arXiv:1809.07585v1.

\bibitem{dancisdavis}
Dancis, J., and Davis, C.  An interlacing theorem for eigenvalues of self-adjoint operators.  \textit{Linear Algebra and its Applications}, {\bf 88/89} (1987), 117--122.

\bibitem{henze18}
Ebner, B., Henze, N., and Yukich, J. E. \, Multivariate goodness-of-fit on flat and curved spaces via nearest neighbor distances. \textit{Journal of Multivariate Analysis}, {\bf 165} (2018), 231--242.

\bibitem{erdelyi1}
Erd\'elyi, A., Magnus, W., Oberhettinger, F., and Tricomi, F. G. \, \textsl{Higher Transcendental Functions}, Volume 1. McGraw-Hill, New York, NY, 1953.

\bibitem{ref5}
Erd\'elyi, A., Magnus, W., Oberhettinger, F., and Tricomi, F. G. \, \textsl{Higher Transcendental Functions}, Volume 2. McGraw-Hill, New York, NY, 1953.


\bibitem{ref20}
G{\=\i}khman, {\u I}. {\=I}., and Skorohod, A. V. \, \textsl{The Theory of Stochastic Processes}, Volume I. Springer, New York, NY, 1980.

\bibitem{GR1983}
Gupta, R. D., and Richards, D. St. P. \ Application of results of Kotz, Johnson and Boyd to the null distribution of Wilks' criterion.  In: \textsl{Contributions to Statistics: Essays in Honour of N. L. Johnson} (editor: P. K. Sen), pp. 205--210, North-Holland, Amsterdam, 1983. 

\bibitem{henze002}
G\"urtler, N., and Henze, N.  Recent and classical goodness-of-fit tests for the Poisson distribution.  \textit{Journal of Statistical Planning and Inference}, {\bf 90} (2000), 207--225.

\bibitem{gut2002}
Gut, A. \, On  the  moment  problem. \textit{Bernoulli}, {\bf 8} (2002), 407--421.

\bibitem{gut}
Gut, A. \, \textsl{Probability: A Graduate Course}, second edition. Springer, New York, NY, 2013.

\bibitem{hadjicosta19}
Hadjicosta, E. \,  \textsl{Integral Transform Methods in Goodness-of-Fit Testing}. Doctoral dissertation, Penn State University, 2018, to appear.

\bibitem{henze82}
Henze, N. \, The limit distribution for maxima of ``weighted'' $r$th-nearest-neighbour distances.  \textit{Journal of Applied Probability}, {\bf 19} (1982), 344--354.

\bibitem{henze83}
Henze, N. \, Ein asymptotischer Satz \"uber den maximalen Minimalabstand von unabh\"angigen Zufallsvektoren mit Anwendung auf einen Anpassungstest im $\bR^p$ und auf der Kugel. \textit{Metrika}, {\bf 30} (1983), 245--259.

\bibitem{henze96}
Henze, N. \, Empirical-distribution-function goodness-of-fit tests for discrete models.  \textit{Canadian Journal of Statistics}, {\bf 24} (1996), 81--93.

\bibitem{henze14}
Henze, N., Hl\'avka, Z., and Meintanis, S. G. \, Testing for spherical symmetry via the empirical characteristic function.  \textit{Statistics}, {\bf 48} (2014), 1282--1296.

\bibitem{henze02}
Henze, N., and Klar, B. \, Goodness-of-fit tests for the inverse Gaussian distribution based on the empirical Laplace transform. \textit{Annals of the Institute of Statistical Mathematics}, {\bf 54} (2002), 425--444.


\bibitem{henze00}
Henze, N., and Nikitin, Ya. Yu. \, A new approach to goodness-of-fit testing based on the integrated empirical process.  \textit{Journal of Nonparametric Statistics}, {\bf 12} (2000), 391--416.

\bibitem{hochstadt}
Hochstadt, H.  One-dimensional perturbations of compact operators.  \textit{Proceedings of the American Mathematical Society}, {\bf 37} (1973), 465--467.

\bibitem{hogg}
Hogg, R.V., and Tanis, E. A. \textit{Probability and Statistical Inference}, eighth edition. Pearson, Upper Saddle River, NJ, 2009.

\bibitem{Johnson}
Johnson, R. A., and Wichern, D.W. \, \textsl{Applied Multivariate Statistical Analysis}, fifth edition. Prentice-Hall, Upper Saddle River, NJ, 1998.

\bibitem{Karlin1964}
Karlin, S. \, The existence of eigenvalues for integral operators. \textit{Transactions of the American Mathematical Society}, {\bf 113} (1964), 1--17. 

\bibitem{Karlin}
Karlin, S. \, \textsl{Total Positivity}. Stanford University Press, Stanford, CA, 1968.

\bibitem{kass}
Kass, R. E., Eden, U. T., and Brown, E. N. \, \textsl{Analysis of Neural Data}.  Springer, New York, 2014.

\bibitem{KJB}
Kotz, S., Johnson, N. L., and Boyd, D. W. \, Series representations of distributions of quadratic forms in normal variables. I. Central case.  \textit{Annals of Mathematical Statistics}, {\bf 38} (1967), 823--837.

\bibitem{lemaitreknio}
Le Ma{\^\i}tre, O. P., and Knio, O. M.  \textsl{Spectral Methods for Uncertainty Quantification}.  Springer, New York, 2010. 

\bibitem{ref19}
Ledoux, M., and Talagrand M. \, \textsl{Probability in Banach Spaces}.  Springer, New York, 1991.

\bibitem{ref1}
Olver, F. W., Lozier, D. W., Boisvert, R. F., and Clark, C. W., editors. \, \textsl{NIST Handbook of Mathematical Functions}.  Cambridge University Press, New York, 2010.


\bibitem{postan}
Postan, M. Y., and Poizner, M. B. \, Method of assessment of insurance expediency of quay structures' damage risks in sea ports.  In: \textsl{Marine Navigation and Safety of Sea Transportation: Maritime Transport \& Shipping} (A. Weintrit and T. Neumann, editors), pp. 123--127.  CRC Press, Boca Raton, FL, 2013.


\bibitem{ref12}
Severini, T. A. \, \textsl{Elements of Distribution Theory}. Cambridge University Press, New York, 2005. 


\bibitem{ref18}
Sneddon, I. N. \, \textsl{The Use of Integral Transforms}. McGraw-Hill, New York, 1972.

\bibitem{sturgul}
Sturgul, J. R. \, \textsl{Discrete Simulation and Animation for Mining Engineers}.  CRC Press, Boca Raton, FL, 2015.

\bibitem{ref23}
Sunder, V. S. \, \textsl{Operators on Hilbert Space}. Hindustan Book Agency, New Delhi, 2015.

\bibitem{ref16}
Szeg\"o, G. \, \textsl{Orthogonal Polynomials}.  American Mathematical Society, Providence, RI, 1939.

\bibitem{ref27}
Taherizadeh, F. \, \textsl{Empirical Hankel Transform and Statistical Goodness-of-Fit Tests for Exponential Distributions}, Ph.D. Thesis, University of Hannover, 2009.


\bibitem{ref22}
Vakhania, N. N. \, \textsl{Probability Distributions on Linear Spaces}. North-Holland, New York, 1981. 

\bibitem{wieand}
Wieand, H. S. \, A condition under which the Pitman and Bahadur approaches to efficiency coincide. \textit{Annals of Statistics}, {\bf 4} (1976), 1003--1011.

\bibitem{wittingmuller}
Witting, H., and M\"uller-Funk, U. \, \textsl{Mathematische Statistik II}. B. G. Teubner, Stuttgart, 1995.

\bibitem{ref17}
Young, N. \, \textsl{An Introduction to Hilbert Space}. Cambridge University Press, New York, 1998.

\bibitem{ref10}
Zemanian, A. H. \, \textsl{Generalized Integral Transforms}.  Wiley-Interscience, New York, 1968.

\bibitem{zolotarev}
Zolotarev, V. M. \, Concerning a certain probability problem. \textit{Theory of Probability \& Its Applications}, {\bf 6} (1961), 201--204.

\end{thebibliography}

\end{document}